\documentclass{amsart}
\usepackage{graphicx}
\usepackage{subfigure,psfrag}

\theoremstyle{plain}

\numberwithin{equation}{section}
\input{tcilatex}

\begin{document}
\title[Minimizers for the p-area]{Existence and uniqueness for p-area minimizers
in the Heisenberg group}
\author{Jih-Hsin Cheng}
\address[Jih-Hsin Cheng and Jenn-Fang Hwang]{ Institute of Mathematics,
Academia Sinica, Taipei, Taiwan, R.O.C.}
\email[Jih-Hsin Cheng]{ cheng@math.sinica.edu.tw}
\urladdr{http://www.math.sinica.edu.tw}
\thanks{}
\author{Jenn-Fang Hwang}
\email[Jenn-Fang Hwang]{ majfh@math.sinica.edu.tw}
\urladdr{http://www.math.sinica.edu.tw}
\thanks{}
\author{Paul Yang}
\address[Paul Yang]{ Department of Mathematics, Princeton University,
Princeton, NJ 08544, U.S.A.}
\email{yang@Math.Princeton.EDU}
\thanks{}

\begin{abstract}
In \cite{CHMY04}, we studied $p$-mean curvature and the associated $p$-minimal
surfaces in the Heisenberg group from the viewpoint of PDE and differential geometry. In this paper,
we look into the problem through the variational formulation. We study a
generalized $p$-area and associated ($p$-) minimizers in general dimensions.

We prove the existence and investigate the uniqueness of minimizers. Since this is reduced 
to solving a degenerate elliptic equation, we need to consider the effect of
the singular set and this requires a careful study. We define the notion of
weak solution and prove that in a certain Sobolev space, a weak solution is
a minimizer and vice versa.
We also give many interesting examples in dimension 2. An intriguing point is that,
in dimension 2, a $C^2$-smooth solution from the PDE viewpoint may not be a 
minimizer. However, this statement is true for higher dimensions due to
the relative smallness of the size of the singular set.      
\end{abstract}

\maketitle

\section{Introduction and statement of the results}

The $p$-minimal (or X-minimal, H-minimal in the terminology of some
authors, e.g., \cite{FSS01}, \cite{GN96}, \cite{Pau01}) surfaces have been
studied extensively in the framework of geometric measure theory. Starting
from the work \cite{CHMY04}, we studied the subject from the viewpoint of
partial differential equations and that of differential geometry (we use the term 
$p$-minimal since this is the notion of minimal surfaces in pseudohermitian
geometry; "$p$" stands for "pseudohermitian").

Let $\Omega $ be a bounded domain in $R^{2n}.$ Let $\vec{X}$ $=$ $%
(x_{1}, $ $x_{1^{\prime }},$ $x_{2},$ $x_{2^{\prime }},$ $..,$ $x_{n},$ $%
x_{n^{\prime }})$ $\in \Omega .$ For a graph $(\vec{X},$ $u(\vec{X}))$ in
the Heisenberg group of dimension $2n+1$ with prescribed $p$-mean curvature $%
H$ $=$ $H(\vec{X}),$ the equation for $u$ $:$ $\Omega \subset R^{2n}$ $%
\rightarrow $ $R$ reads

\begin{equation}
div\frac{\nabla u-\vec{X}^{\ast }}{|\nabla u-\vec{X}^{\ast }|}=H
\label{eqn1.1}
\end{equation}

\noindent where $\vec{X}^{\ast }$ $=$ $(x_{1^{\prime }},$ $-x_{1},$ $%
x_{2^{\prime }},$ $-x_{2},...,$ $x_{n^{\prime }},$ $-x_{n})$ (see (\ref{eqn2.10}) in
Section 2 for a geometric interpretation)$.$ In general, for a vector field $\vec{G}$ $=$ $(g_{1},g_{2},...,g_{2n})$ on $%
\Omega \subset R^{2n},$ we define $\vec{G}^{\ast }$ $\equiv $ $(g_{2},$ $%
-g_{1},$ $g_{4},$ $-g_{3},$ $...,$ $g_{2n},$ $-g_{2n-1}).$
The equation (\ref{eqn1.1}) is
the Euler-Lagrange equation (away from the singular set {$\nabla u-\vec{X}^{\ast }=0$})
of the following 
energy functional (called the
$p$-area of the graph defined by $u$ if $H$ $=$ $0,$ see Section 2):

\begin{equation}
\mathcal{X}(u)=\int_{\Omega }\{|\nabla u-\vec{X}^{\ast }|+Hu\}dx_{1}\wedge
dx_{1^{\prime }}\wedge ...\wedge dx_{n}\wedge dx_{n^{\prime }}.
\label{eqn1.2}
\end{equation}

Since we consider the variation over the whole domain, the singular set will
cause the main difficulty in the study. In order to explain this, we generalize $\mathcal{X}(\cdot )$
by considering an arbitrary vector
field $\vec{F}$ $=$ $\vec{F}(\vec{X})$ instead of $-\vec{X}^{\ast }$ in the
following form:

\begin{equation}
\mathcal{F}_{q}(u)\equiv \int_{\Omega }\{|\nabla u+\vec{F}%
|^{q}+qHu\}dx_{1}\wedge dx_{2}\wedge ...\wedge dx_{m}  \label{eqn1.3}
\end{equation}

\noindent for $1\leq q<\infty$, where $\Omega$ $\subset$ $R^{m}$. Let $S(u)$ denote the singular set of $u$,
consisting of the points where $\nabla u+\vec{F}$ $=$ $0.$ Let $%
u_{\varepsilon }$ $=$ $u+\varepsilon \varphi .$ It is easy to compute (see
Section 3 for the case $q$ $=$ $1$) the first variation of $\mathcal{F}_{q}:$
(omitting the Euclidean volume element)

\begin{eqnarray}
&&\frac{d\mathcal{F}_{q}(u_{\varepsilon })}{d\varepsilon }|_{\varepsilon
=0\pm } \label{eqn1.3''}\\
&=&c_{q}\int_{S(u)}|\nabla \varphi |^{q}+\int_{\Omega \backslash
S(u)}q|\nabla u+\vec{F}|^{q-2}(\nabla u+\vec{F})\cdot \nabla \varphi
+\int_{\Omega }qH\varphi \notag
\end{eqnarray}%

\noindent where $c_{q}$ $=$ $\pm 1$ for $q$ $=$ $1$ and $c_{q}$ $=$ $0$ for $%
1<q<\infty .$ 

For $q$ $=$ $1$, can we ignore the term $\pm \int_{S(u)}|\nabla \varphi |$? A recent paper of Balogh answered this question completely. In \cite{Bal03} Balogh studied the size of the singular set $S(u)$ (called the characteristic set in \cite{Bal03}). He showed (Theorem 3.1(2) in \cite{Bal03}) that for $\vec{F}$ $=$ $-\vec{X}^{\ast }$ in $R^{2n}$, $S(u)$ has locally finite $n$-dimensional Hausdorff measure if $u$ $\in $ $C^{2}$. We obtained the same result as Lemma 5.4 in \cite{CHMY04} by a different argument (we used only elementary linear algebra and the implicit function theorem in the proof; also we were not aware of \cite{Bal03} at the time \cite{CHMY04} was written). In this paper, we generalize this result to the situation of general $\vec{F}$ (see Theorem D below and its proof in Section 6). For $u$ $\in $ $C^{1,1}$ and $\vec{F}$ $=$ $-\vec{X}^{\ast }$ in $R^{2n}$, Balogh showed (Theorem 3.1(1) in \cite{Bal03}) that $dim_{E}S(u)$ $<$ $2n-\delta$ where $dim_{E}$ denotes the Hausdorff dimension with respect to the Euclidean metric and $\delta$ depends on the Lipschitz constant of $\nabla u$. He also proved the existence of $u$ $\in$ $\cap_{0<\alpha <1}C^{1,\alpha }$ such that $S(u)$ has positive Lebesgue measure for any $\vec{F}$ $\in $ $C^{1}(\Omega )$ where $\Omega$ $\subset$ $R^{m}$ is a given bounded domain (Theorem 4.1(2) in \cite{Bal03}). In this paper, we consider functions $u$ of class $W^{1,1}$ so that the size of $S(u)$ may be large according to Balogh. Therefore for the case of $q$ $=$ $1$ in (\ref{eqn1.3''}), we can not neglect the
contribution of the singular set to define the weak solutions (see Definition
3.2) to the Euler-Lagrange equation of $\mathcal{F}_{q}$:

\begin{equation}
div\frac{\nabla u+\vec{F}}{|\nabla u+\vec{F}|^{2-q}}=H.  \label{eqn1.4}
\end{equation}

Equation (\ref{eqn1.4}) has been studied in various situations. For $%
\vec{F}$ $=$ $0,$ $H$ $=$ $0,$ (\ref{eqn1.4}) is known to be the $q$%
-harmonic equation for $1<q<\infty $ while it is the equation associated to
the least gradient problem for $q$ $=$ $1$ (see, for instance, \cite{SWZ92}, %
\cite{Juu04}, \cite{JL04}, etc.)$.$ Geometrically there is a dichotomy for
the 1-form $\Theta $ $\equiv $ $dz+F_{I}dx_{I}$ associated to the vector
field $\nabla u+\vec{F},$ where $\vec{F}$ $=$ $(F_{I}).$ Namely, the
hyperplane distribution defined by the kernel of $\Theta $ might be either
integrable or (completely) nonintegrable ($\Theta $ is called a contact form
in this case). When $\vec{F}$ $=$ $0,$ this is the integrable case. For the
nonintegrable case (e.g. $\vec{F}$ $=$ $-\vec{X}^{\ast }),$ the quantity of
the left side in (\ref{eqn1.4}) with $q$ $=$ $1$ can be realized as the $p$%
-mean curvature of the graph defined by $u$ in pseudohermitian geometry (see
(\ref{eqn2.12})).

We study equation (\ref{eqn1.4}) with $q$ $=$ $1$:

\begin{equation}
div\frac{\nabla u+\vec{F}}{|\nabla u+\vec{F}|}=H  \label{eqn1.4'}
\end{equation}

\textbf{Definition 1.1. }Let $\Omega$ be a domain in $R^{m}$, $m\geq 1$. We say $u\in C^{2}(\Omega )$ is a $C^{2}$ smooth solution to (%
\ref{eqn1.4'}) if and only if (\ref{eqn1.4'}) holds in $\Omega \backslash
S(u).$

\bigskip

In \cite{CHMY04} and \cite{CH04}, we considered $C^{2}$-smooth solutions $u$ to (\ref{eqn1.1}) (i.e., (\ref{eqn1.4'}) with $\vec{F}$ $=$ $-\vec{X}^{\ast }$) with $H$ $=$ $0$ in dimension $2$ and,
among other things, we proved a Bernstein-type theorem. Later in \cite{GP05} the authors obtained a similar Bernstein-type theorem through a
different approach. The description of the singular set for a $C^{2}$-smooth solution to (\ref{eqn1.1}) occupies a central position in \cite{CHMY04}. As a geometric application, we can show the nonexistence of $C^{2}$-smooth, closed surfaces of genus $\geq 2$ with bounded $p$-mean curvature in any pseudohermitian 3-manifold. In \cite{CHMY04} we also proved a uniqueness theorem for $C^{2}$-smooth solutions for the
Dirichlet problem of (\ref{eqn1.4'}) in $R^{2n}$. Recently Ritor\'{e} and Rosales 
proved a rigidity result for $C^{2}$-smooth surfaces of nonzero constant $p$-mean curvature and an Alexandrov-type theorem in the 3-dimensional Heisenberg group (see Theorem 6.1 and Theorem 6.10 in \cite{RR05}, respectively).

In this paper we consider $W^{1,1}$ minimizers for

\begin{equation}
\mathcal{F}(u)\equiv \int_{\Omega }\{|\nabla u+\vec{F}|+Hu\}dx_{1}\wedge
dx_{2}\wedge ...\wedge dx_{m}  \label{eqn1.3'}
\end{equation}

\noindent ((\ref{eqn1.3}) with $q$ $=$ $1).$ In Section 3 we define and show that in the space $W^{1,1}$, a minimizer for (\ref{eqn1.3'}) is a weak solution to the equation (\ref{eqn1.4'}) and
vice versa (see Theorem 3.3). In order to overcome the trouble caused by
singular sets which are not negligible, we introduce the notion of ''regular value''. Suppose $u$ $\in 
$ $W^{1,1},$ $\varphi $ $\in $ $W_{0}^{1,1}.$ Define $u_{\varepsilon }$ $%
\equiv $ $u$ $+$ $\varepsilon \varphi $ for $\varepsilon $ $\in $ $R.$ We
prove that there are at most countably many $\varepsilon $'s for which 

\begin{equation*}
\int_{S(u_{\varepsilon })}\mid \nabla \varphi \mid \neq 0
\end{equation*}

\noindent (cf. (\ref{eqn1.3''}) for $q$ $=$ $1).$ We call such an $\varepsilon $ singular,
otherwise regular. That is, the above integral vanishes for almost all
(regular) $\varepsilon $ (see Lemma 3.1)$.$ So we do not need to worry about
the size of the singular set for regular $\varepsilon $'s. The idea of
considering regular values plays a central role both in the proof of the
equivalence between minimizers and weak solutions and in the proof of the
uniqueness theorems in Section 5.

In Section 4 we prove the existence of a Lipschitz continuous minimizer for $%
\mathcal{F(\cdot )}$ with a given boundary value in the case of $H$ $=$ $0$
under the following condition on $\vec{F}$:

\begin{equation}
\partial _{K}F_{I}=\partial _{I}f_{K},\ \ I,K=1,...,m  \label{eqn1.5}
\end{equation}

\noindent for $C^{1}$-smooth functions $f_{K}$'s (cf. (\ref{eqn4.11})). We require
$\Omega $ to be a p-convex domain (see Definition 4.1).

\bigskip

\textbf{Theorem A. }\textit{Let }$\Omega $ \textit{be a p-convex bounded
domain in }$R^{m},m\geq 2$, \textit{with }$\partial \Omega \in C^{2,\alpha }$\textit{%
\ }$(0<\alpha <1)$\textit{.} \textit{Let }$\varphi \in C^{2,\alpha }(\bar{%
\Omega}).$\textit{\ Suppose }$\vec{F}$\textit{\ }$\in $\textit{\ }$%
C^{1,\alpha }(\bar{\Omega})$\textit{\ satisfies the condition (\ref{eqn1.5}) (or
(\ref{eqn4.11})) for }$C^{1,\alpha }$\textit{-smooth and bounded }$f_{K}$\textit{'s
in }$\Omega .$\textit{\ Then there exists a Lipschitz continuous minimizer }$%
u$\textit{\ }$\in $ $C^{0,1}(\bar{\Omega})$ \textit{for }$\mathcal{F}(\cdot
) $\textit{\ with }$H$ $=$ $0$ \textit{such that }$u=\varphi $\textit{\ on }$%
\partial \Omega .$

\bigskip

We note that a $C^{2}$-smooth bounded domain with positively curved
(positive principal curvatures) boundary is p-convex. Also condition
(\ref{eqn1.5}) includes the case $\vec{F}$ $=$ $-\vec{X}^{\ast }$. We can actually
find out all the solutions to (\ref{eqn1.5}) (see (\ref{eqn4.11''})). We notice that, for $n=1$%
, Pauls (\cite{Pau01}) proved the existence of a continuous $W^{1,p}$
minimizer for $\mathcal{X}(\cdot )$ under the assumption that the graph of
the prescribed boundary function $\varphi $ satisfies the bounded slope
condition (see \cite{GT83}).

The idea of the proof of Theorem A is to invoke Theorem 11.8 in \cite{GT83}
for a family of elliptic approximating equations (see also \cite{Pau01}). Namely we first solve the Dirichlet 
problem for the following equations:

\begin{eqnarray*}
Q_{\varepsilon }u &\equiv &div(\frac{\nabla u+\vec{F}}{\sqrt{\varepsilon
^{2}+|\nabla u+\vec{F}|^{2}}})=0\text{ \ \ in }\Omega ,\\
u &=&\varphi \text{ \ on }\partial \Omega  
\end{eqnarray*}

\noindent (see (\ref{eqn4.1})). We end up obtaining a
uniform $C^{1}$ bound for solutions to the above equations, and a subsequence of
solutions converges to a Lipschitz continuous minimizer as $\varepsilon$ $\rightarrow$
$0$. In
Section 4 we give the details of the proof.

In Section 5 we tackle the problem of uniqueness of minimizers in the Heisenberg
group of arbitrary dimension (see Theorem B). We also generalize the
comparison principle in \cite{CHMY04} (cf. Theorem C, Theorem C$^{\prime }$
there) to a weak version and for a wide class of $\vec{F}$'s (see
Theorem C below).

\bigskip

\textbf{Theorem B.} \textit{Let }$\Omega $ \textit{be a bounded domain in }$%
R^{2n}$\textit{. Let }$u,v$\textit{\ }$\in $\textit{\ }$W^{1,2}(\Omega )$%
\textit{\ be two minimizers for }$\mathcal{F}(\cdot )$\textit{\ such that }$%
u-v$\textit{\ }$\in $\textit{\ }$W_{0}^{1,2}(\Omega )$\textit{. Suppose }$H$ 
$\in $ $L^{\infty }(\Omega )$ and $\vec{F}$\textit{\ }$\in W^{1,2}(\Omega
)$\textit{\ satisfying }$div\vec{F}^{\ast }$\textit{\ }$>$%
\textit{\ }$0$\textit{\ (a.e.).}\textit{\ Then }$u\equiv v$\textit{\ in }$\Omega 
$\textit{\ (a.e.).}

\bigskip

We remark that in the specific case $\vec{F}$ $=$ $-\vec{X}^{\ast }$ the 
assumptions in Theorem B are satisfied. On the other hand, the condition $div \vec{F}^{\ast }>0$ is essential in Theorem B. Let $\Omega =B_{2}-{\bar B}_{1}\subset R^{2}$ 
where $B_{r}$ denotes the open ball of radius $r$. Consider the case $\vec{F}=0$ 
and $H=\frac{1}{r}$. Let $u=f(r),v=g(r)$, and $f\neq g$ with
the properties that $f(1)=g(1)$, $f(2)=g(2)$, and $f^{\prime}>0,g^{\prime}>0$
for $1\leq r \leq 2$. Then it is easy to see that $u$ and $v$ are two minimizers 
for the associated $\mathcal{F}(\cdot )$ (see also page 162 in \cite{CHMY04}).

\bigskip

\textbf{Theorem C. }\textit{Let }$\Omega $\textit{\ be a bounded domain in }$%
R^{2n}.$\textit{\ Let }$\vec{F}$\textit{\ (a vector field) }$\in
W^{1,2}(\Omega )$\textit{\ satisfy }$div\vec{F}^{\ast }$%
\textit{\ }$>$\textit{\ }$0$\textit{\ (a.e.).}\textit{ Suppose }$u,v\in
W^{1,2}(\Omega )$\textit{\ satisfy the following conditions:}

\begin{eqnarray*}
divN(u) &\geq &divN(v)\text{ in }\Omega \text{ (in the weak sense);} \\
u &\leq &v\text{ on }\partial \Omega .
\end{eqnarray*}

\noindent \textit{Then }$u\leq v$\textit{\ in }$\Omega .$

\bigskip

In Sections 6, we study the relation between $C^{2}$-smooth 
solutions and minimizers. In \cite{CHMY04} (Theorem B there), we proved that 
if $u$ is a $C^{2}$-smooth solution to (\ref{eqn1.1}) in dimension 2 with $H$ bounded near a singular point $p_{0}$, then either $p_{0}$ is isolated in $S(u)$ or there exists a small 
neighborhood $B$ of $p_{0}$ which intersects with $S(u)$ in exactly a $C^{1}$-smooth curve $\Gamma$ through $p_{0}$ (the condition on $H$ can be weaker). Moreover, $\Gamma$ divides $B$ into two disjoint nonsingular domains 
$B^{+}$ and $B^{-}$, and $N(u)(p_{0}^{+})$ $\equiv $ $\lim_{p\in B^{+}\rightarrow
p_{0}}N(u)(p)$ and $N(u)(p_{0}^{-})$ $\equiv $ $\lim_{p\in
B^{-}\rightarrow p_{0}}N(u)(p)$ exist. Also $N(u)(p_{0}^{+})$ $=$ 
$-N(u)(p_{0}^{-})$ (see Proposition 3.5 in \cite{CHMY04}). In Section 6 and
the first part of Section 7
(see Proposition $6.2$, Theorem 6.3, (\ref{eqn6.4}), and (\ref{eqn6.5})), we will generalize such a situation and give a criterion for $u$ to be a minimizer. In particular, suppose 
$u$ is $C^{2}$-smooth. Then Proposition $6.2$ or Theorem 6.3
gives a criterion for $u$ to be a minimizer in the situation $H_{m-1}(S(u))$ $>$ $0$ while if 
$H_{m-1}(S(u))$ $=$ $0$, $u$ must be a minimizer (see Lemma 6.1). 
 
In \cite{Pau01}, Pauls constructed two different $C^{2}$ (in fact $C^{\infty }$%
) smooth solutions to the $p$-minimal surface equation ((\ref{eqn1.1}) with $H$
$=$ $0$ or (\ref{eqn1.4'}) with $H$ $=$ $0$ and $\vec{F}$ $=$ $-\vec{X}%
^{\ast })$ with the same $C^{\infty }$-smooth boundary value and the same $p$%
-area in $\Omega$ $\subset$ $R^{2}$. These two solutions do not satisfy the criterion in Proposition $6.2$ or Theorem 6.3, 
hence none of them is a minimizer. We can also see this fact according to 
Theorem B (uniqueness of minimizers). In Section 7 we construct the actual
minimizer for Pauls' example (see Example 7.3). 

In dimensions higher than 2, the situation is quite different. The size of
the singular set can be relatively small under a suitable condition on $\vec{%
F}$ $=$ $(F_{I}).$ For $x\geq 0$, let $[x]$ denote the largest
integer less or equal than $x.$ In Section 6 we obtain an estimate for the size of the
singular set and a condition on $\vec{F}$ for a $C^{2}$-smooth
solution to (\ref{eqn1.4'}) to be a minimizer (see Theorems D and E below).
Recall that $dim_{E}$ denotes the Hausdorff dimension with respect to the Euclidean metric.

\bigskip

\textbf{Theorem D. }\textit{Let }$\Omega $\textit{\ be a domain in }$R^{m}.$%
\textit{\ Suppose }$u$\textit{\ }$\in $\textit{\ }$C^{2}(\Omega )$\textit{\
and }$F_{I}$\textit{\ }$\in $\textit{\ }$C^{1}(\Omega ).$\textit{\ Then for
any }$p$\textit{\ }$\in $\textit{\ }$\Omega ,$\textit{\ there exists a
neighborhood }$V$\textit{\ of }$p$\textit{\ in }$\Omega $\textit{\ such that 
}$S(u)$\textit{\ }$\cap $\textit{\ }$V$\textit{\ is a submanifold of }$V$%
\textit{ satisfying}

\begin{equation}
dim_{E}\mathit{(S(u)\cap V)\leq m-[}\frac{rank\text{ }(\partial _{J}F_{I}-\partial _{I}F_{J})(p)+1}{2}\mathit{%
].}  \label{eqn1.5'}
\end{equation}
 
\bigskip

\textbf{Theorem E. }\textit{Let }$\Omega $\textit{\ be a bounded domain in }$%
R^{m}$, $m\geq 2$. \textit{\ Suppose }$u$\textit{\ }$\in $\textit{\ }$C^{2}(\Omega )$%
\textit{\ }$\cap $\textit{\ }$C^{0}(\bar{\Omega})$\textit{\ is a }$C^{2}$%
\textit{-smooth solution to (\ref{eqn1.4'}) with }$H$\textit{\ }$\in $%
\textit{\ }$C^{0}(\Omega \backslash S(u))$\textit{\ }$\cap $\textit{\ }$%
L^{\infty }(\Omega )$\textit{ and }$F_{I}$ $\in $ $C^{1}(\Omega ).$\textit{\ Suppose there holds}

\begin{equation}
\mathit{\lbrack }\frac{rank\text{ }(\partial _{J}F_{I}-\partial _{I}F_{J})+1}{2}\mathit{]\geq 2}
\label{eqn1.6}
\end{equation}

\noindent \textit{for all }$p$ $\in $ $\Omega .$\textit{ Then }$u$%
\textit{\ is a weak solution to (\ref{eqn1.4'}) and a minimizer for (\ref{eqn1.3'}) 
if in addition }$u$ $\in $ $W^{1,1}(\Omega ).$

\bigskip


\textbf{Corollary F.} \textit{Let }$\Omega $\textit{\ be a bounded domain in }$%
R^{2n}.$\textit{\ Suppose }$u$\textit{\ }$\in $\textit{\ }$C^{2}(\Omega )$%
\textit{\ }$\cap $\textit{\ }$C^{0}(\bar{\Omega})$\textit{\ is a }$C^{2}$%
\textit{-smooth solution to the }$p$\textit{-minimal surface equation ((\ref{eqn1.1}) with} $H$ $=$ $0$%
\textit{). Then in dimension }$\geq $ $4$ $(n\geq 2),$ $u$\textit{\ is a weak solution to the }$p$\textit{-minimal surface equation and a minimizer for (\ref{eqn1.2}) with} $H$
$=$ $0$\textit{ if in addition }$u$ $\in $ $W^{1,1}(\Omega ).$

\bigskip

In Section 8 we study the uniqueness of solutions to elliptic approximating equations $Q_{\varepsilon }u$ $=$ $H$ (see (\ref{eqn4.1})), ${\varepsilon }>0$. Since this is an elliptic equation for a given ${\varepsilon }>0$, the uniqueness of solutions follows essentially from the known elliptic theory (see e.g. \cite{GT83}). But for the reader's convenience, we include a proof here as the Appendix.

We were aware of the paper \cite{Pau05} while this work was being done. After this paper was submitted, we were informed of the work \cite{RR05}. 
Some problems related to this paper were studied in \cite{Pau05} and \cite{RR05}.
We are grateful to Andrea Malchiodi for many discussions, in particular, 
in the study of Example 7.3. We would also like to thank the referee for stimulating comments and pointing out many grammatical errors.
   
  
\bigskip

\section{Hypersurfaces in the Heisenberg group}

In this section we introduce some basic notions for a hypersurface in a
pseudohermitian manifold. By viewing the Heisenberg group or $R^{2n+1}$ as a
suitable pseudohermitian manifold, we give geometric interpretations of $%
(1.1)$ and $(1.4)$.

Let $(M,J,\Theta )$ be a $(2n+1)$-dimensional pseudohermitian manifold
with an integrable $CR$ structure $J$ and a global contact form $\Theta $
such that the bilinear form $G\equiv \frac{1}{2}d\Theta (\cdot ,J\cdot )$ is
positive definite on the contact bundle $\xi \equiv \ker \Theta $ (\cite%
{Lee86}). The metric $G$ is usually called the Levi metric. Consider a
hypersurface $\Sigma $ $\subset $ $M.$ A point $p$ $\in $ $\Sigma $ is
called singular if $\xi $ coincides with $T\Sigma $ at $p.$ Otherwise, $p$
is called nonsingular and $\mathcal{V}$ $\equiv $ $\xi \cap T\Sigma $ is $%
2n-1$ dimensional in this case. There is a unique (up to sign) unit vector $%
N $ $\in $ $\xi $ that is perpendicular to $\mathcal{V}$ with respect to the
Levi metric $G.$ We call $N$ the Legendrian normal or the $p$-normal (''$p$%
'' stands for ''pseudohermitian''). Suppose that $\Sigma $ bounds a domain $%
\Omega $ in $M.$ We define the $p$-area $2n$-form $\mathcal{A}$ by computing
the first variation ($\mathcal{A}$ will be computed below for the case of the
Heisenberg group), away from the singular set, of the standard volume in
the $p$-normal $N:$

\begin{equation}
\delta _{fN}\int_{\Omega }\Theta \wedge (d\Theta )^{n}=c(n)\int_{\Sigma }f%
\mathcal{A}  \label{eqn2.1}
\end{equation}

\noindent where $f$ is a $C^{\infty }$-smooth function on $\Sigma $ with
compact support away from the singular points, and $c(n)$ $=$ $2^{n}n!$ is a
normalization constant. The sign of $N$ is determined by requiring that $%
\mathcal{A}$ is positive with respect to the induced orientation on $\Sigma
. $ So we can talk about the $p$-area of $\Sigma $ by integrating $\mathcal{A%
}$ over $\Sigma $ (which might not be closed from now on)$.$ Then we define the $p$%
-mean curvature $H$ of $\Sigma $ as the first variation of the $p$-area in
the direction of $N:$ (the support of $f$ now is also assumed to be away
from the boundary of $\Sigma $)%
\begin{equation}
\delta _{fN}\int_{\Sigma }\mathcal{A=-}\int_{\Sigma }fH\mathcal{A}.
\label{eqn2.2}
\end{equation}

Consider the Heisenberg group viewed as a (flat) pseudohermitian manifold $%
(R^{2n+1},$ $\Theta _{0},$ $J_{0}).$ Here $\Theta _{0}$ $\equiv $ $dz+$ $%
\sum_{j=1}^{n}(x_{j}dx_{j^{\prime }}-x_{j^{\prime }}dx_{j})$ at a point $(%
\vec{X},$ $z)$ $\equiv $ $(x_{1},$ $x_{1^{\prime }},$ ..., $x_{n},$ $%
x_{n^{\prime }},$ $z)$ $\in $ $R^{2n+1}$ and $J_{0}(\mathring{e}_{j})$ $%
\equiv $ $\mathring{e}_{j^{\prime }},$ $J_{0}(\mathring{e}_{j^{\prime }})$ $%
\equiv $ $-\mathring{e}_{j}$ where 
\begin{equation}
\mathring{e}_{j}\equiv \frac{\partial }{\partial x_{j}}+x_{j^{\prime }}\frac{%
\partial }{\partial z},\text{ \ }\mathring{e}_{j^{\prime }}\equiv \frac{%
\partial }{\partial x_{j^{\prime }}}-x_{j}\frac{\partial }{\partial z}
\label{eqn2.3}
\end{equation}

\noindent $j=1,2,...n$ span $\xi _{0}$ $\equiv $ $\ker \Theta _{0}.$ Let $%
\Sigma $ be a graph defined by $z=u(\vec{X}).$ Note that $\mathring{e}_{j}$
's and $\mathring{e}_{j^{\prime }}$'s form an orthonormal basis with respect
to the Levi metric $G_{0}$ $=$ $(\sum_{j=1}^{n}dx_{j}$ $\wedge $ $%
dx_{j^{\prime }})$ $(\cdot ,$ $J_{0}\cdot ).$ Observe that an element $v$ $=$
$\sum_{j=1}^{n}(a_{j}\mathring{e}_{j}$ $+$ $b_{j^{\prime }}\mathring{e}%
_{j^{\prime }})$ $\in $ $\xi _{0}\cap T\Sigma $ satisfies $d(z-u(\vec{X}))$ $%
(v)$ $=$ $0.$ It follows that 
\begin{equation}
\sum_{j=1}^{n}[(u_{x_{j}}-x_{j^{\prime }})a_{j}+(u_{x_{j^{\prime
}}}+x_{j})b_{j^{\prime }}]=0.  \label{eqn2.4}
\end{equation}

Let $N\equiv $ $-D^{-1}\sum_{j=1}^{n}[(u_{x_{j}}-x_{j^{\prime }})\mathring{%
e}_{j}+(u_{x_{j^{\prime }}}+x_{j})\mathring{e}_{j^{\prime }}]$ where $D$ $%
\equiv $ $(\sum_{j=1}^{n}[(u_{x_{j}}-x_{j^{\prime
}})^{2}\ +\ (u_{x_{j^{\prime }}}+x_{j})^{2}])^{1/2}.$ It is easy to see that $N$ is
perpendicular to $\xi _{0}\cap T\Sigma $ by (2.4), that it is of the unit length
w.r.t. $G_{0}$ and hence $N$ is the $p$-normal (that the associated $\mathcal{%
A}$ is positive will be shown below). We can now compute $\Theta _{0}\wedge
(d\Theta _{0})^{n}$ $=$ $c(n)$ $dz$ $\wedge $ $dx_{1}$ $\wedge $ $%
dx_{1^{\prime }}$ $\wedge $ $...$ $\wedge $ $dx_{n}$ $\wedge $ $%
dx_{n^{\prime }}$ and

\begin{equation}
\iota _{N}\{\Theta _{0}\wedge (d\Theta
_{0})^{n}\}=-c(n)D^{-1}\{(I)+(II)+(III)\}  \label{eqn2.5}
\end{equation}

\noindent where $\iota _{N}$ means taking the interior product with $N$ and (%
$d\hat{x}_{I}$ deleted)

\begin{eqnarray*}
(I) &=&\sum_{j=1}^{n}[(u_{x_{j}}-x_{j^{\prime }})x_{j^{\prime
}}-(u_{x_{j^{\prime }}}+x_{j})x_{j}]dx_{1}\wedge dx_{1^{\prime }}\wedge
...\wedge dx_{n}\wedge dx_{n^{\prime }} \\
(II) &=&-\sum_{j=1}^{n}(u_{x_{j}}-x_{j^{\prime }})dz\wedge dx_{1}\wedge
dx_{1^{\prime }}...d\hat{x}_{j}\wedge dx_{j^{\prime }}...\wedge dx_{n}\wedge
dx_{n^{\prime }} \\
(III) &=&\sum_{j=1}^{n}(u_{x_{j^{\prime }}}+x_{j})dz\wedge dx_{1}\wedge
dx_{1^{\prime }}...dx_{j}\wedge d\hat{x}_{j^{\prime }}...\wedge dx_{n}\wedge
dx_{n^{\prime }}.
\end{eqnarray*}

\noindent It follows that

\begin{eqnarray}
&&\delta _{fN}\int_{\Omega }\Theta _{0}\wedge (d\Theta _{0})^{n}
\label{eqn2.6} \\
&=&\int_{\Omega }L_{fN}\{\Theta _{0}\wedge (d\Theta _{0})^{n}\}=\int_{\Omega
}d(\iota _{fN}\{\Theta _{0}\wedge (d\Theta _{0})^{n}\})  \notag \\
&=&\int_{\Sigma }f\iota _{N}\{\Theta _{0}\wedge (d\Theta _{0})^{n}\}  \notag
\end{eqnarray}

\noindent by the formula $L_{v}$ $=$ $\iota _{v}\circ d$ $+$ $d\circ \iota
_{v}$ and Stokes' theorem. Substituting (2.5) into (2.6) and comparing (2.6)
with (2.1) gives

\begin{equation}
\mathcal{A=-}D^{-1}\{(I)+(II)+(III)\}  \label{eqn2.7}
\end{equation}

\noindent which simplifies to $Ddx_{1}\wedge dx_{1^{\prime }}\wedge ...\wedge
dx_{n}\wedge dx_{n^{\prime }}$ on $\Sigma $ ($z=u(\vec{X})).$ Next we compute

\begin{equation}
\delta _{fN}\int_{\Sigma }\mathcal{A=}\int_{\Sigma }L_{fN}\mathcal{A=}%
\int_{\Sigma }\iota _{fN}\circ d\mathcal{A}.  \label{eqn2.8}
\end{equation}

\noindent Here we have used Stokes' theorem and the condition that the
support of $f$ is away from the singular set and the boundary of $\Sigma .$
Noting that $D$ $=$ $|\nabla u-\vec{X}^{\ast }|$ where $\vec{X}^{\ast }$ $=$ 
$(x_{1^{\prime }},$ $-x_{1},$ $x_{2^{\prime }},$ $-x_{2},...,$ $x_{n^{\prime
}},$ $-x_{n})$, we can easily deduce that $d(D^{-1}(I))$ $=$ $0$ and

\begin{equation*}
d\{D^{-1}[(II)+(III)]\}=(div\frac{\nabla u-\vec{X}^{\ast }}{|\nabla u-\vec{X}%
^{\ast }|})dz\wedge dx_{1}\wedge dx_{1^{\prime }}\wedge ...\wedge
dx_{n}\wedge dx_{n^{\prime }}.
\end{equation*}

\noindent It follows from (2.7) and (2.5) that

\begin{equation}
\iota _{N}\circ d\mathcal{A}=-(div\frac{\nabla u-\vec{X}^{\ast }}{|\nabla u-%
\vec{X}^{\ast }|})\mathcal{A}.  \label{eqn2.9}
\end{equation}

\noindent Substituting (2.9) into (2.8) and comparing (2.8) with (2.2), we
obtain the following expression for the $p$-mean curvature $H_{\Sigma }$ of the 
graph $\Sigma $ $=$ $\{(%
\vec{X},$ $u(\vec{X}))\}:$%
\begin{equation}
H_{\Sigma }=div\frac{\nabla u-\vec{X}^{\ast }}{|\nabla u-\vec{X}^{\ast }|}.
\label{eqn2.10}
\end{equation}

Next we consider a general vector field $\vec{F}$ $=$ $(F_{I})$ instead of 
$-\vec{X}^{\ast }.$ Let $\Theta _{\vec{F}}$ $\equiv $ $dz$ $+$ $%
\sum_{I}F_{I}dx_{I}$ where $I$ ranges over $1,$ $1^{\prime },$ $...,$ $n,$ $%
n^{\prime }.$ Assume that $\Theta _{\vec{F}}$ is a contact form, i.e., $%
\Theta _{\vec{F}}$ $\wedge $ $(d\Theta _{\vec{F}})^{n}$ $\neq $ $0$
everywhere (satisfied for $\vec{F}$ $=$ $-\vec{X}^{\ast }$ as shown
previously). For instance, the condition
is equivalent to $\partial F_{1^{\prime }}/\partial x_{1}$ $-$ $\partial
F_{1}/\partial x_{1^{\prime }}$ $\neq $ $0$ in the case $n$ $=$ $1.$ Define

\begin{equation}
e_{I}=\frac{\partial }{\partial x_{I}}-F_{I}\frac{\partial }{\partial z},%
\text{ \ }I=1,1^{\prime },...,n,n^{\prime }.  \label{eqn2.11}
\end{equation}

\noindent It is easy to see that $\Theta _{\vec{F}}$ annihilates the $e_{I}$'s.
Define the $CR$ structure $J_{\vec{F}}$ on the contact bundle $\ker \Theta _{%
\vec{F}}$ by $J_{\vec{F}}(e_{j})$ $=$ $e_{j^{\prime }}$ and $J_{\vec{F}%
}(e_{j^{\prime }})$ $=$ $-e_{j}$ for $j$ $=$ $1,$ $2,$ $...,$ $n.$ For the $%
2 $-dimensional case ($n$ $=$ $1$), we can find a nonvanishing scalar
function $\lambda $ ($=$ $2(\partial F_{1^{\prime }}/\partial x_{1}$ $-$ $%
\partial F_{1}/\partial x_{1^{\prime }})^{-1}$) such that $\{e_{1},$ $%
e_{1^{\prime }}\}$ forms an orthonormal basis with respect to the Levi
metric $G_{\vec{F}} $ associated to $(J_{\vec{F}},$ $\lambda \Theta _{\vec{F}%
}).$ Let $\psi $ $\equiv $ $z$ $-$ $u(x_{1},x_{1^{\prime }})$ be a defining
function for the graph of $u.$ By a formula in Section 2 of \cite{CHMY04},
we can compute the $p$-mean curvature $H_{\vec{F}}$ with respect to the
pseudohermitian structure $(J_{\vec{F}},$ $\lambda \Theta _{\vec{F}})$ as
follows:

\begin{eqnarray}
H_{\vec{F}} &=&-div_{b}\frac{\nabla _{b}\psi }{|\nabla _{b}\psi |_{G_{\vec{F}%
}}}  \label{eqn2.12} \\
&=&-e_{1}(\frac{e_{1}\psi }{D_{\vec{F}}})-e_{1^{\prime }}(\frac{e_{1^{\prime
}}\psi }{D_{\vec{F}}})  \notag \\
&=&\frac{\partial }{\partial x_{1}}(\frac{u_{x_{1}}+F_{1}}{|\nabla u+\vec{F}|%
})+\frac{\partial }{\partial x_{1^{\prime }}}(\frac{u_{x_{1^{\prime
}}}+F_{1^{\prime }}}{|\nabla u+\vec{F}|})  \notag \\
&=&div\frac{\nabla u+\vec{F}}{|\nabla u+\vec{F}|}.  \notag
\end{eqnarray}

\noindent Here we have used $|\nabla _{b}\psi |_{G_{\vec{F}}}$ $=$ $\sqrt{%
(e_{1}\psi )^{2}+(e_{1^{\prime }}\psi )^{2}}$ $=$ $|\nabla u+\vec{F}|$ by (%
\ref{eqn2.11}).

\bigskip

\section{Minimizers in the Heisenberg group}

In this section we deduce some properties of a minimizer in the
Heisenberg group. In fact we consider a more general area functional (this
is just (\ref{eqn1.3'})):

\begin{equation}
\mathcal{F}(u)\equiv \int_{\Omega }\{\mid \nabla u+\vec{F}\mid +Hu\}
\label{eqn3.1}
\end{equation}

\noindent where $\Omega \subset R^{m}$ is a bounded domain, $\vec{F}$ is an
arbitrary (say, $L^{1})$ vector field on $\Omega ,$ and $H$ $\in $ $%
L^{\infty }(\Omega )$ (we omit the Euclidean volume element). 

\bigskip

\textbf{Definition 3.1.} $u\in W^{1,1}(\Omega )$ is called a minimizer for $%
\mathcal{F}(u)$ $\equiv $ $\int_{\Omega }\{|\nabla u+\vec{F}|$ $+$ $Hu\}$ if 
$\mathcal{F}(u)$ $\leq $ $\mathcal{F}(u+\varphi )$ for any $\varphi \in
W_{0}^{1,1}(\Omega )$, where $\vec{F}$ $\in$ $L^{1}(\Omega )$ and $H$ $\in $ $%
L^{\infty }(\Omega )$.

\bigskip

We are going
to investigate the first variation of $\mathcal{F}$. Let $u,\varphi \in
W^{1,1}(\Omega )$ and $u_{\varepsilon }\equiv u+\varepsilon \varphi $ for $%
\varepsilon \in R.$ It follows that $u_{\varepsilon }-u_{\hat{\varepsilon}%
}=(\varepsilon -\hat{\varepsilon})\varphi .$ Let $S(u_{\varepsilon })$, the
singular set of $u_{\varepsilon }$, denote the set of points where $\nabla
u_{\varepsilon }+\vec{F}=0.$ So from (\ref{eqn3.1}) (noting that 
$\mid \nabla u_{\varepsilon }+\vec{F}\mid$ $=$ $\mid \varepsilon -\hat{\varepsilon}\mid$
$\mid \nabla \varphi \mid$ on $S(u_{\hat{\varepsilon}})$) we have

\begin{eqnarray}
\mathcal{F}(u_{\varepsilon }) &=&\mid \varepsilon -\hat{\varepsilon}\mid
\int_{S(u_{\hat{\varepsilon}})}\mid \nabla \varphi \mid +\int_{\Omega
\backslash S(u_{\hat{\varepsilon}})}\mid \nabla u_{\varepsilon }+\vec{F}\mid
\label{eqn3.2} \\
&&+\int_{\Omega }Hu_{\hat{\varepsilon}}+\int_{\Omega }(\varepsilon -\hat{%
\varepsilon})H\varphi .  \notag
\end{eqnarray}

\noindent Since $\mid \nabla u_{\varepsilon }+\vec{F}\mid ^{2}-\mid \nabla
u_{\hat{\varepsilon}}+\vec{F}\mid ^{2}=2(\varepsilon -\hat{\varepsilon}%
)(\nabla u_{\hat{\varepsilon}}+\vec{F})\cdot \nabla \varphi +(\varepsilon -%
\hat{\varepsilon})^{2}\mid \nabla \varphi \mid ^{2},$ we compute from (\ref%
{eqn3.2})

\begin{eqnarray*}
\frac{\mathcal{F}(u_{\varepsilon })-\mathcal{F}(u_{\hat{\varepsilon}})}{%
\varepsilon -\hat{\varepsilon}} &=&\frac{|\varepsilon -\hat{\varepsilon}|}{%
\varepsilon -\hat{\varepsilon}}\int_{S(u_{\hat{\varepsilon}})}|\nabla
\varphi |+\int_{\Omega \backslash S(u_{\hat{\varepsilon}})}\frac{2(\nabla u_{%
\hat{\varepsilon}}+\vec{F})\cdot \nabla \varphi +(\varepsilon -\hat{%
\varepsilon})|\nabla \varphi |^{2}}{|\nabla u_{\varepsilon }+\vec{F}%
|+|\nabla u_{\hat{\varepsilon}}+\vec{F}|} \\
&&+\int_{\Omega }H\varphi .
\end{eqnarray*}

Note that the integrand of the middle term in the right-hand side of the above 
formula actually
equals $(\varepsilon -\hat{\varepsilon})^{-1}$ $(\mid \nabla u_{\varepsilon
}+\vec{F}\mid -$ $\mid \nabla u_{\hat{\varepsilon}}+\vec{F}\mid )$ whose
absolute value is less than or equal to $\mid \nabla \varphi \mid .$
Therefore by Lebesque's dominated convergence theorem, we can easily take
the limit as $\varepsilon \rightarrow \hat{\varepsilon}\pm $ ($+$: the
right-hand limit; $-$: the left-hand limit)$,$ and obtain

\begin{equation}
\frac{d\mathcal{F}(u_{\hat{\varepsilon}\pm })}{d\varepsilon }=\pm \int_{S(u_{%
\hat{\varepsilon}})}\mid \nabla \varphi \mid +\int_{\Omega \backslash S(u_{%
\hat{\varepsilon}})}N(u_{\hat{\varepsilon}})\cdot \nabla \varphi
+\int_{\Omega }H\varphi  \label{eqn3.3}
\end{equation}

\noindent where $N(v)\equiv $ $\frac{\nabla v+\vec{F}}{|\nabla v+\vec{F}|}$
is defined on $\Omega \backslash S(v).$ Note that $N(u_{\hat{\varepsilon}%
})\cdot \nabla \varphi $ $\in L^{1}(\Omega \backslash S(u_{\hat{\varepsilon}%
}))$ since $|N(u_{\hat{\varepsilon}})\cdot \nabla \varphi |$ $\leq $ $|N(u_{%
\hat{\varepsilon}})|$ $|\nabla \varphi |$ $=$ $|\nabla \varphi |$ and $%
\nabla \varphi $ $\in $ $L^{1}(\Omega )$ by the assumption. Also from the
above argument, we have the estimate

\begin{equation*}
\frac{\mid \mathcal{F}(u_{\varepsilon })-\mathcal{F}(u_{\hat{\varepsilon}})\mid }
{\mid \varepsilon -\hat{\varepsilon}\mid }\leq \int_{\Omega }\mid \nabla
\varphi \mid +||H||_{\infty }\int_{\Omega }|\varphi |.
\end{equation*}

\noindent Namely, $\mathcal{F}(u_{\varepsilon })$ is Lipschitz continuous in 
$\varepsilon $ for $\varphi \in W^{1,1}(\Omega ).$ Let $\kappa (\varepsilon
) $ denote the Lebesque measure of the set $S(u_{\varepsilon })\cap \{\nabla
\varphi \neq 0\}.$ We claim that there are at most countably many $%
\varepsilon$'s with $\kappa (\varepsilon )>0$ for a fixed $\varphi
.$ First observe that $S(u_{\varepsilon _{1}})\cap S(u_{\varepsilon _{2}})$ $%
\subset $ $\{\nabla \varphi =0\},$ and hence ($S(u_{\varepsilon _{1}})\cap
\{\nabla \varphi \neq 0\})$ $\cap $ ($S(u_{\varepsilon _{2}})\cap \{\nabla
\varphi \neq 0\})$ = $\emptyset$ (empty). Let $|\Omega |$ denote the volume of
the bounded domain $\Omega .$ So the number of $\varepsilon $ such that $%
\kappa (\varepsilon )>\frac{1}{n}$ for any positive integer is at most $%
[n|\Omega |]+1$ where $[x]$ denotes the largest integer less than or equal
to $x.$ Therefore there are at most countably many $\varepsilon$'s
with $\kappa (\varepsilon )>0.$ We call such an $\varepsilon $ singular,
otherwise regular (i.e., $\kappa (\varepsilon )=0)$. By (\ref{eqn3.3}), we
obtain (\ref{eqn3.4}) in the following Lemma.

\bigskip

\textbf{Lemma 3.1}. \textit{(1) }$\mathcal{F}(u_{\varepsilon })$\textit{\ is Lipschitz
continuous in }$\varepsilon $\textit{\ for }$\varphi \in W^{1,1}(\Omega ).$%
\textit{\ (2) There are at most countably many singular }$\varepsilon$\textit{'s.}\textit{\ (3) For a regular }$\varepsilon ,$\textit{\ }$%
\int_{S(u_{\varepsilon })}\mid \nabla \varphi \mid =0,$\textit{\ }$\frac{d%
\mathcal{F}(u_{\varepsilon })}{d\varepsilon }$\textit{\ exists, and}

\begin{equation}
\frac{d\mathcal{F}(u_{\varepsilon })}{d\varepsilon }\mathit{=}\int_{\Omega
\backslash S(u_{\varepsilon })}\mathit{N(u}_{\varepsilon }\mathit{)\cdot
\nabla \varphi +}\int_{\Omega }H\varphi \mathit{.}  \label{eqn3.4}
\end{equation}

\bigskip

Next for $\varepsilon _{2},\varepsilon _{1}$ regular with $\varepsilon
_{2}>\varepsilon _{1},$ we compute the difference of $\frac{d\mathcal{F}%
(u_{\varepsilon })}{d\varepsilon }$ for $\varepsilon =\varepsilon
_{2},\varepsilon _{1}$ by (\ref{eqn3.4}). Using $\kappa (\varepsilon _{j})$ $%
=$ $0,$ $j=1,2$ to shrink the domain of the integral, we obtain

\begin{equation}
\frac{d\mathcal{F}(u_{\varepsilon _{2}})}{d\varepsilon }-\frac{d\mathcal{F}%
(u_{\varepsilon _{1}})}{d\varepsilon }=\int_{\Omega \backslash \lbrack
S(u_{\varepsilon _{2}})\cup S(u_{\varepsilon _{1}})]}[N(u_{\varepsilon
_{2}})-N(u_{\varepsilon _{1}})]\cdot \nabla \varphi \geq 0.  \label{eqn3.5}
\end{equation}

\noindent Here we have used Lemma $5.1^{\prime}$ (also holds for $u,v\in W^{1,1}$) in %
\cite{CHMY04} to conclude the last inequality in (\ref{eqn3.5}) by noting
that $\nabla \varphi $ $=$ ($\varepsilon _{2}-\varepsilon _{1})^{-1}(\nabla
u_{\varepsilon _{2}}-\nabla u_{\varepsilon _{1}}).$ We have the following
result.

\bigskip

\textbf{Lemma 3.2. }\textit{(1) }$\frac{d\mathcal{F}(u_{\varepsilon })}{%
d\varepsilon }$\textit{\ is an increasing function of }$\varepsilon $\textit{%
\ for }$\varepsilon $\textit{\ regular. (2) Let }$\varepsilon _{j}$%
\textit{, }$j=1,2,...,$\textit{\ be a sequence of decreasing (increasing,
respectively) regular numbers tending to }$\hat{\varepsilon}$\textit{\ (}$%
\hat{\varepsilon}$\textit{\ may be singular) as }$j\rightarrow \infty $%
\textit{. Then we have}%
\begin{equation}
\lim_{j\rightarrow \infty }\frac{d\mathcal{F}(u_{\varepsilon _{j}})}{%
d\varepsilon }=\frac{d\mathcal{F}(u_{\hat{\varepsilon}+})}{d\varepsilon }%
\text{ \ (}=\frac{d\mathcal{F}(u_{\hat{\varepsilon}-})}{d\varepsilon },\text{
respectively).}  \label{eqn3.6}
\end{equation}

\bigskip

Note that we have the precise expressions for the right-hand limit $%
\frac{d\mathcal{F}(u_{\hat{\varepsilon}+})}{d\varepsilon }$ and the
left-hand limit $\frac{d\mathcal{F}(u_{\hat{\varepsilon}-})}{d\varepsilon }$
at $\hat{\varepsilon}$ in (\ref{eqn3.3}).

\bigskip

\textbf{Proof. }(1) follows from (\ref{eqn3.5}). To prove (2), first observe
that $\int_{S(u_{\varepsilon _{j}})}|\nabla \varphi |$ $=0$ by the
definition of $\varepsilon _{j}$ being regular. Therefore we have

\begin{equation}
\int_{\cup _{j=1}^{\infty }S(u_{\varepsilon _{j}})}|\nabla \varphi |=0.
\label{eqn3.7}
\end{equation}

Let $S_{\infty }\equiv \cup _{j=1}^{\infty }S(u_{\varepsilon _{j}}).$
Since $|N(u_{\varepsilon _{j}})|$ $\leq $ $1,$ we estimate $\mid
\int_{S_{\infty }}N(u_{\varepsilon _{j}})\cdot \nabla \varphi \mid \leq
\int_{S_{\infty }}|\nabla \varphi |=0$ by (\ref{eqn3.7}). So we obtain

\begin{equation}
\int_{S_{\infty }}N(u_{\varepsilon _{j}})\cdot \nabla \varphi =0.
\label{eqn3.8}
\end{equation}

\noindent It then follows from (\ref{eqn3.4}) and (\ref{eqn3.8}) that

\begin{eqnarray}
\frac{d\mathcal{F}(u_{\varepsilon _{j}})}{d\varepsilon } &=&\int_{\Omega
\backslash S(u_{\varepsilon _{j}})}N(u_{\varepsilon _{j}})\cdot \nabla
\varphi +\int_{\Omega }H\varphi  \label{eqn3.9} \\
&=&\int_{\Omega \backslash S_{\infty }}N(u_{\varepsilon _{j}})\cdot \nabla
\varphi +\int_{\Omega }H\varphi .  \notag
\end{eqnarray}

On the other hand, observe that $\lim_{j\rightarrow \infty
}N(u_{\varepsilon _{j}})=N(u_{\hat{\varepsilon}})$ in $\Omega \backslash
\lbrack S_{\infty }\cup S(u_{\hat{\varepsilon}})]$ and

\begin{equation}
N(u_{\varepsilon _{j}})=\frac{(\nabla u_{\hat{\varepsilon}}+\vec{F}%
)+(\varepsilon _{j}-\hat{\varepsilon})\nabla \varphi }{|(\nabla u_{\hat{%
\varepsilon}}+\vec{F})+(\varepsilon _{j}-\hat{\varepsilon})\nabla \varphi |}=%
\frac{(\varepsilon _{j}-\hat{\varepsilon})\nabla \varphi }{|\varepsilon _{j}-%
\hat{\varepsilon}||\nabla \varphi |}  \label{eqn3.10}
\end{equation}

in $S(u_{\hat{\varepsilon}})\backslash S_{\infty }.$ Now we compute

\begin{eqnarray}
&&\int_{\Omega \backslash S_{\infty }}N(u_{\varepsilon _{j}})\cdot \nabla
\varphi  \label{eqn3.11} \\
&=&(\int_{S(u_{\hat{\varepsilon}})\backslash S_{\infty }}+\int_{\Omega
\backslash \lbrack S_{\infty }\cup S(u_{\hat{\varepsilon}})]})N(u_{%
\varepsilon _{j}})\cdot \nabla \varphi  \notag \\
&=&\frac{\varepsilon _{j}-\hat{\varepsilon}}{|\varepsilon _{j}-\hat{%
\varepsilon}|}\int_{S(u_{\hat{\varepsilon}})\backslash S_{\infty }}|\nabla
\varphi |+\int_{\Omega \backslash \lbrack S_{\infty }\cup S(u_{\hat{%
\varepsilon}})]}N(u_{\varepsilon _{j}})\cdot \nabla \varphi  \notag \\
&\rightarrow &\pm \int_{S(u_{\hat{\varepsilon}})}|\nabla \varphi
|+\int_{\Omega \backslash S(u_{\hat{\varepsilon}})}N(u_{\hat{\varepsilon}%
})\cdot \nabla \varphi  \notag
\end{eqnarray}

\noindent as $j\rightarrow \infty $ ($+$ for decreasing $\varepsilon _{j}$; $%
-$ for increasing $\varepsilon _{j}$)$.$ Here we have used (\ref{eqn3.10})
and Lebesque's dominated convergence theorem. By (\ref{eqn3.9}), (\ref%
{eqn3.11}), and in view of (\ref{eqn3.3}), we have proved (\ref{eqn3.6}).

\begin{flushright}
Q.E.D.
\end{flushright}

\bigskip

\textbf{Definition 3.2.} Let $\Omega \subset R^{m}$ be a bounded domain. Let 
$\vec{F}$ be an $L^{1}_{loc}$ vector field on $\Omega .$ Let $H$ $%
\in $ $L^{1}_{loc}(\Omega ).$ We say $u\in W^{1}(\Omega )$ is a weak
solution to the equation (\ref{eqn1.4'}), i.e., $divN(u)$ $=$ $H$ in $\Omega $ if and only if for
any $\varphi \in C^{\infty}_{0}(\Omega ),$ there holds%
\begin{equation}
\int_{S(u)}|\nabla \varphi |+\int_{\Omega \backslash S(u)}N(u)\cdot \nabla
\varphi +\int_{\Omega }H\varphi \geq 0.  \label{eqn3.12}
\end{equation}

\bigskip

Recall that $N(u)\equiv $ $\frac{\nabla u+\vec{F}}{|\nabla u+\vec{F}|%
}$, $S(u)\equiv \{\nabla u+\vec{F}=0\},$ and $N(u)\cdot \nabla \varphi $ $%
\in L^{1}(\Omega \backslash S(u))$ since $|N(u)\cdot \nabla \varphi |$ $\leq 
$ $|N(u)|$ $|\nabla \varphi |$ $=$ $|\nabla \varphi |$ and $\nabla \varphi $ 
$\in $ $L^{1}(\Omega )$ by assumption.$.$ Note that with $\varphi $
replaced by $-\varphi $ in (\ref{eqn3.12}), we also have $%
-\int_{S(u)}|\nabla \varphi |$ $+$ $\int_{\Omega \backslash S(u)}N(u)\cdot
\nabla \varphi $ $+$ $\int_{\Omega }H\varphi $ $\leq $ $0.$ Moreover, if the
$(m-1)$-dimensional Hausdorff measure of $S(u)$ vanishes, then the equality
holds in (\ref{eqn3.12}). We remark that in 
Definition 3.2 for the case $H$ $\in $ $L^{\infty}(\Omega )$, the space 
$C^{\infty}_{0}(\Omega )$ of test functions can be replaced by 
$W^{1,1}_{0}(\Omega )$ since the former is dense in the latter in the $W^{1,1}$ 
norm (\cite{GT83}). Note that in the definition of a minimizer, we require $u$ 
$\in$ $W^{1,1}(\Omega )$,
$\vec{F}$ $\in$ $L^{1}(\Omega )$, and $H$ $\in $ $L^{\infty}(\Omega )$ while 
for the definition of a weak solution, $u$ can be in a larger space $W^{1}(\Omega )$,
$\vec{F}$ $\in$ $L^{1}_{loc}(\Omega )$, and $H$ $\in $ $L^{1}_{loc}(\Omega )$. 

\bigskip
 

\textbf{Theorem 3.3. }\textit{Let\ }$u\in W^{1,1}(\Omega ),$ $\vec{F}$ $\in$
$L^{1}(\Omega ),$\textit{ and }$H$ $\in $ $L^{\infty}(\Omega )$.\textit{\ Then }$%
u$\textit{\ is a minimizer for }$\mathcal{F}(\cdot )$\textit{\ if and only
if }$u$\textit{\ is a weak solution to the equation }$divN(u)$\textit{\ }$=$%
\textit{\ }$H$\textit{.}

\bigskip

\textbf{Proof. }Suppose $u$ is a minimizer for $\mathcal{F}(u).$ Then $\frac{%
d\mathcal{F}(u_{0+})}{d\varepsilon }$ $\geq $ $0,$ and hence (\ref{eqn3.12})
follows from (\ref{eqn3.3}) (letting $\hat{\varepsilon}=0$ in (\ref{eqn3.3}%
)). So $u$ is a weak solution. Conversely, suppose $u$ is a weak solution.
Since $\mathcal{F}(u_{\varepsilon })$ is Lipschitz continuous in $%
\varepsilon $ by Lemma 3.1 (1), $\frac{d\mathcal{F}(u_{\varepsilon })}{%
d\varepsilon }$ exists a.e. (in fact at least for regular $\varepsilon $)
and it is integrable. Moreover, we have%
\begin{equation}
\mathcal{F}(u+\varphi )-\mathcal{F}(u)=\int_{0}^{1}\frac{d\mathcal{F}%
(u_{\varepsilon })}{d\varepsilon }d\varepsilon .  \label{eqn3.13}
\end{equation}

On the other hand, from Lemma 3.2 and the definition of weak
solution (Definition 3.2), we obtain that $\frac{d\mathcal{F}(u_{\varepsilon })}{d\varepsilon 
}$ $\geq $ $0$ for any regular $\varepsilon \in \lbrack 0,1]$ in view of (%
\ref{eqn3.3}) (take $\hat{\varepsilon}=0)$. By Lemma 3.1 (2), $\frac{d%
\mathcal{F}(u_{\varepsilon })}{d\varepsilon }$ $\geq $ $0$ a.e.. It follows
from (\ref{eqn3.13}) that $\mathcal{F}(u+\varphi )$ $\geq $ $\mathcal{F}(u).$
That is to say, $u$ is a minimizer for $\mathcal{F}(u).$

\begin{flushright}
Q.E.D.
\end{flushright}




\bigskip

\section{Existence of minimizers-proof of Theorem A}

Let $\Omega $ be a bounded domain in $R^{m},m\geq 2$. Consider the
following elliptic approximation $u=u_{\varepsilon }$ ($\varepsilon $ $>$ $0)$
(a geometric interpretation can be found in \cite{Pau01}) with given
boundary value $\varphi $ ($\in C^{2,\alpha }(\bar{\Omega}),$ $0$ $<$ $%
\alpha $ $<$ $1,$ say) :

\begin{eqnarray}
Q_{\varepsilon }u &\equiv &div(\frac{\nabla u+\vec{F}}{\sqrt{\varepsilon
^{2}+|\nabla u+\vec{F}|^{2}}})=0\text{ \ \ in }\Omega ,  \label{eqn4.1} \\
u &=&\varphi \text{ \ on }\partial \Omega  \notag
\end{eqnarray}

\noindent where $\vec{F}=(F_{I})$, $I=1,...,m$. In the case of $m=2n$, $I$ ranges over $1$, $1^{\prime}$, ..., $n$, $n^{\prime}$ (e.g., $F_{I}=$ $-x_{I^{\prime }}$ for the case of a p-minimal surface. Here we use the convention that $x_{j^{{\prime}{\prime}}}$ $=$ $-x_{j}$, $j$ $=$ $1$, ..., $n$). We will make use of Theorem 11.8 in %
\cite{GT83} to solve (\ref{eqn4.1}) in $C^{2,\alpha }(\bar{\Omega})$ (then
a subsequence of $u_{\varepsilon }$ will converge to what we
want)$.$ First we check that $Q_{\varepsilon }$ is elliptic. A direct
computation shows that (summation convention applies)

\begin{eqnarray}
Q_{\varepsilon }u &=&\frac{u_{II}(\varepsilon ^{2}+|\nabla u+\vec{F}%
|^{2})-(u_{I}+F_{I})(u_{J}+F_{J})u_{IJ}}{[\varepsilon ^{2}+|\nabla u+\vec{F}%
|^{2}]^{3/2}}  \label{eqn4.2} \\
&&+\frac{(\varepsilon ^{2}+|\nabla u+\vec{F}|^{2})\partial
_{I}F_{I}-(u_{I}+F_{I})(u_{J}+F_{J})\partial _{I}F_{J}}{[\varepsilon
^{2}+|\nabla u+\vec{F}|^{2}]^{3/2}}  \notag \\
&=&a_{IJ}({\varepsilon },x,\nabla u)u_{IJ}+b({\varepsilon },x,\nabla u)  \notag
\end{eqnarray}

\noindent where 
\begin{equation}
a_{IJ}({\varepsilon },x,\nabla u)=\frac{\delta _{IJ}(\varepsilon ^{2}+|\nabla u+\vec{F}%
|^{2})-(u_{I}+F_{I})(u_{J}+F_{J})}{[\varepsilon ^{2}+|\nabla u+\vec{F}%
|^{2}]^{3/2}}  \label{eqn4.3}
\end{equation}

\noindent and

\begin{equation*}
b({\varepsilon },x,\nabla u)=\frac{(\varepsilon ^{2}+|\nabla u+\vec{F}|^{2})\partial
_{I}F_{I}-(u_{I}+F_{I})(u_{J}+F_{J})\partial _{I}F_{J}}{[\varepsilon
^{2}+|\nabla u+\vec{F}|^{2}]^{3/2}}.
\end{equation*}

For $0$ $\neq $ $(p_{I})$ $\in R^{m},$ we compute from (\ref{eqn4.3}) that

\begin{eqnarray}
a_{IJ}p_{I}p_{J} &=&\frac{(\varepsilon ^{2}+|\nabla u+\vec{F}%
|^{2})p_{I}^{2}-(u_{I}+F_{I})(u_{J}+F_{J})p_{I}p_{J}}{[\varepsilon
^{2}+|\nabla u+\vec{F}|^{2}]^{3/2}}  \label{eqn4.4} \\
&\geq &\frac{\varepsilon ^{2}p_{I}^{2}}{[\varepsilon ^{2}+|\nabla u+\vec{F}%
|^{2}]^{3/2}}>0.  \notag
\end{eqnarray}

\noindent Here we have used Cauchy's inequality \TEXTsymbol{\vert}($\nabla u+%
\vec{F})\cdot (p_{I})|^{2}$ $\leq $ $|\nabla u+\vec{F}|^{2}p_{I}^{2}$
(noting that $p_{I}^{2}$ means the sum $\Sigma _{I}p_{I}^{2}).$ It follows
from (\ref{eqn4.2}) and (\ref{eqn4.4}) that $Q_{\varepsilon }$ is elliptic.

To apply Theorem 11.8 in \cite{GT83}, we need to get an apriori estimate
in $C^{1}(\bar{\Omega})$ -norm at least. Suppose $u_{\varepsilon }$ is a $%
C^{2,\alpha }(\bar{\Omega})$ solution to the equation $Q_{\epsilon }u$ $=$ $%
0 $ (assuming $\vec{F}\in C^{1,\alpha }(\bar{\Omega});$ later replacing $%
\vec{F}$ by $\sigma \vec{F}$), $u$ $=$ $\sigma \varphi $ on $\partial \Omega
,$ $0$ $\leq $ $\sigma $ $\leq $ $1.$ In the case of $F_{I}$ $=$ $%
-x_{I^{\prime }},$ $\partial _{I}F_{I}$ $=$ $0,$ $(u_{I}+F_{I})(u_{J}+F_{J})%
\partial _{I}F_{J}$ $=$ $-(u_{J^{\prime }}+F_{J^{\prime }})(u_{J}+F_{J})$ $=$
$0,$ and hence $b(\varepsilon ,x,\nabla u)$ $=$ $0.$ Since $Q_{\varepsilon }$ is
elliptic, it follows from the maximum principle (see e.g. Problem 10.1 in \cite%
{GT83}) that

\begin{equation}
\sup_{\Omega }\text{ }|u_{\varepsilon }|\leq \sup_{\partial \Omega }\text{ }%
|u_{\varepsilon }|=\sup_{\partial \Omega }\text{ }|\sigma \varphi |\leq
\sup_{\partial \Omega }\text{ }|\varphi |.  \label{eqn4.5}
\end{equation}

\noindent Note that the right hand side is independent of $\varepsilon .$
For a general $\vec{F},$ we will invoke the comparison principle for a
second order, quasilinear operator with a ''tail'' term (namely, Theorem
10.1 in \cite{GT83}). First we can find the comparison functions as shown
below. Let $||\ \ ||_{\infty}$ denote the supremum norm. Let $B_{R}$ denote
the ball of radius $R$, centered at the origin.

\bigskip

\textbf{Lemma 4.1}. \textit{Let }$\Omega $ $\subset B_{R}$ $\subset R^{m}$\textit{ be a bounded domain.
Suppose }$\vec{F}$ $\in $ $%
C^{1}(\Omega )$\textit{ be such that }$F_{I}$\textit{ and }$\partial _{I}F_{J}$\textit{ are all bounded
in }$\Omega .$\textit{ Then there are }$C^{\infty }$\textit{-smooth functions} 

$$w=e^{x_{1}+\kappa R}+ e^{x_{2}+\kappa R},\ w^{\prime}=-e^{x_{1}+{\kappa}^{\prime}R}- e^{x_{2}+{\kappa}^{\prime}R}
$$ 

\noindent \textit{in }$R^{m}$\textit{, where }$\kappa =\ \kappa (\varepsilon ,\ R,\ ||F_{I}||_{\infty},$ $||\partial _{I}F_{J}||_{\infty})$ $>$ $0$\textit{ and }${\kappa}^{\prime} =\ {\kappa}^{\prime} (\varepsilon ,\ R,\ ||F_{I}||_{\infty},$
$||\partial _{I}F_{J}||_{\infty})$ $>$ $0$\textit{, such that }$Q_{\varepsilon }w$ $>$ $0$\textit{ and }$%
Q_{\varepsilon }w^{\prime }$ $<$ $0$\textit{ in }$\Omega .$\textit{ Moreover, we can choose }$\kappa$\textit{ and }${\kappa}^{\prime}$\textit{ independent of }$\varepsilon$\textit{ (but depending on }$\varepsilon_{0}$\textit{) for }$0<\varepsilon \leq {\varepsilon}_{0}$\textit{, a positive constant.} 

\bigskip

\textbf{Proof}. Let $w$ have the above expression with $\kappa $ to be determined later. Let $%
w_{1}$ $\equiv $ $\partial _{x_{1}}w,$ $w_{11}$ $\equiv $ $\partial
_{x_{1}}^{2}w,$ $w_{12}$ $\equiv $ $\partial _{x_{2}}\partial _{x_{1}}w,$ and so on. It follows that

\begin{eqnarray}
w_{11} &=&w_{1}=e^{x_{1}+\kappa R},\text{ }w_{22}=w_{2}=e^{x_{2}+\kappa R},\text{ and}  \label{eqn4.6}
\\
w_{IJ} &=&0,\text{ otherwise.}  \notag
\end{eqnarray}

\noindent In view of (\ref{eqn4.2}) with $u$ replaced by $w,$ we compute the
dominating (will be clear soon) term in the numerator, which is cubic in $w$
as follows:

\begin{eqnarray}
&&(\sum_{J}w_{J}^{2})(\sum_{I}w_{II})-\sum_{I,J}w_{I}w_{J}w_{IJ}
\label{eqn4.7} \\
&=&w_{1}^{2}w_{22}+w_{2}^{2}w_{11}\text{ \ (by
(\ref{eqn4.6}))}  \notag \\
&=&e^{2x_{1}+x_{2}+3\kappa R}+e^{2x_{2}+x_{1}+3\kappa R}%
\text{ \ (by (\ref{eqn4.6})).}  \notag
\end{eqnarray}

\noindent It is easy to see that any other term in the expansion of the
numerator is bounded by either $c_{1}e^{2\kappa R},$ $c_{2}e^{\kappa R}$ or $%
c_{3}$ for $\kappa $ large$.$ Here $c_{i}\ =\ c_{i}(\varepsilon ,\ R,\ ||F_{I}||_{\infty},\ ||\partial _{I}F_{J}||_{\infty})$, $i=1,\ 2,\ 3$, are independent of $\kappa .$ Therefore we have $Q_{\varepsilon }w$ $>$ $0$ in $\Omega $ by (\ref{eqn4.7}) for a large $\kappa$ $=$ $\kappa (\varepsilon ,\ R,\ ||F_{I}||_{\infty},\ ||\partial _{I}F_{J}||_{\infty})$. Moreover, $\kappa $ is independent
of $\varepsilon$ for $0<\varepsilon \leq {\varepsilon}_{0}$, a positive constant. Similarly, we can find ${\kappa}^{\prime} =\ {\kappa}^{\prime} (\varepsilon ,\ R,\ ||F_{I}||_{\infty},\ ||\partial _{I}F_{J}||_{\infty})\ >0$ such that $Q_{\varepsilon }w^{\prime }$ $<$ $0$ in $\Omega $.

\begin{flushright}
Q.E.D.
\end{flushright}

\bigskip

\textbf{Proposition 4.2. }\textit{Let }$\Omega \subset B_{R} \subset R^{m}$\textit{ be a bounded
domain. Let }$\vec{F}$ $\in $ $C^{1}(\Omega )$\textit{ such that }$F_{I}$\textit{ and }$%
\partial _{I}F_{J}$\textit{ are all bounded in }$\Omega .$\textit{ Suppose }$u_{\varepsilon }$ 
$\in $ $C^{2}(\Omega )\cap C^{0}(\bar{\Omega})$\textit{ satisfies (\ref{eqn4.1}),
i.e., }$Q_{\varepsilon }u_{\varepsilon }$ $=$ $0$\textit{ in }$\Omega $\textit{ and }$%
u_{\varepsilon }$ $=$ $\sigma \varphi $ $\in $ $C^{0},$ $0$ $\leq $ $\sigma $
$\leq $ $1,$\textit{ on }$\partial \Omega $.\textit{ Then there exists a constant }$C$ $=$ $ C(\varepsilon ,\ R,$ $||F_{I}||_{\infty},$ $||\partial _{I}F_{J}||_{\infty},$ $||\varphi ||_{\infty})$\textit{ (independent of }$\sigma$\textit{) such that}

\begin{equation}
\sup_{\Omega }|u_{\varepsilon }|\leq C.  \label{eqn4.8}
\end{equation}

\noindent \textit{Moreover, }\textit{the bounds hold uniformly}\textit{ for }$0<\varepsilon \leq {\varepsilon}_{0}$\textit{, a positive constant.}
 
\bigskip

\textbf{Proof.} Let $w$, $w^{\prime }$ be the comparison functions as in Lemma 4.1. On $\partial \Omega ,$ $w$ $\leq $ $\sigma \varphi +C_{1}$ $%
= $ $u_{\varepsilon }+C_{1}$ for some constant $C_{1}$ $=$ $C_{1}(\varepsilon ,\ R,$  $||F_{I}||_{\infty},$ $||\partial _{I}F_{J}||_{\infty},$ $||\varphi ||_{\infty})$ (independent of $\sigma$) independent of $\varepsilon$ for $0<\varepsilon \leq {\varepsilon}_{0}$, a positive constant. On the other hand, we have $Q_{\varepsilon }w$ $>$ $0=Q_{\varepsilon }(u_{\varepsilon }+C_{1})$ in $\Omega $ by Lemma
4.1 and the observation that $Q_{\varepsilon }(u_{\varepsilon }+C_{1})$ $=$ $%
Q_{\varepsilon }u_{\varepsilon }$. Now we apply the comparison principle for
quasilinear operators (e.g. Theorem 10.1 in \cite{GT83}) to conclude that

\begin{equation}
w\leq u_{\varepsilon }+C_{1}\text{ in }\Omega .  \label{eqn4.9}
\end{equation}

Similarly, there is a constant $C_{2}$ $=$ $C_{2}(\varepsilon ,\ R,$ $||F_{I}||_{\infty},$  $||\partial _{I}F_{J}||_{\infty},$ $||\varphi ||_{\infty})$ (independent of $\sigma$) independent of $\varepsilon$ for $0<\varepsilon \leq {\varepsilon}_{0}$, such that $w^{\prime }$ $\geq $ $\sigma \varphi -C_{2}$ $=$ $%
u_{\varepsilon }-C_{2}$ on $\partial \Omega $ and $Q_{\varepsilon }w^{\prime
}$ $<$ $0$ $=$ $Q_{\varepsilon }(u_{\varepsilon }-C_{2})$ in $\Omega .$
So we obtain from the comparison principle that

\begin{equation}
w^{\prime }\geq u_{\varepsilon }-C_{2}\text{ in }\Omega .  \label{eqn4.10}
\end{equation}

\noindent Thus (\ref{eqn4.8}) follows from (\ref{eqn4.9}) and (\ref{eqn4.10}%
).

\begin{flushright}
Q.E.D.
\end{flushright}

\bigskip

For the gradient estimate, we will reduce the problem to a gradient estimate at the boundary. 
We need to require a condition on $\vec{F}.$ Suppose
there are $C^{1}$-smooth functions $f_{K}$'s ($K$ $=$ $1,$ $...,$
$m)$ in $\Omega $ such that

\begin{equation}
\partial _{K}F_{I}=\partial _{I}f_{K}.  \label{eqn4.11}
\end{equation}

We remark that if both $F_{I}$ and $G_{I}$ satisfy the condition (\ref%
{eqn4.11}), so does $F_{I}$ $+$ $G_{I}$. In fact, we can write down all the
(local) solutions to (\ref{eqn4.11}). It is easy to see from (\ref{eqn4.11})
that $\partial _{K}$ $(\partial _{J}F_{I}$ $-$ $\partial _{I}F_{J})$ $=$ $0$
for all $I,$ $J,$ $K$ $=$ $1,$ $...,$ $m.$ It
follows that

\begin{equation}
\partial _{J}F_{I}-\partial _{I}F_{J}=C_{IJ}  \label{eqn4.11'}
\end{equation}

\noindent where the constants $C_{IJ}$ satisfy the skew- symmetric relation: 
$C_{IJ}=-C_{JI}.$ Since the left-hand side of (\ref{eqn4.11'}) is linear in $\vec{%
F},$ the general solutions are the solutions to $\partial _{J}F_{I}$ $-$ $%
\partial _{I}F_{J}$ $=$ $0$ plus a special solution. Let $\omega $ $\equiv $ 
$\sum_{I}F_{I}dx_{I}.$ Then $d\omega $ $=$ $0$ if $\partial _{J}F_{I}$ $-$ $%
\partial _{I}F_{J}$ $=$ $0.$ So locally there is a function $g$ such that $%
\omega $ $\equiv $ $dg.$ Hence $F_{I}$ $=$ $\partial _{I}g.$ On the other
hand, we observe that $\tilde{F}_{I}$ $\equiv $ $\frac{1}{2}%
\sum_{K}C_{IK}x_{K}$ is a special solution to (\ref{eqn4.11'}). So the
general solutions to (\ref{eqn4.11'}) are

\begin{equation}
F_{I}=\partial _{I}g+\frac{1}{2}\sum_{K}C_{IK}x_{K}.  \label{eqn4.11''}
\end{equation}

\noindent It is then easy to verify that $\vec{F}$ $=$ $(F_{I})$ having the
form (\ref{eqn4.11''}) are also solutions to (\ref{eqn4.11}) for $f_{K}$ $%
=$ $\partial _{K}g$ $+$ $\frac{1}{2}\sum_{J}C_{JK}x_{J}.$

\bigskip

\textbf{Proposition 4.3.}\textit{ Let }$\Omega \subset R^{m}$\textit{ be a bounded
domain. Let }$\vec{F}$ $=$ $(F_{I})$ $\in $ $C^{1}(\Omega )$\textit{ satisfy the
condition (\ref{eqn4.11}), for all }$I,K$ $=$ $1,$ $...,$ $%
m$\textit{, where all }$f_{K}$\textit{'s are bounded. Suppose }$u_{\varepsilon }$ $%
\in $ $C^{2}(\bar{\Omega})$\textit{ satisfies the equation }$Q_{\varepsilon
}u_{\varepsilon }$ $=$ $H_{0}$\textit{, a constant,}\textit{ in }$\Omega $.\textit{ Then we have}

\begin{equation}
\sup_{\Omega }|\partial _{K}u_{\varepsilon }|\leq \sup_{\partial \Omega
}|\partial _{K}u_{\varepsilon }|+2||f_{K}||_{\infty}.  \label{eqn4.12}
\end{equation}

\bigskip

\textbf{Proof.} Write $\nabla u+\vec{F}=(u_{I}+F_{I}).$ Let $D_{\varepsilon
}(u)\equiv \sqrt{\varepsilon ^{2}+|\nabla u+\vec{F}|^{2}}.$ Compute (summing
over $J$ while fixing $I$ and $K$)

\begin{eqnarray}
&&\partial _{K}\frac{u_{I}+F_{I}}{D_{\varepsilon }(u)}  \label{eqn4.13} \\
&=&\frac{u_{IK}+\partial _{K}F_{I}}{D_{\varepsilon }(u)}-\frac{%
(u_{I}+F_{I})(u_{J}+F_{J})(u_{JK}+\partial _{K}F_{J})}{D_{\varepsilon
}^{3}(u)}  \notag \\
&=&\frac{\delta _{IJ}-\nu _{I}(u)\nu _{J}(u)}{D_{\varepsilon }(u)}\partial
_{J}(u_{K}+f_{K})  \notag
\end{eqnarray}

\noindent where $\nu _{I}(u)\equiv (u_{I}+F_{I})/D_{\varepsilon }(u)$ and we
have used the condition (\ref{eqn4.11}). Now for $v\in C_{0}^{2}(\Omega ),$
we compute

\begin{eqnarray}
0&=&\int_{\Omega }(Q_{\varepsilon }u_{\varepsilon }-H_{0})\partial
_{K}v=\int_{\Omega }\partial _{I}\frac{(u_{\varepsilon })_{I}+F_{I}}{%
D_{\varepsilon }(u_{\varepsilon })}\partial _{K}v\text{ (summing over }I%
\text{)}  \label{eqn4.14} \\
&=&-\int_{\Omega }\frac{(u_{\varepsilon })_{I}+F_{I}}{D_{\varepsilon
}(u_{\varepsilon })}\partial _{I}\partial _{K}v  \notag \\
&=&\int_{\Omega }\partial _{K}\frac{(u_{\varepsilon })_{I}+F_{I}}{%
D_{\varepsilon }(u_{\varepsilon })}\partial _{I}v\text{ \ (}\partial
_{I}\partial _{K}=\partial _{K}\partial _{I}\text{)}  \notag \\
&=&\int_{\Omega }\{a_{IJ}(\varepsilon ,x,\nabla u_{\varepsilon })\partial
_{J}[(u_{\varepsilon })_{K}+f_{K}]\}\partial _{I}v\text{ (summing over }I%
\text{ and }J\text{)}  \notag
\end{eqnarray}

\noindent by (\ref{eqn4.13}) with $u$ replaced by $u_{\varepsilon }.$ Here $%
a_{IJ}(\varepsilon ,x,\nabla u_{\varepsilon })$ $=$ $[\delta _{IJ}-\nu
_{I}(u_{\varepsilon })\nu _{J}(u_{\varepsilon })]/D_{\varepsilon
}(u_{\varepsilon })$ (cf. (\ref{eqn4.3})). It is then easy to see that (\ref%
{eqn4.14}) holds also for $v\in C_{0}^{1}(\Omega )$ (use the regularization $%
v_{h}$ of (7.13) in \cite{GT83} to approximate $v$). So $(u_{\varepsilon
})_{K}+f_{K}$ is a weak solution to the equation $Lw$ $\equiv $ $\partial
_{I}\{a_{IJ}(\varepsilon ,x,\nabla u_{\varepsilon })\partial _{J}w\}$ $=$ $0$ (cf. (8.2)
in \cite{GT83})$.$ By (\ref{eqn4.4}), this is an elliptic equation in
divergence form$.$ So by the maximum principle (e.g. Theorem 8.1 in \cite%
{GT83} with $b^{i}=$ $c^{i}=$ $d=$ $0$ and $a_{IJ}$ bounded), we have%
\begin{equation*}
\sup_{\Omega }|(u_{\varepsilon })_{K}+f_{K}|\leq \sup_{\partial \Omega
}|(u_{\varepsilon })_{K}+f_{K}|.
\end{equation*}

\noindent Then (\ref{eqn4.12}) follows.

\begin{flushright}
Q.E.D.
\end{flushright}

\bigskip

For a general $\vec{F},$ the bound for $\nabla u_{\varepsilon }$ may
depend on $\varepsilon $ if we invoke the maximum principle for a more
general situation (for instance, Theorem 8.16 in \cite{GT83}).

To perform the boundary gradient estimate, we need a comparison function to
apply the comparison principle. Let $\Omega $ $\subset $ $R^{m}$ be a
bounded domain with coordinates denoted by $x_{1},$ $x_{2},$ $%
...,$ $x_{m}.$ We call a coordinate system orthonormal if
it is obtained by a translation and a rotation from $x_{1},$ $x_{2},$ $...,$ $x_{m}.$ We define a certain notion of
convexity for $\Omega $ as follows.

\bigskip

\textbf{Definition 4.1.} We call $\Omega $ $\subset $ $R^{m}$ parabolically convex or p-convex
in short if for any $p\in \partial \Omega ,$ there exists an orthonormal
coordinate system $(\tilde{x}_{1},$ $\tilde{x}_{2},$ $...,$ $\tilde{x}_{m})$ with the origin at $p$ and $%
\Omega \subset \{a\tilde{x}_{1}^{2}-\tilde{x}_{2}<0\}$ where $a>0$
is independent of $p.$

\bigskip

Note that a $C^{2}$-smooth bounded domain with positively curved
(positive principal curvatures) boundary is p-convex. 

\bigskip

\textbf{Proposition 4.4}.\textit{ Let }$\Omega $ $\subset $ $R^{m}$\textit{ be a p-convex
bounded domain. Suppose }$u_{\varepsilon }$ $\in $ $C^{2}(\Omega )\cap C^{1}(%
\bar{\Omega})$\textit{ satisfies }$Q_{\varepsilon }u_{\varepsilon }$ $=$ $0$\textit{ in }$%
\Omega $\textit{ and }$u_{\varepsilon }$ $=$ $\sigma \varphi $ $\in $ $C^{2}(\bar{%
\Omega})$\textit{ on }$\partial \Omega $\textit{ with }$\vec{F}$ $\in $ $C^{1}(\bar{\Omega})$\textit{
for }$0$ $\leq $ $\sigma $ $\leq $ $1$.\textit{ Then there exists a constant }$C$ $=$ 
$C(\varepsilon ,\ a,$ $||F_{I}||_{\infty},$ $||\partial _{I}F_{J}||_{\infty},$ $||\partial _{I}\varphi ||_{\infty},$ $||\partial _{I}\partial _{J}\varphi ||_{\infty})$\textit{ (independent of }$\sigma$\textit{)}\textit{
such that}

\begin{equation}
\sup_{\partial \Omega }|\nabla u_{\varepsilon }|\leq C.  \label{eqn4.15}
\end{equation}

\noindent \textit{Moreover, }\textit{the bounds hold uniformly}\textit{ for }$0<\varepsilon \leq {\varepsilon}_{0}$\textit{, a positive constant.}
 
\bigskip

\textbf{Proof}. Given $p$ $\in $ $\partial \Omega ,$ we have an orthonormal
coordinate system $(\tilde{x}_{1},$ $\tilde{x}_{2},$ $...,$ $\tilde{x}_{m})$ as in the definition of
p-convexity. Consider the comparison function $w$ $=$ $\alpha G$ $+$ $\sigma
\varphi $ where $G$ is the function $\tilde{G}\equiv $ $a\tilde{x}_{1}^{2}-%
\tilde{x}_{2}$ viewed as a function of $(x_{I}),$ $I=1,\ 2,\ ...,\ m,$ for large $\alpha $ to be determined. In view of the
invariance of $Q_{\varepsilon }(u)$ under the coordinate changes of
translations and rotations, we compute ($\tilde{Q}_{\varepsilon },$ $\tilde{D%
}_{\varepsilon }$ being the corresponding operator, quantity of $%
Q_{\varepsilon },$ $D_{\varepsilon }$ with respect to $(\tilde{x}_{I}),$
respectively)

\begin{eqnarray}
Q_{\varepsilon }(w) &=&\tilde{Q}_{\varepsilon }(\tilde{w})\text{ \ (}\tilde{w%
}\text{ is }w\text{ viewed as a function of }(\tilde{x}_{I}))
\label{eqn4.16} \\
&=&\frac{P(\tilde{G})\alpha ^{3}+A\alpha ^{2}+B\alpha +E}{\tilde{D}%
_{\varepsilon }(\tilde{w})^{3}}  \notag
\end{eqnarray}

\noindent by (\ref{eqn4.2}) where $P(\tilde{G})$ is the corresponding
quantity of $P(G)$ $\equiv $ $G_{x_{1}}^{2}G_{x_{2}x_{2}}$ $-$ $2G_{x_{1}}G_{x_{2}}$$\ $$G_{x_{1}x_{2}}$ $+$ $%
G_{x_{2}}^{2}G_{x_{1}x_{1}}$ with respect to $(\tilde{x}_{I}),$
and $A$ is a function of $a,\ F_{I},\ \partial _{I}F_{J},\ \partial _{I}\varphi ,\ \partial _{I}\partial _{J}\varphi $
while $B,$ $E$ are functions of $\varepsilon ,\ a,\ F_{I},\ \partial _{I}F_{J},\ \partial _{I}\varphi ,\ \partial _{I}\partial _{J}\varphi .$ Moreover, a direct computation shows that $%
P(\tilde{G})$ $=$ $2a.$ Since $a>0,$ $Q_{\varepsilon }(w)$ $\geq $ ($\leq $,
respectively) $0$ $=$ $Q_{\varepsilon }(u_{\varepsilon })$ for positive
(negative, respectively) large $\alpha$ $=$ $\alpha (\varepsilon ,\ a,$  $||F_{I}||_{\infty},$ $||\partial _{I}F_{J}||_{\infty},$ $||\partial _{I}\varphi ||_{\infty},$ $||\partial _{I}\partial _{J}\varphi ||_{\infty})$ by (\ref{eqn4.16}). Note that $\alpha $ is independent of $%
\sigma $ and independent of $\varepsilon $ for $0<\varepsilon \leq {\varepsilon}_{0}.$ On the other hand, $w$ $=$ $\alpha G$ $+$ $%
\sigma \varphi $ $\leq $ ($\geq $, respectively) $\sigma \varphi $ $=$ $%
u_{\varepsilon }$ on $\partial \Omega $ since $G$ $\leq $ $0$ on $\bar{\Omega%
}$ by the p-convexity. Therefore $w$ $\leq $ ($\geq $, respectively)$%
u_{\varepsilon }$ in $\Omega $ by the comparison principle for second order
quasilinear operators (see e.g. Theorem 10.1 in \cite{GT83}). Noting that $%
G(p)=0 $ and hence $w(p)$ $=$ $\sigma \varphi (p)$ $=$ $u_{\varepsilon }(p)$%
, we then have

\begin{equation}
\frac{\partial u_{\varepsilon }}{\partial \nu }\leq (\geq ,\text{
respectively})\frac{\partial w}{\partial \nu }  \label{eqn4.17}
\end{equation}

\noindent where $\nu $ $=$ $-{\partial}_{{\tilde{x}}_{2}}$ at $p.$ Observe that $%
\frac{\partial w}{\partial \nu }$ (for either positive or negative $\alpha )$
is bounded by a constant depending on $\varepsilon ,\ a,$ $||F_{I}||_{\infty},$ $||\partial _{I}F_{J}||_{\infty},$ $||\partial _{I}\varphi ||_{\infty},$ $||\partial _{I}\partial _{J}\varphi ||_{\infty}$, but independent of $\sigma $, $p$ (moreover, the bounds hold for $0<\varepsilon \leq {\varepsilon}_{0}$, a positive constant),   
so is $\frac{\partial u_{\varepsilon }}{\partial \nu }$ by (\ref{eqn4.17}).
Since $u_{\varepsilon }$ $=$ $\sigma \varphi $ on $\partial \Omega $ and $u_{\varepsilon }$ $-$ 
$\sigma\varphi$ $\in$ $C^{1}(\bar{\Omega})$, we can easily show that the derivatives of $u_{\varepsilon }$ $-$ $\sigma\varphi$ in the $\tilde{x}_{1}$, $\tilde{x}_{3}$, $...,$ $\tilde{x}_{m})$ (except $\tilde{x}_{2}$) directions all vanish at $p$. It follows that in the $\tilde{x}_{j}$ $(j\neq 2)$ direction, the derivative of $u_{\varepsilon }$ is the same as the derivative of $\sigma \varphi $. So of
course it is bounded by $||\nabla \varphi ||_{\infty}$ (note that $0$ $\leq $ $\sigma $ $\leq $ $1$)$.$ Altogether we have proved
(\ref{eqn4.15}).

\begin{flushright}
Q.E.D.
\end{flushright}

\bigskip

\textbf{Proof of Theorem A.}

In order to apply Theorem 11.8 in \cite{GT83} to solve the Dirichlet
problem (\ref{eqn4.1}), we consider a family of equations:

\begin{eqnarray}
Q_{\varepsilon ,\sigma }u &\equiv &div\frac{\nabla u+\sigma \vec{F}}{\sqrt{%
\varepsilon ^{2}+|\nabla u+\sigma \vec{F}|^{2}}}=0\text{ in }\Omega
\label{eqn4.18} \\
u &=&\sigma \varphi \text{ \ on }\partial \Omega ,\text{ \ }0\leq \sigma
\leq 1.  \notag
\end{eqnarray}

\noindent Express $Q_{\varepsilon ,\sigma }u$ $=$ $a_{IJ}(\varepsilon ,x,\nabla u;\sigma
)u_{IJ}$ $+$ $b(\varepsilon ,x,\nabla u;\sigma )$ where $a_{IJ}(\varepsilon ,x,\nabla u;\sigma )$ and $%
b(\varepsilon ,x,\nabla u;\sigma )$ are given by (\ref{eqn4.3}) with $\vec{F}$ replaced
by $\sigma \vec{F}.$ It is then easy to check that the conditions (i), (ii),
(iii) on page 287 of \cite{GT83} are satisfied. To have an apriori H\"{o}%
lder estimate for $\nabla u,$ we invoke Theorem 13.2 in \cite{GT83}.
Comparing (\ref{eqn4.18}) with (13.2) in \cite{GT83} gives

\begin{equation*}
\mathbf{A(}x,u,\nabla u)=\frac{\nabla u+\sigma \vec{F}}{\sqrt{\varepsilon
^{2}+|\nabla u+\sigma \vec{F}|^{2}}},B(x,u,\nabla u)=0.
\end{equation*}

\noindent Following pages 319-320 of \cite{GT83}, we find $\bar{a}^{IJ}$ $%
\equiv $ $D_{p_{J}}A^{I}$ $=$ $a_{IJ}(\varepsilon ,x,\nabla u;\sigma )$ and $\lambda
(\varepsilon ,x,u,$ $\nabla u)$ $=$ $\varepsilon ^{2}/[\varepsilon ^{2}+|\nabla u+\sigma 
\vec{F}|^{2}]^{3/2}$ by (\ref{eqn4.4})$.$ Therefore we can take $\lambda
_{K} $ $=$ $\varepsilon ^{2}/[\varepsilon ^{2}+(K+C)^{2}]^{3/2}$ in (13.4)
of \cite{GT83}, in which $K$ $\equiv $ $|u|_{1;\Omega }$ (see page 53 in \cite{GT83} for the notation) and $C$ $\equiv $ $||\vec{F}||_{\infty}$. Similarly we estimate

\begin{equation*}
|D_{p_{J}}A^{I}|=|a_{IJ}(\varepsilon ,x,\nabla u;\sigma )|\leq \frac{1}{\sqrt{\varepsilon
^{2}+|\nabla u+\sigma \vec{F}|^{2}}}\leq \frac{1}{\varepsilon }.
\end{equation*}%
\noindent So we can take $\Lambda _{K}$ $=$ $\varepsilon ^{-1}.$ Since both $%
D_{z}A^{I}$ and $B$ vanish, we compute

\begin{eqnarray*}
|\delta _{J}A^{I}|+|B| &=&|D_{x_{J}}A^{I}| \\
&=&\frac{|(\varepsilon ^{2}+|\nabla u+\sigma \vec{F}|^{2})\partial
_{J}(\sigma F_{I})-(u_{I}+\sigma F_{I})(u_{L}+\sigma F_{L})\partial
_{J}(\sigma F_{L})|}{[\varepsilon ^{2}+|\nabla u+\sigma \vec{F}|^{2}]^{3/2}}
\\
&\leq &\frac{(\frac{3}{2}+n)\sup_{K,I}|\partial _{K}F_{I}|}{\sqrt{%
\varepsilon ^{2}+|\nabla u+\sigma \vec{F}|^{2}}}.
\end{eqnarray*}%
\noindent Therefore we can take an upper bound $\mu _{K}$ $=$ $\varepsilon
^{-1}(\frac{3}{2}+n)\sup_{J,I}|\partial _{J}F_{I}|.$ Now by Theorem 13.2 in %
\cite{GT83}, we have an apriori H\"{o}lder bound for $\nabla u$ in terms of $%
n,$ $K$ $(\equiv $ $|u|_{1;\Omega }),$ $\Lambda _{K}/\lambda _{K},$ $\mu
_{K}/\lambda _{K},$ size of $\Omega ,$ and $|\varphi |_{2;\Omega }.$ On the
other hand, we observe that Lemma 4.1, Propositions 4.2-4.4 still hold for $%
Q_{\varepsilon ,\sigma }$ instead of $Q_{\varepsilon }.$ So we have an
apriori $C^{1}$ bound for solutions of (\ref{eqn4.18}), independent of $%
\sigma $ and $\varepsilon $ (for $0<\varepsilon \leq {\varepsilon}_{0}$). Altogether we have obtained an apriori $%
C^{1,\beta }(\bar{\Omega})$ ($\beta >0)$ bound for solutions of (\ref%
{eqn4.18}), independent of $\sigma $ (but depend on $\varepsilon $). By
Theorem 11.8 in \cite{GT83}, we obtain 

\bigskip

\textbf{Theorem 4.5}. \textit{Let }$\Omega $ \textit{be a p-convex bounded
domain in }$R^{m},m\geq 2,$ \textit{with }$\partial \Omega \in C^{2,\alpha }$\textit{%
\ }$(0<\alpha <1)$\textit{.} \textit{Let }$\varphi \in C^{2,\alpha }(\bar{%
\Omega}).$\textit{\ Suppose }$\vec{F}$\textit{\ }$\in $\textit{\ }$%
C^{1,\alpha }(\bar{\Omega})$\textit{\ satisfies the condition (\ref{eqn4.11}) for }$C^{1,\alpha }$\textit{-smooth and bounded }$f_{K}$\textit{'s
in }$\Omega .$\textit{ Then there exists  
a solution }$u_{\varepsilon }$ $\in $ $%
C^{2,\alpha }(\bar{\Omega})$\textit{ of the Dirichlet problem: }$Q_{\varepsilon
}(u)=0 $\textit{ in }$\Omega ,$ $u=\varphi $\textit{ on }$\partial \Omega $\textit{ for given }$%
\varepsilon $ $>$ $0.$ 

\bigskip

\noindent \textbf{(Proof of Theorem A Continued)}
 
Propositions 4.2-4.4 tell us that there exists a
constant $C$ $=$ $C(\varepsilon ,\ a,\ R,$ $||F_{I}||_{\infty},$ $||\partial _{I}F_{J}||_{\infty},$ $||\varphi ||_{\infty},$ $||\partial _{I}\varphi ||_{\infty},$ $||\partial _{I}\partial _{J}\varphi ||_{\infty},$ $||f_{I}||_{\infty})$ such that

\begin{equation}
\sup_{\Omega }|u_{\varepsilon }|+\sup_{\Omega }|\nabla u_{\varepsilon }|\leq
C.  \label{eqn4.19}
\end{equation}

\noindent Moreover, the bounds hold uniformly for $0<\varepsilon \leq {\varepsilon}_{0}$, a positive constant.

In view of (\ref{eqn4.19}) we can find a subsequence 
$u_{\varepsilon _{j}}$ $(0<\varepsilon _{j}\leq \varepsilon _{0},\ \varepsilon _{j}\rightarrow 0)$ converging to $u_{0}$ in $C^{0}$ by the Arzela-Ascoli theorem.
Then the Lipschitzianity of $u_{0}$ follows by taking the limit
of ratios: ($x\neq y$)

\begin{eqnarray*} 
\left| \frac{u_{\varepsilon _{j}}(x)-u_{\varepsilon _{j}}(y)}{x-y}\right|
\  (\leq C).	
\end{eqnarray*}


Next we claim that $u_{0}$ is a minimizer for $\mathcal{F(\cdot )}$ (see (%
\ref{eqn1.3})) such that $u_{0}$ $=$ $\varphi $ on $\partial \Omega .$
Observe that $W^{1,q}(\Omega )$ is compactly imbedded in $L^{1}(\Omega )$
(e.g., Theorem 7.26 in \cite{GT83}). So we may as well assume that $%
u_{\varepsilon _{j}}$ converges to $u_{0}$ in $L^{1}(\Omega ).$ Also note
that $|\vec{p}+\vec{F}|$ is convex in $\vec{p}$ since $|\lambda \vec{p}_{1}$ 
$+$ $(1-\lambda )\vec{p}_{2}$ $+$ $\vec{F}|$ $=$ $|\lambda (\vec{p}_{1}+\vec{%
F})$ $+$ $(1-\lambda )(\vec{p}_{2}+\vec{F})|$ $\leq $ $\lambda |\vec{p}_{1}+%
\vec{F}|$ $+$ $(1-\lambda )|\vec{p}_{2}+\vec{F}|$ for $0$ $\leq $ $\lambda $ 
$\leq $ $1.$ We can therefore apply Theorem 4.1.2 in \cite{Morrey66} to
conclude the lower semicontinuity of $\mathcal{F(\cdot )}$ (see (\ref{eqn1.3}%
))$:$

\begin{equation}
\mathcal{F(}u_{0}\mathcal{)\leq }\lim \inf_{j\rightarrow \infty }\mathcal{F(}%
u_{\varepsilon _{j}}).  \label{eqn4.20}
\end{equation}

\noindent Now for $v$ $\in $ $W^{1,1}$ with $v$ $-$ $\varphi $ $\in $ $%
W_{0}^{1,1},$ we estimate

\begin{eqnarray}
\mathcal{F(}u_{\varepsilon _{j}}) &\equiv &\int_{\Omega }\mid \nabla
u_{\varepsilon _{j}}+\vec{F}\mid \text{ (omitting volume element)}
\label{eqn4.21} \\
&\leq &\int_{\Omega }\sqrt{\varepsilon _{j}^{2}+|\nabla u_{\varepsilon _{j}}+%
\vec{F}|^{2}}  \notag \\
&\leq &\int_{\Omega }\sqrt{\varepsilon _{j}^{2}+|\nabla v+\vec{F}|^{2}} 
\notag \\
&\leq &\varepsilon _{j}\text{ }vol(\Omega )+\int_{\Omega }|\nabla v+\vec{F}|
\notag
\end{eqnarray}

\noindent where we have used the fact that the Dirichlet solution $%
u_{\varepsilon _{j}}$ $\in$ $C^{2}(\bar{\Omega})$ is also a minimizer for $\mathcal{F}_{\varepsilon
_{j}}(u)$ $\equiv $ $\int_{\Omega }\sqrt{\varepsilon _{j}^{2}+|\nabla u+\vec{%
F}|^{2}}$. Taking the limit infimum of (\ref{eqn4.21}) and making use
of (\ref{eqn4.20}), we finally obtain that $\mathcal{F(}u_{0}\mathcal{)}$ $%
\mathcal{\leq }$ $\int_{\Omega }|\nabla v+\vec{F}|$ $\equiv $ $\mathcal{F(}v%
\mathcal{)}$. That is to say, $u_{0}$ is a minimizer for $\mathcal{F(}\cdot 
\mathcal{)}$.

\begin{flushright}
Q.E.D.
\end{flushright}

\bigskip

\section{Uniqueness of minimizers-proof of Theorems B and C}

Recall (see Section 3) that $\Omega \subset R^{m}$ denotes a bounded
domain and $\mathcal{F}(u)$ $\equiv $ $\int_{\Omega }\{|\nabla u+\vec{F}|$ $%
+ $ $Hu\}$ for $u\in W^{1,1}(\Omega )$, $\vec{F}$ $\in $ $L^{1}(\Omega )$,
and $H$ $\in$ $L^{\infty}(\Omega )$. We will prove two ($W^{1,1}$) minimizers 
for $\mathcal{F}(u)$ with the same
''boundary value'' have the same normal vector ''almostly''.

\bigskip

\textbf{Theorem 5.1. }\textit{Let }$u,v$\textit{\ }$\in $\textit{\ }$%
W^{1,1}(\Omega )$\textit{\ be two minimizers for }$\mathcal{F}(u)$\textit{\
such that }$u-v$\textit{\ }$\in $\textit{\ }$W_{0}^{1,1}(\Omega ).$\textit{\
Let }$u_{\varepsilon }$\textit{\ }$\equiv $\textit{\ }$u+\varepsilon (v-u).$%
\textit{\ Then for any pair of regular }$\varepsilon _{1},$\textit{\ }$%
\varepsilon _{2}$\textit{\ }$\in $\textit{\ }$[0,1],$\textit{\ there holds }$%
N(u_{\varepsilon _{1}})=N(u_{\varepsilon _{2}})$\textit{\ in }$\Omega
\backslash \lbrack S(u_{\varepsilon _{1}})\cup S(u_{\varepsilon _{2}})](a.e.).$

\bigskip

\textbf{Proof. }By (\ref{eqn3.13}) with $\varphi =v-u,$ we have

\begin{equation}
0=\mathcal{F}(v)-\mathcal{F}(u)=\int_{0}^{1}\frac{d\mathcal{F}%
(u_{\varepsilon })}{d\varepsilon }d\varepsilon .  \label{eqn5.1}
\end{equation}

As in the proof of Theorem 3.3, the same argument shows that $\frac{d%
\mathcal{F}(u_{\varepsilon })}{d\varepsilon }\geq 0$ for any regular $%
\varepsilon $ $\in $ $[0,1].$ In view of (\ref{eqn5.1}) and Lemma 3.2(1), $\frac{d\mathcal{F}%
(u_{\varepsilon })}{d\varepsilon }$ $=$ $0$ for any regular $\varepsilon $ $%
\in $ $[0,1].$ It follows from (\ref{eqn3.4}) that $\int_{\Omega \backslash
S(u_{\varepsilon })}N(u_{\varepsilon })\cdot \nabla (v-u)$ $=$ $0$.
Therefore for any pair of regular $\varepsilon _{1},$ $\varepsilon _{2}$ $%
\in $ $[0,1],$ there holds

\begin{equation}
\int_{\Omega \backslash \lbrack S(u_{\varepsilon _{1}})\cup S(u_{\varepsilon
_{2}})]}[N(u_{\varepsilon _{2}})-N(u_{\varepsilon _{1}})]\cdot \nabla
(v-u)=0.  \label{eqn5.2}
\end{equation}

\noindent Here we have used $\int_{S(u_{\varepsilon _{1}})\backslash
S(u_{\varepsilon _{2}})}N(u_{\varepsilon _{2}})\cdot \nabla (v-u)$ $=$ $0$
and $\int_{S(u_{\varepsilon _{2}})\backslash S(u_{\varepsilon
_{1}})}N(u_{\varepsilon _{1}})\cdot \nabla (v-u)$ $=$ $0$ by observing that
for $j=1,2,$ $|N(u_{\varepsilon _{j}})\cdot \nabla (v-u)|$ $\leq $ $|\nabla
(v-u)|$ and $\int_{S(u_{\varepsilon _{j}})}|\nabla (v-u)|=0$ from the
definition of $\varepsilon _{j}$ being regular. Write $v-u$ $=$ $%
(u_{\varepsilon _{2}}-u_{\varepsilon _{1}})/(\varepsilon _{2}-\varepsilon
_{1})$ for $\varepsilon _{2}\neq \varepsilon _{1}.$ By Lemma 5.1' in \cite%
{CHMY04}, the integrand in (\ref{eqn5.2}) is

\begin{equation*}
\frac{|\nabla u_{\varepsilon _{2}}+\vec{F}|+|\nabla u_{\varepsilon _{1}}+%
\vec{F}|}{2(\varepsilon _{2}-\varepsilon _{1})}|N(u_{\varepsilon
_{2}})-N(u_{\varepsilon _{1}})|^{2}.
\end{equation*}

It then follows that $N(u_{\varepsilon _{1}})=N(u_{\varepsilon _{2}})$ in $%
\Omega \backslash \lbrack S(u_{\varepsilon _{1}})\cup S(u_{\varepsilon
_{2}})].$

\begin{flushright}
Q.E.D.
\end{flushright}

\bigskip

For a vector field $\vec{G}$ $=$ $(g_{1},g_{2},...,g_{2n})$ on $%
\Omega \subset R^{2n},$ we recall that $\vec{G}^{\ast }$ $\equiv $ $(g_{2},$ $%
-g_{1},$ $g_{4},$ $-g_{3},$ $...,$ $g_{2n},$ $-g_{2n-1}).$

\bigskip

\textbf{Lemma 5.2. }\textit{Let }$u,v$\textit{\ }$\in $\textit{\ }$%
W^{1}(\Omega )$\textit{\ where the domain }$\Omega $\textit{\ is contained
in }$R^{2n}.$\textit{\ Let }$u_{\varepsilon }$\textit{\ }$\equiv $\textit{\ }%
$u+\varepsilon (v-u).$\textit{\ Suppose }$N(u_{\varepsilon
_{1}})=N(u_{\varepsilon _{2}})$\textit{\ in }$\Omega \backslash \lbrack
S(u_{\varepsilon _{1}})\cup S(u_{\varepsilon _{2}})]$\textit{\ for a pair 
}$\varepsilon _{1},$\textit{\ }$\varepsilon _{2}$\textit{\ such that }$%
\varepsilon _{1}$\textit{\ }$\neq $\textit{\ }$\varepsilon _{2}$\textit{.
Then for }$j=1,2,$\textit{\ there holds}

\begin{equation}
(\nabla u_{\varepsilon _{j}}+\vec{F})^{\ast }\cdot (\nabla v-\nabla u)=0%
\text{ in }\Omega \text{ (a.e.)}\mathit{.}  \label{eqn5.3}
\end{equation}

\bigskip

\textbf{Proof. }We will prove (\ref{eqn5.3}) only for $j=1$ (similar
argument works also for $j=2)$ For $p$ $\in $ $S(u_{\varepsilon _{1}}),$ $%
\nabla u_{\varepsilon _{1}}+\vec{F}=0.$ So (\ref{eqn5.3}) holds obviously.
For $p$ $\in $ $S(u_{\varepsilon _{2}}),$ (\ref{eqn5.3}) also holds by
observing that $\nabla v-\nabla u$ $=$ $[(\nabla u_{\varepsilon _{1}}+\vec{F}
$ $)-$ $(\nabla u_{\varepsilon _{2}}+\vec{F})]/(\varepsilon _{1}-\varepsilon
_{2})$ $=$ $(\nabla u_{\varepsilon _{1}}+\vec{F}$ $)/(\varepsilon
_{1}-\varepsilon _{2})$ and $\vec{G}^{\ast }\cdot \vec{G}$ $=$ $0.$ For the
remaining case: $p$ $\in $ $\Omega \backslash \lbrack S(u_{\varepsilon
_{1}})\cup S(u_{\varepsilon _{2}})],$ we observe that for $j=1,2,$

\begin{equation}
N(u_{\varepsilon _{j}})^{\ast }\cdot \nabla u_{\varepsilon _{j}}=\frac{\vec{F%
}^{\ast }\cdot \nabla u_{\varepsilon _{j}}}{|\nabla u_{\varepsilon _{j}}+%
\vec{F}|}=\vec{F}^{\ast }\cdot N(u_{\varepsilon _{j}}).  \label{eqn5.4}
\end{equation}

Here we have used the property $\vec{G}^{\ast }\cdot \vec{G}$ $=$ $0$ twice.
Since $N(u_{\varepsilon _{1}})=N(u_{\varepsilon _{2}})$ in $\Omega
\backslash \lbrack S(u_{\varepsilon _{1}})\cup S(u_{\varepsilon _{2}})]$ by
assumption (hence $N(u_{\varepsilon _{1}})^{\ast }=N(u_{\varepsilon
_{2}})^{\ast }$ also), we take the difference of (\ref{eqn5.4}) for $j=1$
and $j=2$ to obtain

\begin{equation}
N(u_{\varepsilon _{1}})^{\ast }\cdot (\nabla u_{\varepsilon _{2}}-\nabla
u_{\varepsilon _{1}})=0.  \label{eqn5.5}
\end{equation}

Formula (\ref{eqn5.3}) for $j=1$ on $\Omega \backslash \lbrack
S(u_{\varepsilon _{1}})\cup S(u_{\varepsilon _{2}})]$ then follows from (\ref%
{eqn5.5}) by noting that $v-u$ $=$ $(u_{\varepsilon _{2}}-u_{\varepsilon
_{1}})/(\varepsilon _{2}-\varepsilon _{1}).$

\begin{flushright}
Q.E.D.
\end{flushright}

\bigskip

We will use the following general criterion to prove the uniqueness
of minimizers and a comparison principle for weak functions later.

\bigskip

\textbf{Theorem 5.3.} \textit{Let }$\Omega $\textit{\ be a bounded domain in 
}$R^{2n}.$\textit{\ Let $w\in
W_{0}^{1,p}(\Omega )$, $\sigma \in W^{1,q}(\Omega )$, where $1\leq p<\infty$,
$q=\frac{p}{p-1}$ ($q=\infty$ for $p=1$).}
\textit{\ Let $\vec{F}$ (a vector field) 
$\in$ $W^{1,1}(\Omega )\cap L^{q}(\Omega )$ satisfying $div\vec{F}%
^{\ast }$}\textit{\ }$>$\textit{\ }$0$\textit{\ (a.e.) or }$div\vec{F}^{\ast
} $\textit{\ }$<$\textit{\ }$0$\ \textit{(a.e.).}\textit{\ Suppose }$%
(\nabla \sigma +\vec{F})^{\ast }\cdot \nabla w$\textit{\ }$=$\textit{\ }$0$%
\textit{\ in }$\Omega $\textit{\ (a.e.). Then }$w$\textit{\ }$\equiv $%
\textit{\ }$0$\textit{\ in }$\Omega $\textit{\ (a.e.).}

\bigskip

\textbf{Proof. }Take $\omega _{j}$ $\in C_{0}^{\infty }(\Omega )\rightarrow
w $ in $W^{1,p}$ and $\vec{F}_{\bar{k}}\in C^{\infty }(\Omega )$ $\rightarrow$ 
$\vec{F}$ in $W^{1,1}\cap L^{q}$. Suppose $\omega _{j}$ does not vanish identically. Then
there exists a decreasing sequence of positive numbers $a_{i}$ converging to 
$0$ such that $\Omega _{j,i}$ $\equiv $ $\{|\omega _{j}|$ $>a_{i}\}$ $%
\subset \subset $ $\Omega $ is not empty for large $i$ and $\partial \Omega
_{j,i}$ is $C^{\infty }$ smooth (by Sard's theorem; note that $|\omega _{j}|$
is $C^{\infty }$ smooth where $\omega _{j}$ $\neq $ $0).$ Also we take $%
v_{k} $ $\in C^{\infty }(\Omega )\rightarrow \sigma $ in $W^{1,2}.$ Consider

\begin{equation}
I_{j,i,k,\bar{k}}\equiv \int_{\partial \Omega _{j,i}}|\omega _{j}|\text{ }(\nabla
v_{k}+\vec{F}_{\bar{k}})^{\ast }\cdot \nu  \label{eqn5.6}
\end{equation}

where $\nu $ denotes the boundary normal. We first compute

\begin{eqnarray}
\int_{\partial \Omega _{j,i}}|\omega _{j}|\text{ }(\nabla v_{k}+\vec{F}_{\bar{k}}%
)^{\ast }\cdot \nu &=&a_{i}\int_{\partial \Omega _{j,i}}(\nabla v_{k}+\vec{F}_{\bar{k}}%
)^{\ast }\cdot \nu  \label{eqn5.7} \\
&=&a_{i}\int_{\Omega _{j,i}}div[(\nabla v_{k})^{\ast }+\vec{F}^{\ast }_{\bar{k}}] 
\notag \\
&=&a_{i}\int_{\Omega _{j,i}}div\vec{F}^{\ast }_{\bar{k}}  \notag
\end{eqnarray}

Here we have used Green's theorem for the second equality and $div(\nabla
v_{k})^{\ast }=0$ for the third equality in (\ref{eqn5.7}). It
follows from (\ref{eqn5.7}) that

\begin{equation}
\lim_{i\rightarrow \infty }I_{j,i,k,{\bar{k}}}=0.  \label{eqn5.8}
\end{equation}

On the other hand, a similar reasoning gives

\begin{eqnarray}
I_{j,i,k,{\bar{k}}} &=&\int_{\Omega _{j,i}}\nabla |\omega _{j}|\cdot (\nabla v_{k}+%
\vec{F}_{\bar{k}})^{\ast }+|\omega _{j}|\text{ }div[(\nabla v_{k})^{\ast }+\vec{F}%
^{\ast }_{\bar{k}}]  \label{eqn5.9} \\
&=&\int_{\Omega _{j,i}}\nabla |\omega _{j}|\cdot (\nabla v_{k}+\vec{F}_{\bar{k}}%
)^{\ast }+|\omega _{j}|\text{ }div\vec{F}^{\ast }_{\bar{k}}.  \notag
\end{eqnarray}

Observe that $\cup _{i}\Omega _{j,i}$ $=$ $\{|\omega _{j}|$ $>$ $0\}$ $=$ $%
\Omega \backslash \{\omega _{j}$ $=0\},$ $(\Omega \backslash \{\omega _{j}$ $%
=0\})$$\backslash \Omega _{j,i}$ $=$ $\cup _{l=i}^{\infty }(\Omega _{j,l+1}$
$\backslash \Omega _{j,l})$, $\vec{F}_{\bar{k}}$ $\in$ $W^{1,1}(\Omega )$, and hence

\begin{eqnarray}
&&(\int_{\Omega _{j,i}}-\int_{\Omega \backslash \{\omega _{j}=0\}})\{\nabla
|\omega _{j}|\cdot (\nabla v_{k}+\vec{F}_{\bar{k}})^{\ast }+|\omega _{j}|\text{ }div%
\vec{F}^{\ast }_{\bar{k}}\}  \label{eqn5.10} \\
&=&-\Sigma _{l=i}^{\infty }\int_{\Omega _{j,l+1}\backslash \Omega
_{j,l}}\{\nabla |\omega _{j}|\cdot (\nabla v_{k}+\vec{F}_{\bar{k}})^{\ast }+|\omega
_{j}|\text{ }div\vec{F}^{\ast }_{\bar{k}}\}  \notag \\
&=&-\Sigma _{l=i}^{\infty }\text{ }(I_{j,l+1,k,{\bar{k}}}-I_{j,l,k,{\bar{k}}})=I_{j,i,k,{\bar{k}}}  \notag
\end{eqnarray}

\noindent by (\ref{eqn5.8}). It follows from (\ref{eqn5.9}), (\ref{eqn5.10})
that

\begin{eqnarray}
0 &=&\int_{\Omega \backslash \{\omega _{j}=0\}}\nabla |\omega _{j}|\cdot
(\nabla v_{k}+\vec{F}_{\bar{k}})^{\ast }+|\omega _{j}|\text{ }div\vec{F}^{\ast }_{\bar{k}}
 \notag \\
&=&\int_{\Omega }\nabla |\omega _{j}|\cdot (\nabla v_{k}+\vec{F}_{\bar{k}})^{\ast
}+|\omega _{j}|\text{ }div\vec{F}^{\ast }_{\bar{k}}.  \notag
\end{eqnarray}

\noindent Here we have used $\nabla |\omega _{j}|=0$ if $\omega _{j}=0$
(p.152 in \cite{GT83}). Letting ${\bar{k}}$ $\rightarrow $ $\infty $ in the above
formula gives

\begin{eqnarray}
0 &=&\int_{\Omega }\nabla |\omega _{j}|\cdot
(\nabla v_{k}+\vec{F})^{\ast }+|\omega _{j}|\text{ }div\vec{F}^{\ast }.
\label{eqn5.11} 
\end{eqnarray}

Letting $k$ $\rightarrow $ $\infty $ in the first
term of (\ref{eqn5.11}), we then estimate by using the assumption $\nabla
w\cdot (\nabla \sigma +\vec{F})^{\ast }$ $=$ $0$

\begin{eqnarray}
&&\int_{\Omega }\nabla |\omega _{j}|\cdot (\nabla \sigma +\vec{F})^{\ast }
\label{eqn5.12} \\
&=&\int_{\{\omega _{j}>0\}}(\nabla \omega _{j}-\nabla w)\cdot (\nabla \sigma
+\vec{F})^{\ast }-\int_{\{\omega _{j}<0\}}(\nabla \omega _{j}-\nabla w)\cdot
(\nabla \sigma +\vec{F})^{\ast }  \notag \\
&\longrightarrow &0\text{ \ \ \ \ }as\text{ \ }j\rightarrow \infty .  \notag
\end{eqnarray}

\noindent Here we have used $\omega _{j}$ $\rightarrow w$ in $W^{1,p}$ and $%
(\nabla \sigma +\vec{F})^{\ast }$ $\in L^{q}(\Omega )$ by assumption. For
the second term of (\ref{eqn5.11}), we have 
\begin{equation}
\lim_{j\rightarrow \infty }\int_{\Omega }|\omega _{j}|\text{ }div\vec{F}%
^{\ast }=\int_{\Omega }|w|\text{ }div\vec{F}^{\ast }>0\text{ or }<0
\label{eqn5.13}
\end{equation}

\noindent if $w\neq 0$ (noting that $div\vec{F}^{\ast }>0$ or $<0$ by
assumption)$.$ By (\ref{eqn5.11}), (\ref{eqn5.12}), and (\ref{eqn5.13}), we
reach a contradiction. Therefore $w\equiv 0$ in $\Omega $ (a.e.).

\begin{flushright}
Q.E.D.
\end{flushright}

\bigskip

\textbf{Remark.} If $\vec{F}$ does not satisfy the condition in Theorem
5.3, then the theorem may not hold as shown by the following examples. Let $%
\Omega $ $=$ $(0,\pi )$ $\times $ $(0,\pi )$ $\subset $ $R^{2}.$ Let $w$ $=$ 
$\sin x\sin y$ $\in $ $W_{0}^{1,2}.$ Then $\nabla w$ $=$ $(\cos x\sin y,$ $%
\sin x\cos y).$ Take $\sigma $ $=$ $0$ and $\vec{F}$ $=$ $(\cos x\sin y,$ $%
\sin x\cos y).$ It is easy to see that $\vec{F}^{\ast }$ $=$ $(\sin x\cos y,$
$-\cos x\sin y),$ $div\vec{F}^{\ast }$ $=$ $0,$ and $(\nabla \sigma +\vec{F}%
)^{\ast }\cdot \nabla w$\textit{\ }$=$\textit{\ }$\vec{F}^{\ast }\cdot
\nabla w$\textit{\ }$=$ $0.$ With the same $\sigma $ $(=0)$ and $w$ as
above, we can also take $\vec{F}$ $=$ $\sin x$ $(\cos x\sin y,$ $\sin x\cos
y).$ Then still $(\nabla \sigma +\vec{F})^{\ast }\cdot \nabla w$\textit{\ }$%
= $\textit{\ }$\vec{F}^{\ast }\cdot \nabla w$\textit{\ }$=$ $0$ while $div%
\vec{F}^{\ast }$ $=$ $\cos x$ $\sin x$ $\cos y$ has no definite sign in $%
\Omega .$

\bigskip

\textbf{Proof of Theorem B. }

The proof follows from Theorem 5.1, Lemma 5.2, and Theorem 5.3 with $p=q=2$, $%
\sigma $ $=$ $u_{\varepsilon _{1}}$, and $w$ $=$ $v-u$.

\begin{flushright}
Q.E.D.
\end{flushright}

\bigskip

Next we want to prove a comparison principle for weak sub- and super- solutions (a
comparison principle for $C^{2}$-smooth functions has been studied in \cite{CHMY04}%
. See Theorem C and Theorem C' there). First we need to define relevant
differential inequalities in some weak sense. Let $\Omega \subset R^{m}$
denote a bounded domain. Recall that $N(u)\equiv $ $\frac{\nabla u+\vec{F}}{%
|\nabla u+\vec{F}|}$ is defined on $\Omega \backslash S(u)$ ($\vec{F}$, say,
is an $L^{1}_{loc}$ vector field in $\Omega ).$

\bigskip

\textbf{Definition 5.1. }Let $H$ $\in$ $L^{1}_{loc}(\Omega )$. We say $u\in W^{1,1}(\Omega )$ satisfies $%
divN(u)\geq H$ ($\leq H,$ respectively) in the weak sense in $\Omega $ if
and only if for any $\varphi \in C_{0}^{\infty}(\Omega )$ and $\varphi \geq 0,$
there holds

\begin{eqnarray}
-\int_{S(u)}|\nabla \varphi |+\int_{\Omega \backslash S(u)}N(u)\cdot \nabla
\varphi +\int_{\Omega }H\varphi &\leq &0  \label{eqn5.14} \\
(\int_{S(u)}|\nabla \varphi |+\int_{\Omega \backslash S(u)}N(u)\cdot \nabla
\varphi +\int_{\Omega }H\varphi &\geq &0,\text{ respectively).}
\label{eqn5.15}
\end{eqnarray}

\bigskip

Recall that we defined the weak solution to $divN(u)$ $=$ $H$ in
Section 3 (see (\ref{eqn3.12})). The following result justifies the above
definitions.

\bigskip

\textbf{Proposition 5.4. }\textit{Let }$H$ $\in$ $L^{\infty}(\Omega ).$\textit{ Then }$u\in W^{1,1}(\Omega )$\textit{\ satisfies }$%
divN(u)$\textit{\ }$\geq $ $H$\textit{\ and }$divN(u)$\textit{\ }$\leq $ $H$%
\textit{\ in the weak sense if and only if }$u\in W^{1,1}(\Omega )$\textit{\
is a weak solution to the equation }$divN(u)$\textit{\ }$=$\textit{\ }$H.$

\bigskip

\textbf{Proof. }Since $C_{0}^{\infty}(\Omega )$ is dense in $W_{0}^{1,1}(\Omega )$,
(\ref{eqn5.14}) and (\ref{eqn5.15}) hold for every $\varphi \in C_{0}^{\infty}(\Omega )$ if and only if they hold for every $\varphi \in W_{0}^{1,1}(\Omega )$. Write $\varphi =\varphi ^{+}-\varphi ^{-}$ for $\varphi \in
W_{0}^{1,1}(\Omega )$ where $\varphi ^{+}$ $\equiv $ $\max \{\varphi ,0\}$
and $\varphi ^{-}$ $\equiv $ $\max \{-\varphi ,0\}.$ Express

\begin{eqnarray}
&&\int_{S(u)}|\nabla \varphi |+\int_{\Omega \backslash S(u)}N(u)\cdot \nabla
\varphi +\int_{\Omega }H\varphi  \label{eqn5.16} \\
&=&\{\int_{S(u)}|\nabla \varphi ^{+}|+\int_{\Omega \backslash S(u)}N(u)\cdot
\nabla \varphi ^{+}+\int_{\Omega }H\varphi ^{+}\}  \notag \\
&&-\{-\int_{S(u)}|\nabla \varphi ^{-}|+\int_{\Omega \backslash
S(u)}N(u)\cdot \nabla \varphi ^{-}+\int_{\Omega }H\varphi ^{-}\}.  \notag
\end{eqnarray}

Note that $\varphi ^{+}\geq 0$ and $\varphi ^{-}\geq 0.$ Now suppose $u$ is
a weak solution to $divN(u)$\textit{\ }$\leq $ $H$\textit{\ and }$divN(u)$%
\textit{\ }$\geq $ $H.$ Then the right-hand side of (\ref{eqn5.16}) is
nonnegative by our definitions. So the left hand side of (\ref{eqn5.16}) is
nonnegative, i.e., (\ref{eqn3.12}) holds. Conversely, suppose $u$ is a weak
solution to $divN(u)$\textit{\ }$=$\textit{\ }$H.$ That is to say, the left
hand side of (\ref{eqn5.16}) is nonnegative (note that $\varphi $ is not
restricted to be nonnegative here). By taking $\varphi \geq 0$ i.e. $\varphi
^{-}=0$ ($\varphi \leq 0$ i.e. $\varphi ^{+}=0,$respectively$)$ in (\ref%
{eqn5.16}), we obtain (\ref{eqn5.15}) ((\ref{eqn5.14}), respectively).

\begin{flushright}
Q.E.D.
\end{flushright}

\bigskip

\textbf{Definition 5.2. }$u,v\in W^{1}(\Omega )$ satisfy $divN(u)\geq
divN(v)$ in $\Omega $ in the weak sense if and only if for any $\varphi \in
W_{0}^{1,1}(\Omega )$ and $\varphi \geq 0,$ there holds

\begin{equation}
-\int_{S(u)}|\nabla \varphi |+\int_{\Omega \backslash S(u)}N(u)\cdot \nabla
\varphi \leq +\int_{S(v)}|\nabla \varphi |+\int_{\Omega \backslash
S(v)}N(v)\cdot \nabla \varphi .  \label{eqn5.17}
\end{equation}

\bigskip

\textbf{Definition 5.3. }$u,v\in W^{1,1}(\Omega )$ satisfy $u\leq v$ on $%
\partial \Omega $ if and only if $(u-v)^{+}$ $\equiv $ $\max (u-v,0)$ $\in
W_{0}^{1,1}(\Omega ).$

\bigskip

\textbf{Theorem 5.5. }\textit{Suppose }$u,v\in W^{1,1}(\Omega )$\textit{\
satisfy the following conditions:}

\begin{eqnarray*}
divN(u) &\geq &divN(v)\text{ in }\Omega \text{ (in the weak sense);} \\
u &\leq &v\text{ on }\partial \Omega .
\end{eqnarray*}

\noindent \textit{Then }$N(u)=N(v)$\textit{\ on }$\{u>v\}\backslash \lbrack S(u)\cup
S(v)].$

\bigskip

\textbf{Proof. }Let $\varphi =(u-v)^{+}.$ The condition $u\leq v$ on $%
\partial \Omega $ implies that $\varphi $ $\in $ $W_{0}^{1,1}(\Omega ).$ Let 
$v_{\varepsilon }\equiv v+\varepsilon \varphi .$ From Lemma 3.2 (1), $\frac{d%
\mathcal{F}(v_{\varepsilon })}{d\varepsilon }$ is increasing in regular $%
\varepsilon .$ It follows that $\frac{d\mathcal{F}(v_{0+})}{d\varepsilon }%
\leq \frac{d\mathcal{F}(v_{1-})}{d\varepsilon }$ by Lemma 3.2 (2). In view
of the formula (\ref{eqn3.3}), we have

\begin{equation}
+\int_{S(v)}|\nabla \varphi |+\int_{\Omega \backslash S(v)}N(v)\cdot \nabla
\varphi \leq -\int_{S(v_{1})}|\nabla \varphi |+\int_{\Omega \backslash
S(v_{1})}N(v_{1})\cdot \nabla \varphi .  \label{eqn5.18}
\end{equation}

Observe that $v_{1}=u$ on $\{u>v\}$ and $\varphi =0$ on $\{u\leq v\}.$ So
the right hand side of (\ref{eqn5.18}) equals the left hand side of (\ref%
{eqn5.17}). It follows that

\begin{equation}
-\int_{S(u)}|\nabla \varphi |+\int_{\Omega \backslash S(u)}N(u)\cdot \nabla
\varphi =+\int_{S(v)}|\nabla \varphi |+\int_{\Omega \backslash
S(v)}N(v)\cdot \nabla \varphi .  \label{eqn5.19}
\end{equation}

Write%
\begin{eqnarray}
&&\int_{\Omega \backslash S(u)}N(u)\cdot \nabla \varphi -\int_{\Omega
\backslash S(v)}N(v)\cdot \nabla \varphi  \label{eqn5.20} \\
&=&\int_{\Omega \backslash \lbrack S(u)\cup S(v)]}(N(u)-N(v))\cdot \nabla
\varphi  \notag \\
&&+\int_{S(v)\backslash S(u)}N(u)\cdot \nabla \varphi -\int_{S(u)\backslash
S(v)}N(v)\cdot \nabla \varphi .  \notag
\end{eqnarray}

We claim 
\begin{equation}
-\int_{S(u)}|\nabla \varphi |-\int_{S(u)\backslash S(v)}N(v)\cdot \nabla
\varphi =0.  \label{eqn5.21}
\end{equation}

Since $\varphi =0$ on $\{u\leq v\},$ we only have to discuss the case that $%
u>v.$ In this case, $\varphi =u-v$ and hence $\nabla \varphi $ $=$ ($\nabla
u $ $+$ $\vec{F})$ $-$ $(\nabla v$ $+$ $\vec{F})$ $=$ $-(\nabla v$ $+$ $\vec{%
F})$ in $S(u)$ (and $=0$ in $S(u)\cap S(v))$. So $N(v)\cdot \nabla \varphi $ 
$= $ $\frac{\nabla v+\vec{F}}{|\nabla v+\vec{F}|}\cdot \lbrack -(\nabla v$ $%
+ $ $\vec{F})]$ $=$ $-|\nabla v$ $+$ $\vec{F}|$ in $S(u)\backslash S(v).$ It
is now clear that (\ref{eqn5.21}) holds. Similarly there also holds

\begin{equation}
-\int_{S(v)}|\nabla \varphi |+\int_{S(v)\backslash S(u)}N(u)\cdot \nabla
\varphi =0.  \label{eqn5.22}
\end{equation}

Combining (\ref{eqn5.19}), (\ref{eqn5.20}), (\ref{eqn5.21}), (\ref{eqn5.22})
gives

\begin{equation}
\int_{\Omega \backslash \lbrack S(u)\cup S(v)]}(N(u)-N(v))\cdot \nabla
\varphi =0.  \label{eqn5.23}
\end{equation}

By Lemma 5.1' in \cite{CHMY04} (which works also for $u,v$ $\in $ $W^{1,1}(\Omega
) $), we have

\begin{equation}
(N(u)-N(v))\cdot \nabla \varphi =\frac{|\nabla u+\vec{F}|+|\nabla v+\vec{F}|%
}{2}|N(u)-N(v)|^{2}  \label{eqn5.24}
\end{equation}

on $\{u>v\}\backslash \lbrack S(u)\cup S(v)]$ (where $\varphi =u-v$)$.$
Noting that $\varphi =0$ on $\{u\leq v\}$ and substituting (\ref{eqn5.24})
into (\ref{eqn5.23}), we finally obtain $N(u)$ $=$ $N(v)$ on $%
\{u>v\}\backslash \lbrack S(u)\cup S(v)].$

\begin{flushright}
Q.E.D.
\end{flushright}

\bigskip

We can now prove the comparison principle for weak sub- and super- solutions.

\bigskip

\textbf{Proof of Theorem C. }

By Theorem 5.5 and Lemma 5.2 (switching the roles of $u$ and $v$ and
taking $\Omega $ $=$ $\{u>v\}$, $\varepsilon _{1}=0,$ $\varepsilon _{2}=1$),
we obtain ($\nabla v+\vec{F})^{\ast }$ $\cdot $ $\nabla (u-v)^{+}$ $=$ $0$.
Then we apply Theorem 5.3 (with $p=q=2$, $\sigma $ $=$ $v$, and $w$ $=$ $(u-v)^{+}$)
to conclude that $(u-v)^{+}$ $=$ $0$ in $\Omega .$ That is to say, $u\leq v$%
\textit{\ }in\textit{\ }$\Omega .$

\begin{flushright}
Q.E.D.
\end{flushright}

\bigskip

\section{When a smooth solution is a minimizer}

In this section we determine when a smooth solution is a minimizer. We will
prove Theorem D, Theorem E, and Corollary F. We first prove a result for the case
$H_{m-1}(S(u))$ $=$ $0$, in which a $C^{2}$-smooth solution must be a weak solution.

\bigskip

\textbf{Lemma 6.1}. \textit{Let }$\Omega $\textit{\ be a bounded domain in }$%
R^{m}.$\textit{\ Suppose }$u$\textit{\ }$\in $\textit{\ }$C^{1}(\Omega )$%
\textit{\ }$\cap $\textit{\ }$C^{2}(\Omega \backslash S(u))$\textit{\ }$\cap 
$\textit{\ }$C^{0}(\bar{\Omega})$\textit{\ satisfies (\ref{eqn1.4'}) in }$%
\Omega \backslash S(u)$\textit{\ with }$\vec{F}$ $\in$ $C^{1}(\Omega \backslash S(u))$\textit{ and }$H$\textit{\ }$\in $\textit{\ }$%
C^{0}(\Omega \backslash S(u))$\textit{\ }$\cap $\textit{\ }$L^{1}_{loc}
(\Omega ).$\textit{\ Suppose }$H_{m-1}(S(u)),$\textit{\ the }$m-1$\textit{\
dimensional Hausdorff measure of }$S(u),$\textit{\ vanishes. Then }$u$%
\textit{\ is a weak solution to (\ref{eqn1.4'}) and a minimizer for (\ref{eqn1.3'}) 
if }$u$ $\in $ $W^{1,1}(\Omega )$\textit{ and }$H$ $\in$ 
$L^{\infty}(\Omega )$\textit{ also.}

\textit{\bigskip }

\textbf{Proof.} By Theorem 3.3, it suffices to prove that for any $\varphi $ 
$\in $ $C_{0}^{\infty }(\Omega )$ (3.12) holds. That is,

\begin{equation*}
\int_{S(u)}|\nabla \varphi |+\int_{\Omega \backslash S(u)}N(u)\cdot \nabla
\varphi +\int_{\Omega }H\varphi \geq 0.
\end{equation*}

Write $\Omega $ $%
=$ $\Omega _{+}$ $\cup $ $\Omega _{0}$ $\cup $ $\Omega _{-}$ where $\Omega
_{+}$ $\equiv $ $\{\varphi $ $>$ $0\},$ $\Omega _{-}$ $\equiv $ $\{\varphi $ 
$<$ $0\},$ and $\Omega _{0}$ $\equiv $ $\{\varphi $ $=$ $0\}.$ If $\Omega
_{+}$ $\neq $ $\emptyset ,$ then there exists a sequence of $\varepsilon
_{j} $ $>$ $0$ approaching $0,$ such that $\Omega _{\varepsilon _{j}}$ $%
\equiv $ $\{\varphi $ $>$ $\varepsilon _{j}\}$ $\neq $ $\emptyset ,$ $\cup
_{j=1}^{\infty }\Omega _{\varepsilon _{j}}$ $=$ $\Omega _{+}$ and $\partial
\Omega _{\varepsilon _{j}}$ are $C^{\infty }$-smooth by Sard's theorem.
Since $u$\textit{\ }$\in $\textit{\ }$C^{1}(\Omega ),$ $S(u)$ $\cap $ $\bar{%
\Omega}_{\varepsilon _{j}}$ is compact. Together with the condition $%
H_{m-1}(S(u))$ $=$ $0,$ for any $\alpha $ $>$ $0,$ we can find a finite
cover of balls $B_{r_{k}}(p_{k})$ of center $p_{k}$ and radius $r_{k},$ $%
k=1, $ $2$, ...,$K$ for $S(u)$ $\cap $ $\bar{\Omega}_{\varepsilon _{j}}$
such that

\begin{equation}
\sum_{k=1}^{K}H_{m-1}(\partial B_{r_{k}}(p_{k}))<\alpha .  \label{eqn7.1}
\end{equation}

On the other hand we compute by the divergence theorem and the equation (\ref%
{eqn1.4'})

\begin{eqnarray}
&&\int_{\partial (\Omega _{\varepsilon _{j}}\backslash \cup
B_{r_{k}}(p_{k}))}(\varphi -\varepsilon _{j})N(u)\cdot \nu  \label{eqn7.2} \\
&=&\int_{\Omega _{\varepsilon _{j}}\backslash \cup B_{r_{k}}(p_{k})}\nabla
\varphi \cdot N(u)+(\varphi -\varepsilon _{j})H.  \notag
\end{eqnarray}

\noindent Since $\varphi -\varepsilon _{j}$ $=$ $0$ on $\partial \Omega
_{\varepsilon _{j}},$ we can estimate the boundary term in (\ref{eqn7.2}) as
follows:

\begin{eqnarray}
&\mid &\int_{\partial (\Omega _{\varepsilon _{j}}\backslash \cup
B_{r_{k}}(p_{k}))}(\varphi -\varepsilon _{j})N(u)\cdot \nu \mid
\label{eqn7.3} \\
&\leq &\{\max_{\Omega }|\varphi -\varepsilon _{j}|\}H_{m-1}(\cup
_{k=1}^{K}\partial B_{r_{k}}(p_{k}))  \notag \\
&\leq &\alpha \max_{\Omega }|\varphi -\varepsilon _{j}|  \notag
\end{eqnarray}

\noindent by (\ref{eqn7.1}) and the fact that $|N(u)|$ $=$ $|\nu |$ $=$ $1.$ Letting $%
\alpha \rightarrow 0$ in (\ref{eqn7.3}) gives

\begin{equation}
\int_{\Omega _{\varepsilon _{j}}\backslash S(u)}\nabla \varphi \cdot
N(u)+(\varphi -\varepsilon _{j})H=0  \label{eqn7.4}
\end{equation}

\noindent in view of (\ref{eqn7.2}). Letting $\varepsilon _{j}$ $\rightarrow 
$ $0$ in (\ref{eqn7.4}), we obtain

\begin{equation}
\int_{\Omega _{+}\backslash S(u)}\nabla \varphi \cdot N(u)+\int_{\Omega
_{+}}\varphi H=0  \label{eqn7.5}
\end{equation}

\noindent by noting that the volume of $\{0$ $<$ $\varphi $ $\leq $ $%
\varepsilon _{j}\}$ tends to $0$ as $\varepsilon _{j}$ $\rightarrow $ $0.$
Similarly we also have%
\begin{equation}
\int_{\Omega _{-}\backslash S(u)}\nabla \varphi \cdot N(u)+\int_{\Omega
_{-}}\varphi H=0.  \label{eqn7.6}
\end{equation}

\noindent On the other hand, it is obvious that the integral of $\varphi H$
over $\Omega _{0}$ vanishes since $\varphi $ $=$ $0$ on $\Omega _{0}.$
Observing that $\nabla \varphi $ $=$ $0$ a.e. on $\Omega _{0}$ in view of
Lemma 7.7 in \cite{GT83}, we conclude that%
\begin{equation}
\int_{\Omega _{0}\backslash S(u)}\nabla \varphi \cdot N(u)=0.  \label{eqn7.7}
\end{equation}

\noindent It now follows from (\ref{eqn7.5}), (\ref{eqn7.6}), and (\ref%
{eqn7.7}) that 
\begin{equation}
\int_{\Omega \backslash S(u)}\nabla \varphi \cdot N(u)+\int_{\Omega }\varphi
H=0  \label{eqn7.8}
\end{equation}

\noindent for $\varphi $ $\in $ $C_{0}^{\infty }(\Omega ).$ Comparing (\ref{eqn7.8}) with
(\ref{eqn3.12}) and noting that the first integral of (\ref{eqn3.12}) is zero by $%
H_{m-1}(S(u)) $ $=$ $0,$ we have completed the proof.

\begin{flushright}
Q.E.D.
\end{flushright}

\bigskip

\textbf{Proof of Theorem D.}

Write $\nabla u$ $+$ $\vec{F}$ $=$ $%
(u_{I}+F_{I})_{I=1}^{m}.$ Consider the map $G:$ $p\in \Omega $ $\rightarrow $
$((u_{I}+F_{I})(p))_{I=1}^{m}.$ Computing the differential $dG$ of $G$ at a
singular point $p$ (where $G(p)$ $=$ $0),$ we obtain $(\partial _{J}u_{I}$ $+
$ $\partial _{J}F_{I})$ in matrix form (note that $G$ $\in$ $C^{1}$). From
elementary linear algebra we compute

\begin{eqnarray}
&&rank\text{ }(\partial _{J}u_{I}+\partial _{J}F_{I})+rank\text{ }(\partial
_{I}u_{J}+\partial _{I}F_{J})  \label{eqn7.9} \\
&\geq &rank\text{ }\{(\partial _{J}u_{I}+\partial _{J}F_{I})-(\partial
_{I}u_{J}+\partial _{I}F_{J})\}  \notag \\
&=&rank\text{ }(\partial _{J}F_{I}-\partial _{I}F_{J}).  \notag
\end{eqnarray}

\noindent Observing that $rank$ $(\partial _{J}u_{I}+\partial _{J}F_{I})$ $=$
$rank$ $(\partial _{I}u_{J}+\partial _{I}F_{J})$ (the transpose has the same
rank)$,$ we can deduce from (\ref{eqn7.9}) that $rank$ $dG(p)$ $\geq $ $%
\mathit{[}\frac{rank\text{ }(h_{JI}(p))+1}{2}\mathit{]}$ where $h_{JI}$ $%
\equiv $ $(\partial _{J}F_{I}-\partial _{I}F_{J}).$ It follows that

\begin{equation}
dim(Ker\text{ }dG(p))\leq m-\mathit{[}\frac{rank\text{ }(h_{JI}(p))+1}{2}\mathit{%
].}  \label{eqn7.10}
\end{equation}

\noindent Then by the implicit function theorem there exists an open
neighborhood $V$ of $p$ in $\Omega$ such that $G^{-1}(0)\cap V=S(u)\cap V$ is a
submanifold of $V,$ having (Euclidean) dimension $dim_{E}$ bounded by the right side of (\ref{eqn7.10}).

\begin{flushright}
Q.E.D.
\end{flushright}

\bigskip 

\textbf{Proof of Theorem E.}

It suffices to prove that $H_{m-1}(S(u))$ $=$ $0$ in view of Lemma 6.1. Combining (\ref{eqn1.5'}) and (\ref{eqn1.6}), we bound
$dim_{E}S(u)$ by $m-2.$ It follows that $H_{m-1}(S(u))$ $=$ $0$. 

\begin{flushright}
Q.E.D.
\end{flushright}

\textbf{Proof of Corollary F.}

For $m=2n$, $\vec{F}=-\vec{X}^{\ast }$, we compute $rank\text{ }(h_{JI})=2n$. 
Therefore $(\ref{eqn1.6})$ is reduced to $n\geq 2$, hence $m\geq 4$. 

\begin{flushright}
Q.E.D.
\end{flushright}

We remark that the condition $(\ref{eqn1.6})$ does not hold in dimension $m=2$. So 
$H_{1}(S(u))$ may not vanish. Therefore a $C^2$-smooth solution may not be a minimizer
in this case (see Example 7.4). We will discuss the general situation that $H_{m-1}(S(u))$ $>$ $0$ below.

First we will give a criterion for, in particular, a $C^{2}$-smooth solution to be a minimizer. 
Let $\Omega $ be a domain in $R^{m}$. Let $\Gamma $ $\subset $ $
\Omega $ be a $m-1$ dimensional, orientable, $C^{1}$-smooth submanifold. 
Let $B$ $\subset\subset $ $\Omega 
$ be an open neighborhood of a point in $\Gamma $ with $C^{1}$-smooth
boundary and $\bar{B}$ being compact. Suppose $\Gamma \cap B$ divides $B$ into two disjoint
parts (note that $\Gamma $ may or may not contain some singular
points). That is, $B\backslash \Gamma $ $=$ $B\backslash $ $(\Gamma \cap B)$ 
$=$ $B^{+}\cup B^{-}$ where $B^{+}$ and $B^{-}$ are disjoint domains (proper
open and connected) (see Figure 1(a) or Figure 1(b) below). Suppose $u$
is $C^{2}$-smooth in $\Omega \backslash \Gamma $ and has no singular points in
$\Omega \backslash \Gamma $. Let $\vec{F}$ $\in$ $C^{1}(\Omega )$ for simplicity. Suppose also $N^{+}(u)$ and $%
N^{-}(u) $ (restrictions of $N(u)$ to $B^{+}$ and $B^{-},$ respectively) are
continuous up to $\Gamma \cap B$, i.e., $N^{+}(u)$ $\in$ $C^{0}(\bar{B^{+}})$,
$N^{-}(u)$ $\in$ $C^{0}(\bar{B^{-}})$, so that $divN^{\pm }(u)$ $=$ $H$ in $B^{\pm
},$ respectively$.$ Let $\nu ^{+}$ and $\nu ^{-}$ denote the outward unit
normals to $\Gamma \cap B$ with respect to $B^{+}$ and $B^{-},$ respectively$%
.$ Note that $\nu ^{+}$ $=$ $-\nu ^{-}.$

\begin{figure}[ht]
\begin{center}
\mbox{
\subfigure[]{\includegraphics[height= 4.4cm]{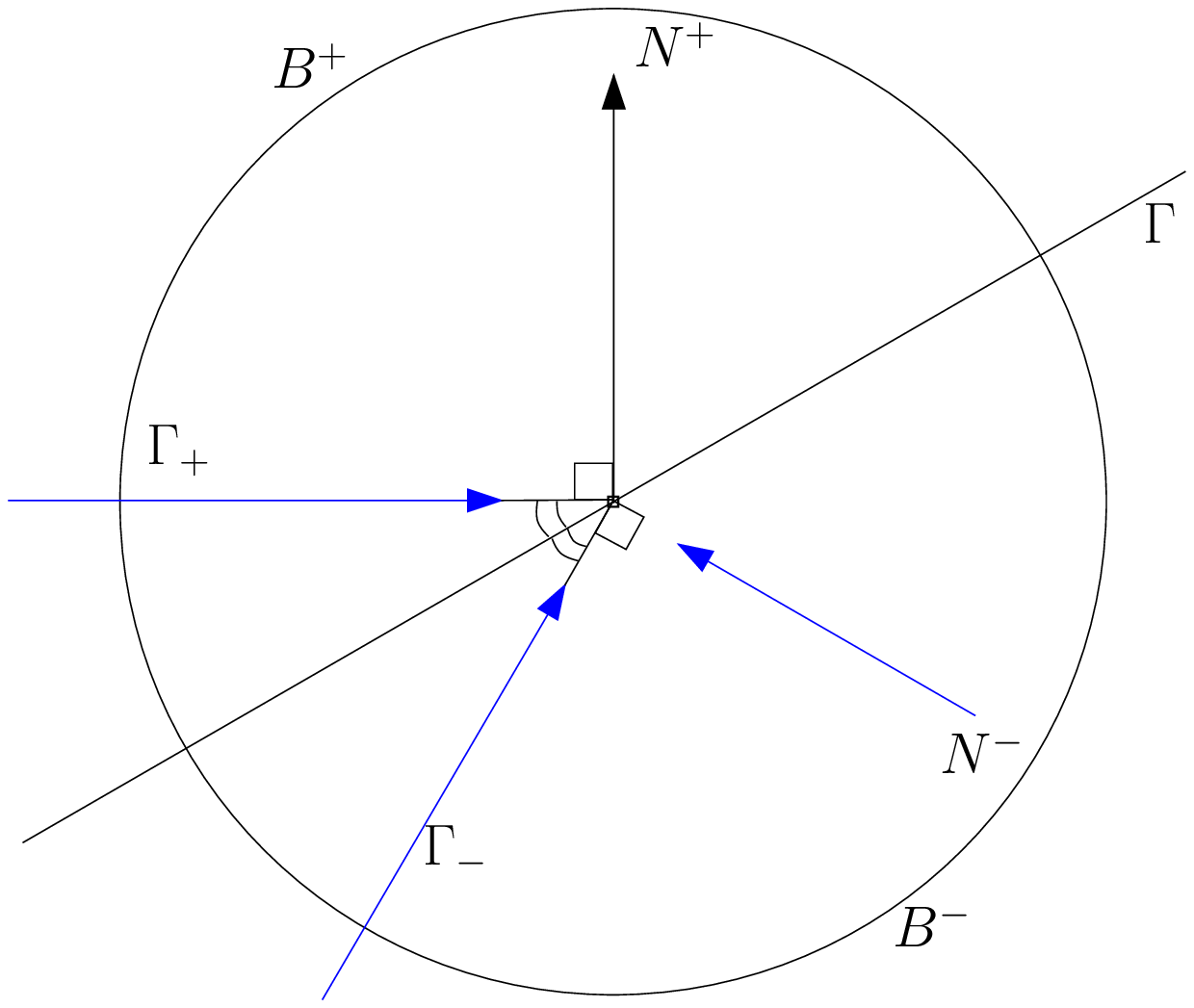}}
\subfigure[]{\includegraphics[height= 4.4cm]{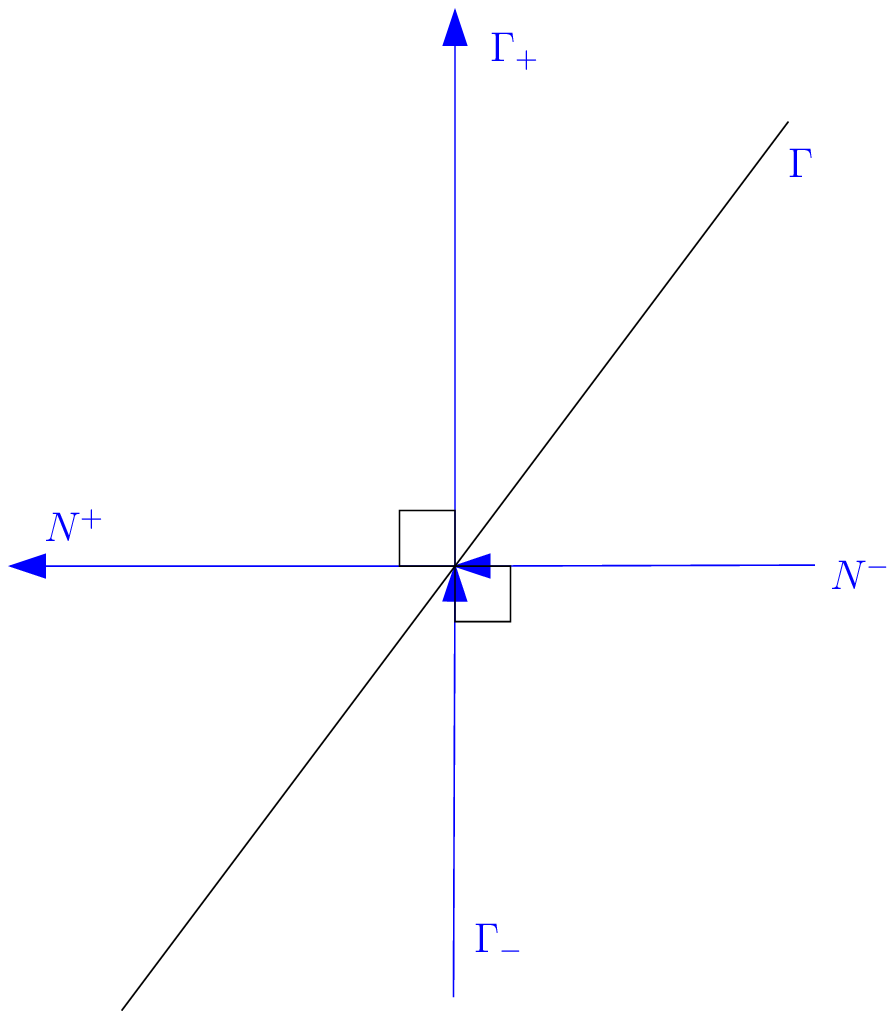}} }
\caption{}
\end{center}
\end{figure}

\textbf{Proposition 6.2.}\textit{\ Suppose we have the situation 
described above. Then }$u$\textit{\ is a weak solution to (\ref{eqn1.4'})}\textit{ on }$B$%
\textit{\ with }$H$\textit{\ }$\in$ $C^{0}(B\backslash \Gamma )$ $\cap$ $L^{\infty}(B)$\textit{\ if
and only if along }$\Gamma \cap B,$ \textit{there holds}

\begin{equation}
\mathit{(N}^{+}\mathit{(u)-N}^{-}\mathit{(u))\cdot \nu }^{+}\mathit{=(N}^{+}%
\mathit{(u)-N}^{-}\mathit{(u))\cdot \nu }^{-}\mathit{=0.}  \label{eqn6.1}
\end{equation}

\bigskip

Note that for $u$ $\in$ $W^{1,1}(B)$, $u$ is a weak solution to (\ref{eqn1.4'}) if
and only if $u$ is a minimizer for (\ref{eqn1.3'}) in view of Theorem 3.3.

\medskip
  
\textbf{Proof}. Using the divergence theorem, we compute

\begin{eqnarray}
&&\int_{B\backslash \Gamma }N(u)\cdot \nabla \varphi +H\varphi 
=(\int_{B^{+}}+\int_{B^{-}})(N(u)\cdot \nabla \varphi +H\varphi ) \label{eqn6.3} \\
&=&\int_{\partial B^{+}}\varphi N^{+}(u)\cdot \nu ^{+}+\int_{\partial
B^{-}}\varphi N^{-}(u)\cdot \nu ^{-}  \notag \\
&=&\int_{\Gamma \cap B}\varphi (N^{+}(u)-N^{-}(u))\cdot \nu ^{+}.  \notag
\end{eqnarray}

\noindent Here we have used $\nu ^{-}$ $=$ $-\nu ^{+}$ and $divN(u)$ $=$ $H$
in both $B^{+}$ and $B^{-}.$ Observing that $H_{m}(S(u)\cap B)$ $=$ $0$ since $H_{m}(S(u)\cap B)$
$\leq$ $H_{m}(\Gamma \cap B)$ $=$ $0$, we conclude from (\ref{eqn3.12}) 
(also $\varphi $ replaced by $-\varphi )$
that $u$ is a weak solution to (\ref{eqn1.4'}) if and only if
 
\begin{equation}
\int_{B\backslash \Gamma }N(u)\cdot \nabla \varphi +H\varphi =0  \label{eqn6.2}
\end{equation}

\noindent for all $\varphi $ $\in $ $C_{0}^{\infty}(B)$. On the other hand, (\ref{eqn6.2}) 
holds if and only if (\ref{eqn6.1}) holds by (\ref{eqn6.3}).

\begin{flushright}
Q.E.D.
\end{flushright}
 
\bigskip

In order to have a criterion for a more general situation, we extend 
Proposition 6.2 as follows. Let $\Omega $ $\subset $ $R^{m}
$ be a bounded domain. Let $A$ $\subset $ $\Gamma $ $\subset $ $\Omega $ such
that $\Gamma $ is relatively closed in $\Omega $, $H_{m-1}(\bar{A})$ $=$ $0$,
and $\Gamma \backslash A$ is a $C^{1}$-smooth $m-1$ 
dimensional manifold. Suppose $\Omega \backslash \Gamma$ $=$ 
${\cup}_{j=1}^{\infty}{\Omega}_{j}$, the union of at most countably many domains 
${\Omega}_{j}$. For each $j$, we have $\partial{\Omega}_{j}$
$\subset$ $\partial\Omega \cup \Gamma$. We can view $\Omega\backslash\Gamma$ 
as domains ${\Omega}_{j}$ obtained by cutting apart along ${\Gamma}$ and 
$\Gamma \backslash A$ as the union of two copies of $\Gamma \backslash A$.
Let $\nu _{j}$ denote the outward unit normal to $\partial{\Omega}_{j}$.
Then $\nu _{j}$ exists for any point $p$ $\in$ $\partial{\Omega}_{j}$ $\cap$
$(\Gamma \backslash A)$. At $p$, there is another $l$ ($l$
may equal $j$) such that $\nu _{l}$ $=$ $-\nu _{j}$. Let $\vec{F}$ $\in$ 
$C^{1}(\Omega \backslash \Gamma)$ for simplicity and $H$ $\in$ 
$C^{0}(\Omega\backslash\Gamma )$ $\cap$ $L^{1}_{loc}(\Omega )$. Suppose $u$ 
$\in$ $C^{1}(\Omega\backslash\Gamma )$ has 
no singular points in $\Omega\backslash\Gamma$. Let $N_{j}(u)$ denote the 
restriction of $N(u)$ on ${\Omega}_{j}$.   

\bigskip 

\textbf{Theorem 6.3}.\textit{ Suppose we have the situation 
described above. Furthermore, suppose }$N_{j}(u)$ $\in$ $C^{0}({\Omega}_{j}$
$\cup$ $(\partial{\Omega}_{j}$ $\cap$ $(\Gamma \backslash A)))$ $\cap$ 
$C^{1}({\Omega}_{j})$\textit{ satisfies }$divN_{j}(u)$ $=$ $H$\textit{ in }
${\Omega}_{j}$ for any $j$.\textit{ Then }$u$\textit{ is a weak solution 
to (\ref{eqn1.4'}) in }$\Omega$\textit{ if and only if
for each }$p$ $\in$ $\Gamma \backslash A$\textit{, there exist }$j$, 
$l$\textit{ as described above, such that at }$p$,\textit{ there holds} 
  
\begin{equation*}
\mathit{(N}_{j}\mathit{(u)-N}_{l}\mathit{(u))\cdot \nu }_{j}\mathit{=(N}_{j}%
\mathit{(u)-N}_{l}\mathit{(u))\cdot \nu }_{l}\mathit{=0}\mathit{.}
\end{equation*}

\bigskip

We should remind the reader that for $u$ $\in$ $W^{1,1}(\Omega )$ and $H$ $\in$
$L^{\infty}(\Omega )$, $u$ is a weak solution 
to (\ref{eqn1.4'}) if and only if $u$ is a minimizer for (\ref{eqn1.3'}) in view of 
Theorem 3.3. 

\medskip

\textbf{Proof}. Let $U$ $\subset \subset $ $\Omega $ have compact
closure in $\Omega $, and suppose that the boundary $\partial U$ is $C^{1}$-smooth. 
For $\varphi $
$\in $ $C_{0}^{\infty}(\Omega )$ with support contained in $U,$ we compute

\begin{eqnarray}
\int_{U}N(u)\cdot \nabla \varphi +H\varphi  &=&\int_{U\backslash \Gamma
}N(u)\cdot \nabla \varphi +H\varphi   \label{eqn6.3'} \\
&=&\int_{\partial (U\backslash \Gamma )}\varphi N(u)\cdot \nu   \notag \\
&=&\sum_{(j,l)}\int_{\partial{\Omega}_{j}\cap
(\Gamma \backslash A)\cap U}\varphi (N_{j}(u)\cdot {\nu}_{j}+N_{l}(u)\cdot 
{\nu}_{l})  \notag \\
&=&\sum_{(j,l)}\int_{\partial{\Omega}_{j}\cap
(\Gamma \backslash A)\cap U}\varphi (N_{j}(u)-N_{l}(u))\cdot {\nu}_{j}.  \notag
\end{eqnarray}

\noindent For the last equality we have used $\nu_{l}$ $=$ $-\nu_{j}$. Now observe
that $u$ is a weak solution in $\Omega $ if and only if the first term of (%
\ref{eqn6.3'}) vanishes for any $\varphi $ $\in $ $C_{0}^{\infty }(\Omega )$
and associated $U.$ On the other hand, this is equivalent to concluding that 
$(N_{j}(u)$ $-$ $N_{l}(u))$ $\cdot $ $\nu_{j}$ $=$ $0$ on $\Gamma\backslash A$ by
(\ref{eqn6.3'}). 

\begin{flushright}
Q.E.D.
\end{flushright}

\bigskip

We remark that it is possible that $u$ $\in$ $C^{1}\backslash C^{2}$ while 
$N(u)$ $\in$ $C^{1}$ in the nonsingular domain. For instance, let $u$ $=$ $xy+g(y)$ 
with $g$ $\in$ $C^{1}\backslash C^{2}$. Take $\vec{F}$ $=$ $-\vec{X}^{\ast }$. We can
then compute $N(u)$ $=$ $(0,\pm 1)$ in the nonsingular domain defined by $2x+g'(y)$
$\neq$ $0$.
 
We will also make a remark on deducing the second equality in (\ref{eqn6.3'}).
First note that at points of $A$ with $H_{m-1}(\bar{A})$ $=$ $0$, $\nu$ may not exist. 
How do we deal with this?
For any $\varepsilon $ $>$ $0,$ we
can find a finite open cover $\cup _{j=1}^{k}D_{j}$ $\supset $ $\bar{A}$
such that $\sum_{j=1}^{k}H_{m-1}(\partial D_{j})$ $<$ $\varepsilon .$ By the divergence
theorem we have 

\begin{equation*}
\int_{(U\backslash \Gamma )\backslash \cup _{j=1}^{k}D_{j}}N(u)\cdot \nabla
\varphi +H\varphi =\int_{\partial \lbrack (U\backslash \Gamma )\backslash
\cup _{j=1}^{k}D_{j}]}\varphi N(u)\cdot \nu .
\end{equation*}

\noindent Passing to the limit as $\varepsilon $ $\rightarrow $ $0$ and observing that
the integrands are bounded (since $|N(u)|=1),$ we obtain 

\begin{equation*}
\int_{U\backslash \Gamma }N(u)\cdot \nabla \varphi +H\varphi =\int_{\partial
(U\backslash \Gamma )}\varphi N(u)\cdot \nu .
\end{equation*}

\noindent The idea of the above argument was used in \cite{CF74}. We have displayed this
idea in the proof of Lemma 6.1. We also
used a similar argument in the proof of Theorem 5.2 in \cite{CHMY04}. 
We remark that Pauls had a similar result (for $m=2$, $\vec{F}=-\vec{X}^{\ast }$,
and $H$ $=$ $0$) as Theorem C in \cite{Pau05}. Ritor\'{e} and Rosales also obtained 
a similar result for $C^{2}$-smooth minimizers (for $m=2$, $\vec{F}=-\vec{X}^{\ast }$,
and $H$ $=$ $constant$) as Theorem 4.15 in \cite{RR05}.    

\bigskip

\section{Examples}

We shall give examples of Lipschitz (continuous) minimizers in dimension
2. 

\bigskip

\textbf{Definition 7.1.} A $p$-area minimizer or a $p$-minimizer in short is a minimizer for (\ref{eqn1.2}) with $H$ $=$ $0$.

\bigskip

Throughout this section, we will always work on the situation that $m=2$, $\vec{F}=-\vec{X}^{\ast }$,
and $H$ $=$ $0$. Recall that the integral curves of $N^{\perp }(u)$ are
straight lines (see Section 4 in \cite{CHMY04}), called the characteristic
lines, segments, or rays. We call the angle between $\Gamma $ (oriented) and
a characteristic ray (with direction $N^{\perp }(u))$ in $B^{+}$ ($B^{-},$
respectively) touching a point $p$ $\in $ $\Gamma $ the incident (reflected,
respectively) angle at $p.$ Therefore geometrically (\ref{eqn6.1}) is
equivalent to saying that at $p$ $\in $ $\Gamma \cap B,$ either $N^{+}(u)=N^{-}(u)$ (see Figure 1(b)) or $N^{+}(u)\neq N^{-}(u)$ which implies

\begin{equation}
\textit{The incident angle}\mathit{=}\textit{The reflected angle. }
\label{eqn6.4}
\end{equation}

\noindent (see Figure 1(a))$.$ Suppose $u$ $\in $ $C^{2}$ at a point $p$ $%
\in $ $\Gamma \cap B$ and $\Gamma $ is a singular curve. Recall that if the characteristic line segments $%
\Gamma _{+}$ and $\Gamma _{-}$ in $B^{+}$ and $B^{-}$ respectively meet at $%
p,$ then $\Gamma _{+}\cup \{p\}\cup \Gamma _{-}$ must form a straight line
segment according to (the proof of) Proposition 3.5 in \cite{CHMY04}.
Therefore by (\ref{eqn6.4}) (note that $N^{+}(u)\ =\ -N^{-}(u)$ at $p$ in this situation), we can conclude that

\begin{equation}
\Gamma _{+}\textit{ and }\Gamma _{-}\textit{ are perpendicular to }\Gamma \textit{
at }p\textit{ if }u\in C^{2}\textit{ at }p.  \label{eqn6.5}
\end{equation}

The constraint (\ref{eqn6.5}) gives a necessary and sufficient condition
for a $C^{2}$-smooth solution of (\ref{eqn1.1}) with $H$ $=$ $0$ to be a $p$-minimizer. We can
have a function $u$ $\in $ $C^{2}(\Omega )$ which satisfies the p-minimal
surface equation $divN(u)$ $=$ $0$ in $\Omega \backslash S(u),$ but is not a
weak solution or a $p$-minimizer.

\bigskip

\textbf{Example 7.1}. Consider $\vec{F}$ $=$ $(-y,x)$ in the following $N(u)$%
's.

(a) By taking $a=\cos \vartheta ,$ $b=\sin \vartheta ,$ and $g(-bx+ay)$ $=$ (%
$\cot \vartheta )$ $(-bx+ay)^{2}$ in (1.2) of \cite{CHMY04} for $0$ $<$ $%
\vartheta $ $<$ $\frac{\pi }{2}$, we obtain $u(x,y)$ $=$ $-xy+y^{2}\cot
\vartheta .$ This is a $C^{2}$ smooth solution to $divN(u)$ $=$ $0$ in $%
R^{2}\backslash S(u)$ for $\vec{F}$ $=$ $(-y,x)$ by a direct computation. We
can easily determine the singular set $S(u)$ $\equiv $ $\{u_{x}-y$ $=$ $0,$ $%
u_{y}+x$ $=$ $0\}$ $=$ $\{y=0\}.$ On the other hand, $N^{\perp }(u)$ $=$ $%
(\cos \vartheta ,$ $\sin \vartheta )$ which is not perpendicular to the $x$%
-axis $\{y=0\}$ (see Figure 2(a))$.$ So in view of (\ref{eqn6.5}), this $u$
is not a $p$-minimizer on any
bounded domain $\Omega $ containing part of the $x$-axis.

(b)\textbf{\ }Let $u(x,y)$ $=$ $-xy+y^{2}\cot \vartheta $ for $y>0;$ $=$ $%
-xy+y^{2}\cot \eta $ for $y<0;$ $=0$ for $y=0$ where $0$ $<$ $\vartheta
,\eta $ $<$ $2\pi ,$ $\vartheta \neq \pi ,$ $\eta \neq \pi .$ We compute

\begin{eqnarray}
N^{\perp }(u) &=&(\frac{\cos \vartheta }{\sin \vartheta }|\sin \vartheta
|,|\sin \vartheta |)\text{ for }y>0;  \label{eqn6.6} \\
N^{\perp }(u) &=&(-\frac{\cos \eta }{\sin \eta }|\sin \eta |,-|\sin \eta |)%
\text{ for }y<0.  \notag
\end{eqnarray}

\noindent Observe that (\ref{eqn6.4}) (or (\ref{eqn6.1})) holds if and only
if $\vartheta +\eta $ $=$ $2\pi $ (see Figure 2(b)) by (\ref{eqn6.6})$.$
Therefore we conclude that $u$ is a ($C^{1,1}$-smooth) $p$-minimizer on any
bounded domain in $R^{2}$ if
and only if $\vartheta +\eta $ $=$ $2\pi $ in view of (\ref{eqn6.4}).

\setcounter{subfigure}{0}
\begin{figure}[ht]
\begin{center}
\mbox{
\subfigure[]{\includegraphics[width= 5.4cm]{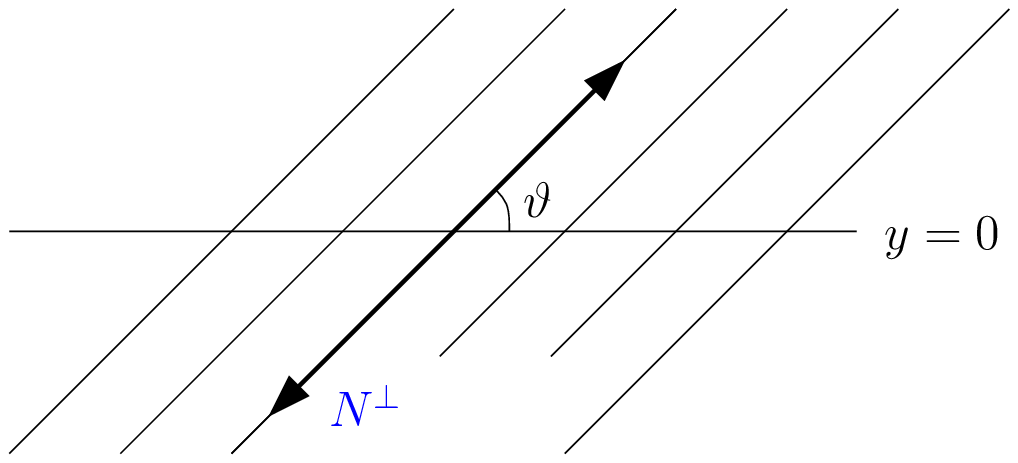}}
\subfigure[]{\includegraphics[width= 5.4cm]{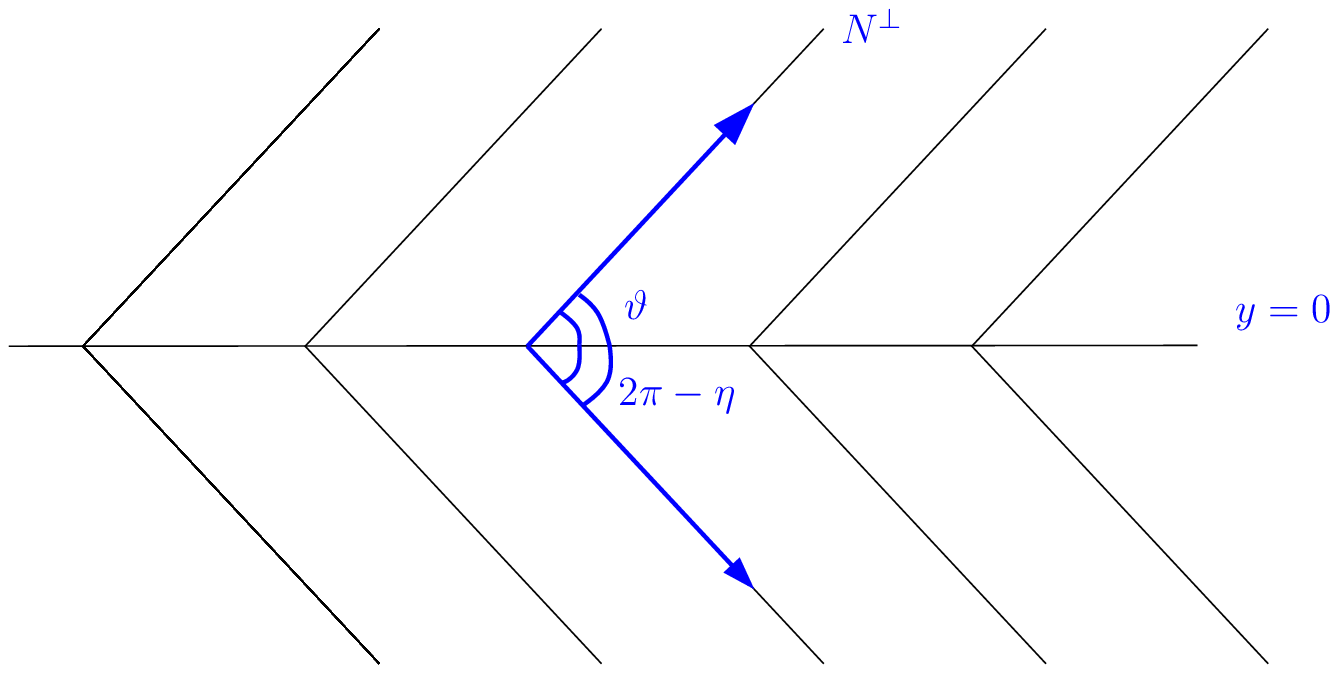}} }
\caption{}
 \end{center}
\end{figure}

\textbf{Example 7.2}. Let $u(x,y)$ $=$ $xy$ for $y>0$, and $u=0$ for $y\leq 0.$
Consider the case of $\vec{F}$ $=$ $(-y,x).$ Compute

\begin{eqnarray*}
N^{\perp }(u) &=&(1,0)\text{ for }x>0,y>0;\text{ }N^{\perp }(u)=(-1,0)\text{
for }x<0,y>0. \\
N^{\perp }(u) &=&\frac{(x,y)}{\sqrt{x^{2}+y^{2}}}\text{ for }y<0
\end{eqnarray*}

\noindent (see Figure 3). Observe that the positive $y$-axis $\{x=0,$ $y>0\}$
is a singular curve where (\ref{eqn6.5}) holds true. Also on the $x$-axis $%
\{y=0\}$ except the origin, $N^{\perp }(u)$ is continuous and hence (\ref%
{eqn6.1}) holds true (note that the $x$-axis is not a singular
curve, but is a curve where $u$ is not $C^{1}$ smooth). Applying Theorem 6.3 with $\Gamma$ $=$ $\{x=0,$ $y > 0\}$ $\cup$ $\{ y=0\}$, we conclude that $u$ is a (Lipschitz) $p$-minimizer on any bounded domain $\Omega$ $\subset$ $R^{2}$. 





\begin{figure}[ht]
\begin{center}
\psfrag{Si}{Singular curve}
\psfrag{Ch}{Characteristic lines}
 \includegraphics[width=7cm]{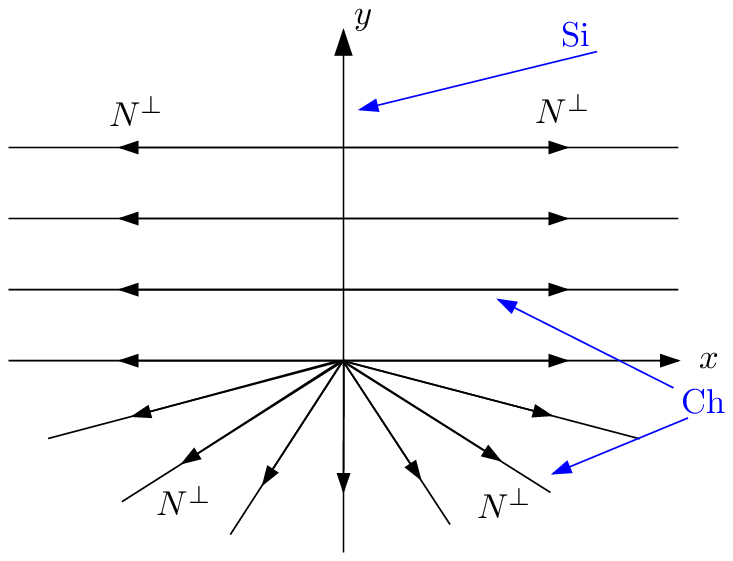}\\
  \end{center}
\caption{}
 \end{figure}

We remark that it is not possible to construct a Lipschitz $p$-minimizer having a loop consisting of characteristic lines (see
Figure 4 for an example). Indeed, by contradiction, suppose that the loop 
consists of three characteristic lines 
$\gamma _{1},$ $\gamma _{2},$ and $\gamma _{3}$ as indicated in Figure 4$.$
Let $\Delta $ denote the region surrounded by $\gamma _{1},$ $\gamma _{2},$
and $\gamma _{3}.$ We integrate the contact form $\Theta $ $\equiv $ $du$ $+$
$xdy$ $-$ $ydx$ over the loop as follows:

\begin{eqnarray*}
0 &=&\int_{\gamma _{1}\cup \gamma _{2}\cup \gamma _{3}}\Theta \text{ (}%
\gamma _{1},\gamma _{2},\text{ and }\gamma _{3}\text{ being Legendrian)} \\
&=&\int_{\Delta }d\Theta \text{ \ (Stokes' Theorem)} \\
&=&2\int_{\Delta }dx\wedge dy=2\text{ Area}(\Delta )\neq 0.
\end{eqnarray*}

\noindent This contradiction confirms our claim.

\begin{figure}[ht]
\begin{center}
 \includegraphics[width=5cm]{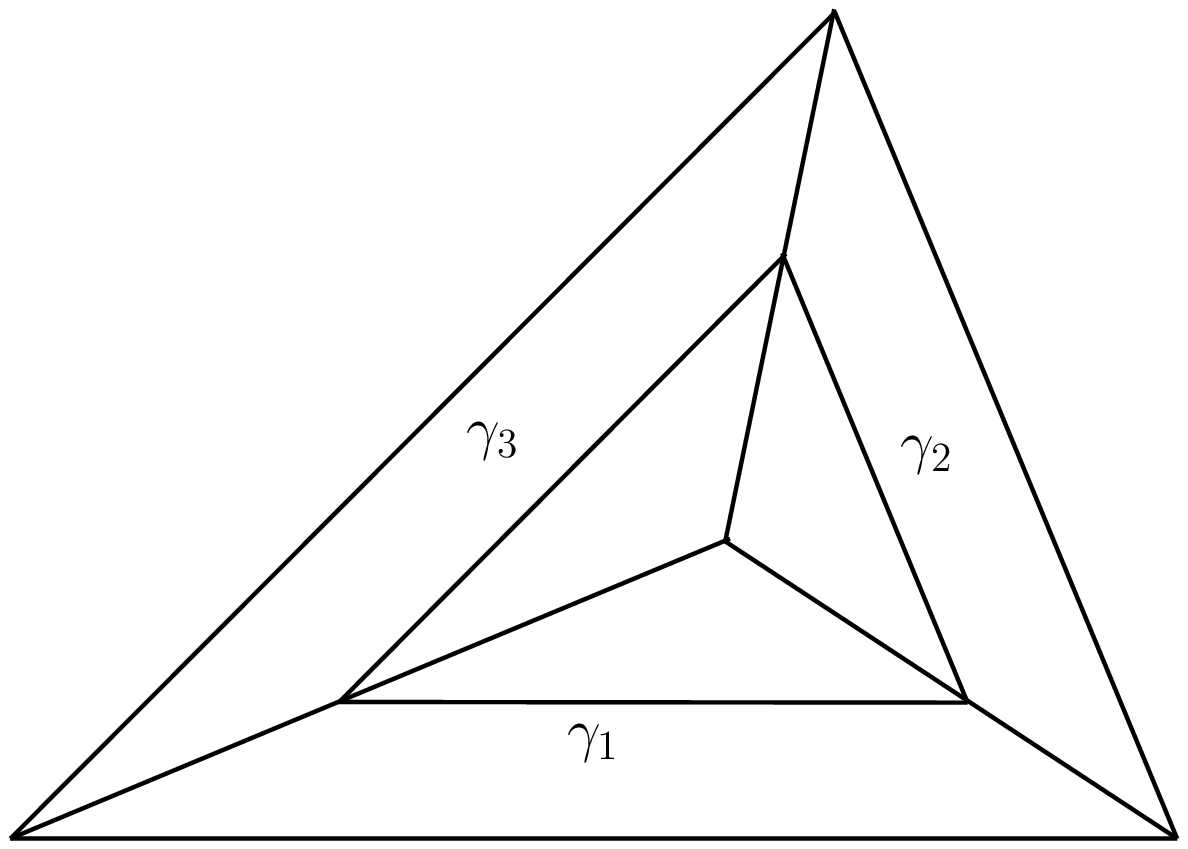}
\caption{}
  \end{center}
 \end{figure}

\bigskip

\textbf{Example 7.3}. There can be two distinct $C^{2}$-smooth p-minimal
graphs (i.e., satisfying (\ref{eqn1.1}) on nonsingular domain) having the
same boundary value and the same p-area, but both of them are not $p$-minimizers. 
Consider $u$
$=$ $x^{2}+xy,$ $v$ $=$ $xy+1-y^{2}$ (first given in \cite{Pau01}). We can
easily verify that $u$ and $v$ satisfy (\ref{eqn1.1}) with $H$ $=$ $0$ on their respective
nonsingular domains and have the same
value on the unit circle in the $xy$ - plane. But they do not satisfy (\ref{eqn6.5}). So by Proposition $6.2$ or Theorem 6.3, neither of them can be a $p$-minimizer. Compute the $p$-area (see (\ref{eqn1.2})) of $u$ and $%
v$ over the unit disc $\Delta $ as follows:

\begin{eqnarray*}
\mathcal{X(}u\mathcal{)} &=&\int_{\Delta }\sqrt{8}|x|dxdy=\frac{8\sqrt{2}}{3}%
, \\
\mathcal{X}(v) &=&\int_{\Delta }2|x-y|dxdy=\frac{8\sqrt{2}}{3}.
\end{eqnarray*}

\noindent So they have the same $p$-area. By the uniqueness of $p$-minimizers
(see Theorem B), we also conclude that neither $u$ nor $v$ can be the $p$-minimizer.

We are going to describe what the (unique) $p$-minimizer looks like on $%
\Delta $ with the boundary value (or curve) $\rho (\theta )$ $\equiv $ $\cos
^{2}\theta $ $+$ $\cos \theta $ $\sin \theta $ ($\theta $ is the standard
angle parameter for $\partial \Delta ).$ Let $(\alpha ,$ $\beta ,$ $\gamma )$
be a point of a line segment $\tilde{L}$ $\subset $ $\bar{\Delta}\times R$
meeting the boundary curve with the projection $L$ $\subset $ $\bar{\Delta}$
passing through the origin$.$ Suppose $\theta ^{\prime }$ is the angle
between the positive $x$-axis and part of $L$, lying in the upper half plane.
Then we have

\begin{eqnarray}
\alpha &=&t\cos \theta ^{\prime },\text{ }\beta =t\sin \theta ^{\prime }
\label{eqn6.19} \\
\gamma &=&\cos ^{2}\theta ^{\prime }+\cos \theta ^{\prime }\sin \theta
^{\prime }  \notag
\end{eqnarray}

\noindent for $-1$ $\leq $ $t$ $\leq $ $1$ (note that $\rho (\pi +\theta )$ $%
=$ $\rho (\theta )).$ Suppose the contact plane passing through $(\alpha ,$ $%
\beta ,$ $\gamma )$ intersects the boundary curve $\rho $ at $(\cos \theta ,$
$\sin \theta ,$ $\rho (\theta ))$. Then we have the following relation:

\begin{equation}
\rho (\theta )-\rho (\theta ^{\prime })+t\sin (\theta -\theta ^{\prime })=0
\label{eqn6.20}
\end{equation}

\noindent by observing that $z-\gamma +x(y-\beta )-y(x-\alpha )=0$ is the
equation for such a contact plane in $R^{3}$ with coordinates $(x,y,z).$ By
elementary trigonometry for the above specific $\rho ,$ we can reduce (\ref%
{eqn6.20}) to

\begin{equation}
\frac{\sqrt{2}}{2}[\sin (2\theta +\frac{\pi }{4})-\sin (2\theta ^{\prime }+%
\frac{\pi }{4})]+t\sin (\theta -\theta ^{\prime })=0.  \label{eqn6.21}
\end{equation}

The idea is to choose $\theta ^{\prime }$ such that $\sin (2\theta
^{\prime }+\frac{\pi }{4})$ $=$ $0.$ Then we solve (\ref{eqn6.21}) for $%
\theta $ (perhaps we have multiple solutions)$.$ Keeping $\tilde{L}$ or $L$ associated to $%
\theta ^{\prime }$ as the singular set in mind, we connect $(\alpha ,$ $\beta ,$
$\gamma )$ $\in $ $\tilde{L}$ to a point of the boundary curve, associated
to $\theta ,$ by a line segment. Since these line segments are Legendrian,
their union forms a Legendrian ruled surface, hence a $p$-minimal surface (%
\cite{CHMY04}). Moreover, if two characteristic lines (i.e., above Legendrian
lines projected to the $xy$-plane) meet at a point of $\tilde{L},$
condition (\ref{eqn6.4}) holds. So in this way we can construct the $p$-minimizer by
Proposition $6.2$ or Theorem 6.3. We give more details below.

First solving $\sin (2\theta ^{\prime }+\frac{\pi }{4})$ $=$ $0$ gives $%
\theta ^{\prime }$ $=$ $(n-\frac{1}{4})\frac{\pi }{2}$ where $n$ is an
integer. There are two such $\theta ^{\prime }$'s modulo an integral
multiple of $\pi ,$ namely $\theta ^{\prime }$ $=$ $\frac{3}{8}\pi $ and $%
\theta ^{\prime }$ $=$ $\frac{7}{8}\pi .$ We take $\theta ^{\prime }$ $=$ $%
\frac{3}{8}\pi $ (it turns out that $\theta ^{\prime }$ $=$ $\frac{7}{8}\pi $
won't give rise to a $p$-minimal graph in the following argument). So (\ref%
{eqn6.21}) is reduced to 
\begin{equation}
\frac{\sqrt{2}}{2}\sin 2(\theta -\frac{3}{8}\pi )=t\sin (\theta -\frac{3}{8}%
\pi ).  \label{eqn6.22}
\end{equation}

\noindent (note that for $\theta ^{\prime }$ $=$ $\frac{7}{8}\pi $ we have $%
-t$ instead of $t$ in (\ref{eqn6.22})). By the double angle formula, we
deduce from (\ref{eqn6.22}) that 
\begin{equation}
(a)\text{ }\cos (\theta -\frac{3}{8}\pi )=\frac{t}{\sqrt{2}}\text{ ;}(b)%
\text{ }\sin (\theta -\frac{3}{8}\pi )=0.  \label{eqn6.23}
\end{equation}

\noindent The solutions to $(b)$ of (\ref{eqn6.23}) are $\frac{3}{8}\pi $ $%
+n\pi $ for any integer $n,$ which we ignore. We have two solutions $\theta
_{1},$ $\theta _{2}$ (modulo an integral multiple of $2\pi )$ to $(a)$ of (%
\ref{eqn6.23}) for a given $t$ with the relation

\begin{equation}
\theta _{1}-\frac{3}{8}\pi =\frac{3}{8}\pi -\theta _{2}.  \label{eqn6.24}
\end{equation}

\noindent When $t$ runs from $-1$ to $1,$ $\theta _{1}$ runs from $\frac{9}{8%
}\pi $ to $\frac{5}{8}\pi $ clockwise while $\theta _{2}$ runs from $-\frac{3%
}{8}\pi $ to $\frac{1}{8}\pi $ counterclockwise (See Figure 5 below).

\begin{figure}[ht]
\begin{center}
 \includegraphics[width=6cm]{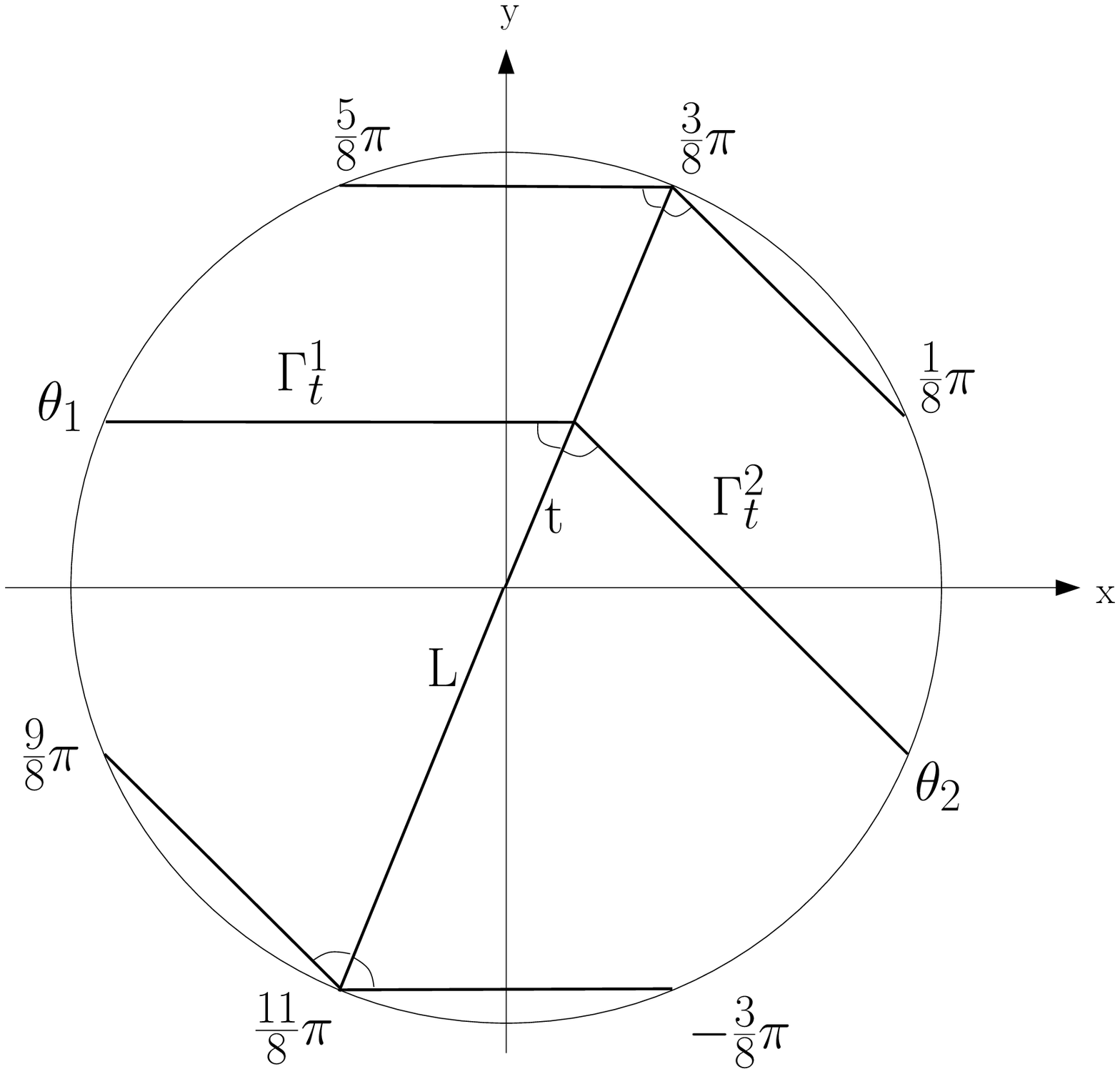}
\caption{}
 \end{center}
 \end{figure}

Denote the line segments between $(\alpha ,$ $\beta )$ $\in $ $L$ ($\theta
^{\prime }$ $=$ $\frac{3}{8}\pi $) and the boundary point $(\cos \theta _{j},
$ $\sin \theta _{j}),$ $j$ $=$ $1,$ $2,$ by $\Gamma _{t}^{1},$ $\Gamma
_{t}^{2},$ respectively. $\Gamma _{t}^{1}$ and $\Gamma _{t}^{2}$ are the $xy$%
-plane projections of two Legendrian lines $\tilde{\Gamma}_{t}^{1},$ $\tilde{%
\Gamma}_{t}^{2}$ connecting $(\alpha ,$ $\beta ,$ $\gamma )$ $\in $ $\tilde{L%
}$ to $(\cos \theta _{j},$ $\sin \theta _{j},$ $\rho (\theta _{j})),$ $j$ $=$
$1,$ $2,$ respectively. We will define a graph $\check{u}$ over $\bar{\Delta}%
,$ whose restriction to the region $\Omega $ $\equiv $ $L$ $\cup $ $(\cup
_{j=1,2;-1\leq t\leq 1}\Gamma _{t}^{j})$ is $\tilde{L}$ $\cup $ $(\cup
_{j=1,2;-1\leq t\leq 1}\tilde{\Gamma}_{t}^{j}).$ We can parametrize $\tilde{%
\Gamma}_{t}^{2},$ say, in the following form (see (4.9) in\cite{CHMY04}):

\begin{eqnarray*}
x &=&s(\sin \eta (t))+\alpha (t) \\
y &=&-s(\cos \eta (t))+\beta (t) \\
z &=&s[\beta (t)\sin \eta (t)+\alpha (t)\cos \eta (t)]+\gamma (t).
\end{eqnarray*}

\noindent Here $\eta (t)$ $=$ $\frac{\pi }{2}$ $+$ $\theta _{2}(t)$ $-$ $%
\delta (t)$ in which $\cos (\theta _{2}(t)-\frac{3}{8}\pi )=\frac{t}{\sqrt{2}%
}$ (see $(a)$ of (\ref{eqn6.23})) and $\tan \delta (t)$ $=$ $t\sqrt{1-t^{2}/2%
}/(1-t^{2}/\sqrt{2})$ by elementary plane geometry (we leave the details to
the reader). On the other hand, (\ref{eqn6.4}) holds along $L$ due to (\ref%
{eqn6.24}). Therefore $\check{u}$ $\in $ $C^{1,1}$ is a weak solution to $%
(1.1)$ with $H$ $=$ $0$ over the region $\Omega .$ The remaining domain $%
\bar{\Delta}\backslash \Omega $ consists of four small fan-shaped regions
(see Figure 5). For each of such regions, we can connect two points on the
boundary curve, indicated by $\theta ^{\prime }$ and $\theta $ which are
related by (\ref{eqn6.21}) with $t$ $=$ $1.$ Thus we obtain a family of
Legendrian line segments whose lengths are getting smaller when both $\theta
^{\prime }$ and $\theta $ tend to some critical value (e.g., for the
fan-shaped region between $\frac{1}{8}\pi $ and $\frac{3}{8}\pi ,$ the
critical value is $\frac{1}{4}\pi $ by solving $\frac{d\theta }{d\theta
^{\prime }}$ $=$ $0).$ These Legendrian line segments form a portion of the
graph $\check{u}$ over $\bar{\Delta}\backslash \Omega .$ So $\check{u}$ is a 
$C^{2}$-smooth $p$-minimal graph over $\bar{\Delta}\backslash \Omega .$
Altogether $\check{u}$ $\in $ $C^{1,1}(\bar{\Delta})$ is the (unique)
$p$-minimizer by Theorem 3.3.

\bigskip

\section{Appendix: uniqueness of solutions to (\ref{eqn4.1})}

The existence of solutions to (\ref{eqn4.1}) is asserted in Theorem 4.5. In
this section we are going to prove the uniqueness. In fact we can obtain
more general results. First we define

\begin{equation}
N_{\varepsilon }(u)\equiv \frac{\nabla u+\vec{F}}{\sqrt{\varepsilon
^{2}+|\nabla u+\vec{F}|^{2}}}.  \label{eqn8.1}
\end{equation}

\noindent Let $\vec{\alpha}$ $\equiv $ $(\varepsilon ,\nabla u+\vec{F})$ $%
\in $ $R\times R^{m}$ $=$ $R^{m+1}.$ Denote $|\vec{\alpha}|$ as $\alpha .$
Similarly let $\vec{\beta}$ $\equiv $ $(\varepsilon ,\nabla v+\vec{F})$ and $%
\beta \equiv |\vec{\beta}|$.

\bigskip

\textbf{Lemma 8.1. }\textit{Let }$u,$\textit{\ }$v$\textit{\ }$\in $\textit{%
\ }$W^{1}(\Omega )$\textit{\ where }$\Omega $ $\subset $ $R^{m}$\textit{\ (}$%
m\geq 1)$ \textit{is an arbitrary domain. Then}

\begin{equation}
\mathit{(N}_{\varepsilon }\mathit{(u)-N}_{\varepsilon }\mathit{(v))\cdot
(\nabla u-\nabla v)\geq }\frac{\alpha +\beta }{2}\mathit{\mid N}%
_{\varepsilon }\mathit{(u)-N}_{\varepsilon }\mathit{(v)\mid }^{2}\mathit{.}
\label{eqn8.2}
\end{equation}

\noindent \textit{Moreover, the equality holds for }$\varepsilon =0.$\textit{%
\ When }$\varepsilon >0,$\textit{\ }$(N_{\varepsilon }(u)-N_{\varepsilon
}(v))$\textit{\ }$\cdot $\textit{\ }$(\nabla u-\nabla v)$\textit{\ }$=$%
\textit{\ }$0$\textit{\ if and only if }$\nabla u$\textit{\ }$=$\textit{\ }$%
\nabla v.$

\bigskip

\textbf{Proof. }We compute

\begin{eqnarray}
&&(N_{\varepsilon }(u)-N_{\varepsilon }(v))\cdot (\nabla u-\nabla v)
\label{eqn8.3} \\
&=&(\frac{\nabla u+\vec{F}}{\alpha }-\frac{\nabla v+\vec{F}}{\beta })\cdot
\{(\nabla u+\vec{F})-(\nabla v+\vec{F})\}  \notag \\
&=&\{\frac{(\varepsilon ,\nabla u+\vec{F})}{\alpha }-\frac{(\varepsilon
,\nabla v+\vec{F})}{\beta }\}\cdot \{(\varepsilon ,\nabla u+\vec{F}%
)-(\varepsilon ,\nabla v+\vec{F})\}  \notag \\
&=&\{\frac{\vec{\alpha}}{\alpha }-\frac{\vec{\beta}}{\beta }\}\cdot \{\vec{%
\alpha}-\vec{\beta}\}=(\alpha +\beta )(1-\cos \theta )  \notag
\end{eqnarray}

\noindent where $\vec{\alpha}\cdot \vec{\beta}=\alpha \beta \cos \theta .$
On the other hand, we can estimate

\begin{eqnarray}
& & \mid N_{\varepsilon }(u)-N_{\varepsilon }(v)\mid ^{2}  \label{eqn8.4} \\
&\leq &\mid \frac{\vec{\alpha}}{\alpha }-\frac{\vec{\beta}}{\beta }\mid
^{2}=2(1-\cos \theta ).  \notag
\end{eqnarray}

\noindent Now (\ref{eqn8.2}) follows from (\ref{eqn8.3}) and (\ref{eqn8.4}).
Observing that the equality in (\ref{eqn8.4}) holds for $\varepsilon $ $=$ $%
0,$ we obtain the equality in (\ref{eqn8.2}) for $\varepsilon $ $=$ $0.$
Suppose $(N_{\varepsilon }(u)$ $-$ $N_{\varepsilon }(v))$ $\cdot $ $(\nabla
u $ $-$ $\nabla v)$ $=$ $0.$ By (\ref{eqn8.2}) we have $N_{\varepsilon }(u)$ 
$= $ $N_{\varepsilon }(v).$ Taking the modulus of this equality gives 
\TEXTsymbol{\vert}$\nabla u+\vec{F}|$ $=$ $|\nabla v+\vec{F}|$ if $%
\varepsilon $ $>$ $0.$ It follows that $\nabla u$ $=$ $\nabla v.$

\begin{flushright}
Q.E.D.
\end{flushright}

\bigskip

We remark that the equality in (\ref{eqn8.2}) for $\varepsilon $ $=$ $0$ has
been obtained as Lemma $5.1^{\prime }$ in \cite{CHMY04}. Recall $%
Q_{\varepsilon }u$ $\equiv $ $divN_{\varepsilon }(u)$ (see (\ref{eqn4.1}), (%
\ref{eqn8.1})). Note that for $\vec{F}$ $=$ $0,$ $\varepsilon $ $=$ $1,$ $Q_{\varepsilon }u$
is the Riemannian mean curvature of the graph defined by $u.$ In this case,
the above inequality has been obtained in \cite{Mik79}, \cite{Hw88}%
, and \cite{CK91} independently.

\bigskip

\textbf{Definition 8.1. }Let $\Omega \subset R^{m}$\ be a bounded domain and
$\varepsilon$ $>$ $0$.
Suppose $u,$ $v$ $\in $ $W^{1}(\Omega )$ and $\vec{F}$ is measurable$.$ We
say $Q_{\varepsilon }u$ $-$ $Q_{\varepsilon }v$ $\geq $ $0$ ($\leq $ $0,$
respectively) weakly if for any $\varphi $ $\in $ $C_{0}^{1}(\Omega ),$ $%
\varphi $ $\geq $ $0,$ there holds

\begin{equation}
\int_{\Omega }(N_{\varepsilon }(u)-N_{\varepsilon }(v))\cdot \nabla \varphi
\leq 0\text{ (}\geq 0,\text{ respectively}).  \label{eqn8.5}
\end{equation}

Note that $N_{\varepsilon }(u)$ and $N_{\varepsilon }(v)$ are integrable
since they are bounded by $1.$ We have the following comparison principle
for $Q_{\varepsilon }.$

\bigskip

\textbf{Theorem 8.2. }\textit{Let }$\Omega \subset R^{m}$\textit{\ be a
bounded domain and }$\varepsilon $ $>$ $0$\textit{. Suppose }$u,$\textit{\ }$%
v$\textit{\ }$\in $\textit{\ }$C^{1}(\Omega )$\textit{\ }$\cap $\textit{\ }$%
C^{0}(\bar{\Omega})$\textit{\ satisfy }$Q_{\varepsilon }u$\textit{\ }$-$%
\textit{\ }$Q_{\varepsilon }v$\textit{\ }$\geq $\textit{\ }$0$\textit{\ (}$%
\leq $\textit{\ }$0,$\textit{\ respectively) weakly and }$u$\textit{\ }$-$%
\textit{\ }$v$\textit{\ }$\leq $\textit{\ }$0$\textit{\ (}$\geq 0,$\textit{\
respectively}$)$\textit{\ on }$\partial \Omega .$\textit{\ Then }$u$\textit{%
\ }$-$\textit{\ }$v$\textit{\ }$\leq $\textit{\ }$0$\textit{\ (}$\geq 0,$%
\textit{\ respectively}$)$\textit{\ in }$\Omega .$

\bigskip

\textbf{Proof. }Given $a$ $>$ $0,$ we choose a function $f_{a}$ $\in $ $%
C^{1}(R)$ with the property that $f_{a}$ $\equiv $ $0$ in $(-\infty ,$ $a],$ 
$f_{a}$ $>$ $0$, and $f_{a}^{\prime }$ $>$ $0$ in $(a,$ $\infty ).$ Observe
that $f_{a}(u-v)$ $\in $ $C_{0}^{1}(\Omega )$ (i.e., $f_{a}(u-v)$ $\in $ $%
C^{1}$ and has compact support in $\Omega $) by the assumption $u$\textit{\ }%
$-$\textit{\ }$v$\textit{\ }$\leq $\textit{\ }$0$ on $\partial \Omega $. It
follows from (\ref{eqn8.5}) that

\begin{eqnarray}
0 &\geq &\int_{\Omega }(N_{\varepsilon }(u)-N_{\varepsilon }(v))\cdot \nabla
(f_{a}(u-v))  \label{eqn8.6} \\
&=&\int_{\{u-v>a\}}(N_{\varepsilon }(u)-N_{\varepsilon }(v))\cdot
f_{a}^{\prime }(u-v)(\nabla u-\nabla v)  \notag \\
&\geq &0\text{ \ (by (\ref{eqn8.2})).}  \notag
\end{eqnarray}

\noindent Therefore we have $(\nabla u$ $-$ $\nabla v)$ $\cdot $ $%
(N_{\varepsilon }(u)$ $-$ $N_{\varepsilon }(v))$ $=$ $0$ in $\{u$ $-$ $v$ $>$
$a\}$ since $f_{a}^{\prime }(u-v)$ $>$ $0$ and $(N_{\varepsilon }(u)$ $-$ $%
N_{\varepsilon }(v))$ $\cdot $ $(\nabla u-\nabla v)\geq $ $0$ in (\ref%
{eqn8.6}). It follows that $\nabla u$ $=$ $\nabla v$ in $\{u$ $-$ $v$ $>$ $%
a\}$ by Lemma 8.1. Thus we obtain $u$ $-$ $v$ $\equiv $ $a$ in $\{u$ $-$ $v$ 
$>$ $a\},$ a contradiction. So $\{u$ $-$ $v$ $>$ $a\}$ is empty. Since $a$ $%
> $ $0$ is arbitrary, we conclude that $\{u$ $-$ $v$ $>$ $0\}$ is empty. So $%
u$\textit{\ }$-$\textit{\ }$v$\textit{\ }$\leq $\textit{\ }$0$ in $\Omega .$

\begin{flushright}
Q.E.D.
\end{flushright}

We remark that basically the above result can be deduced from Theorem 10.7 in \cite{GT83}.

\bigskip

\textbf{Corollary 8.3. }\textit{Let }$\Omega \subset R^{m}$\textit{\ be a
bounded domain. Let }$\varepsilon $\textit{\ }$>$\textit{\ }$0.$
\textit{ Suppose }$u,$\textit{\ }$v$\textit{\ }$\in $\textit{\ }$C^{2}(\Omega )$%
\textit{\ }$\cap $\textit{\ }$C^{0}(\bar{\Omega})$\textit{ and }$\vec{F}\in C^{1}(\Omega )$\textit{\ satisfy }$%
Q_{\varepsilon }u$\textit{\ }$=$\textit{\ }$Q_{\varepsilon }v$\textit{\ in }$%
\Omega $ \textit{and }$u$\textit{\ }$=$\textit{\ }$v$\textit{\ on }$\partial
\Omega .$\textit{\ Then }$u\equiv v$\textit{\ in }$\Omega .$

\textit{\bigskip }

\textbf{Theorem 8.4. }\textit{Let }$\Omega \subset R^{m}$\textit{\ be a
bounded domain and }$\varepsilon $\textit{\ }$>$\textit{\ }$0$\textit{.
Suppose }$u,$\textit{\ }$v$\textit{\ }$\in $\textit{\ }$W^{1,1}(\Omega )$%
\textit{\ satisfy }$Q_{\varepsilon }u$\textit{\ }$-$\textit{\ }$%
Q_{\varepsilon }v$\textit{\ }$\geq $\textit{\ }$0$\textit{\ (}$\leq $\textit{%
\ }$0,$\textit{\ respectively) weakly and }$(u$\textit{\ }$-$\textit{\ }$%
v)^{+}$\textit{\ (}$(u$\textit{\ }$-$\textit{\ }$v)^{-},$\textit{\
respectively) }$\in $\textit{\ }$W_{0}^{1,1}(\Omega ).$\textit{\ Then }$u$%
\textit{\ }$-$\textit{\ }$v$\textit{\ }$\leq $\textit{\ }$0$\textit{\ (}$%
\geq $\textit{\ }$0,$\textit{\ respectively) in }$\Omega .$

\bigskip

\textbf{Proof. }First\textbf{\ }we observe that (\ref{eqn8.5}) still holds
for $\varphi $ $\in $ $W_{0}^{1,1}(\Omega ),$ $\varphi $ $\geq $ $0.$ It
follows that

\begin{eqnarray}
0 &\geq &\int_{\Omega }(N_{\varepsilon }(u)-N_{\varepsilon }(v))\cdot \nabla
(u\mathit{\ }-\mathit{\ }v)^{+}  \label{eqn8.7} \\
&=&\int_{\{u-v>0\}}(N_{\varepsilon }(u)-N_{\varepsilon }(v))\cdot \nabla (u%
\mathit{\ }-\mathit{\ }v)  \notag \\
&\geq &0  \notag
\end{eqnarray}

\noindent by (\ref{eqn8.2}). Therefore if $u$ $-$ $v$ $>$ $0,$ $\nabla (u\mathit{\ }-\mathit{\ }%
v)^{+}$ $=$ $\nabla (u\mathit{\ }-\mathit{\ }v)$ $=$ $0$ by Lemma 7.6 in %
\cite{GT83}, (\ref{eqn8.7}), and Lemma 8.1. Also if $u$ $-$ $v$ $\leq $ $0,$ 
$\nabla (u\mathit{\ }-\mathit{\ }v)^{+}$ $=$ $0$ by Lemma 7.6 in \cite{GT83}%
. Altogether we have shown that $\nabla (u\mathit{\ }-\mathit{\ }v)^{+}$ $=$ 
$0$ in $\Omega .$ Now applying the Sobolev inequality to $(u\mathit{\ }-%
\mathit{\ }v)^{+}$ $\in $ $W_{0}^{1,1}(\Omega ),$ we obtain $(u\mathit{\ }-%
\mathit{\ }v)^{+}$ $=$ $0$ in $\Omega .$ That is, $u$\textit{\ }$-$\textit{\ 
}$v$\textit{\ }$\leq $\textit{\ }$0$ in $\Omega .$

\begin{flushright}
Q.E.D.
\end{flushright}

\bigskip

We remark that the proof of Theorem 8.4 is based on the idea of the proof of
Theorem 8.1 in \cite{GT83}.

\end{document}